# Floer homologies for Lagrangian intersections and instantons

Ronnie Lee and Weiping Li

Dedicated to the memory of Andreas Floer

# Contents

# 1  Introduction

In 1985 lectures at MSRI, A. Casson introduced an interesting integer valued invariant for any oriented integral homology 3-sphere $Y$ via beautiful constructions on representation spaces (see [1] for an exposition). The Casson invariant $\lambda(Y)$ is roughly defined by measuring the oriented number of irreducible representations of the fundamental group $\pi_1(Y)$ in $SU(2)$. Such an invariant generalized the Rohlin invariant and gives surprising corollaries in low dimensional topology.

More precisely, consider a Heegaard splitting of $Y$ into two handle bodies $Y_1, Y_2$ by a compact Riemann surface $\Sigma$. Let $\mathcal{M}(\Sigma)$ denote the moduli space of $SU(2)$-representations of $\pi_1(\Sigma)$. As is well known, $\mathcal{M}(\Sigma)$ has the structure of a complex variety and in particular a stratified symplectic space (c.f. [4], [21]). The representations of $Y_1, Y_2$ give rise to Lagrangian subspaces $L_1, L_2$ in $\mathcal{M}(\Sigma)$. Under the assumption that $Y$ is an integral homology sphere and away from singularities of $\mathcal{M}(\Sigma)$, the intersection $L_1 \cap L_2$ of these Lagrangian subspaces $L_1, L_2$ is compact. Therefore we can perturb them into general position and count the intersection number $\#(L_1 \cap L_2)$,

$$2\lambda(Y) = \#L_1 \cap L_2. \tag{1.1}$$

In [34], C. Taubes showed that the Casson invariant can be interpreted as an infinite dimensional generalization of the classical Euler characteristic number in gauge theory of 3-manifolds. Indeed, it is one half of the Euler characteristic of Floer instanton homology theory $HF_*(Y)$ (see §2.1 and [13]). The Casson invariant can be regarded as the Lagrangian intersections of $L_1, L_2$ in the symplectic space $\mathcal{M}(\Sigma)$. In view of Taubes' work, Atiyah in [3] posed the problem whether there is a way of computing $HF_*(Y)$ via the symplectic approach.

Since $\mathcal{M}(\Sigma)$ is a Kähler manifold with **singularties**, in particular symplectic stratified space. Atiyah and Bott in [4] have shown that the symplectic structure is canonical and independent of the metrics on $\Sigma$. Furthermore $L_1$ and $L_2$ are stratified Lagrangian subspaces, i.e. at each corresponding stratum the subspaces of middle dimension on which the symplectic 2-form $\omega$ on the respective stratum of $\mathcal{M}(\Sigma)$ restricted is identically zero. Floer [15] studied the problem of intersections of Lagrangian submanifolds of compact symplectic manifold, and for the purpose of proving Aronld's conjecture, developed symplectic Floer homology theories ([15], [16]). Atiyah in [3] then conjectured that the Floer homology defined in the symplectic context from Lagrangian intersections ([15]) coincides with the Floer homology defined in the anti-self-dual context ([13]). In addition, he outlined an idea for proving this conjecture and in [18] Floer also listed this as one of his problems.

S. Dostoglou and D. Salamon [11] took the first step toward Atiyah conjecture in the special case when the underlying 3-manifold is a mapping torus $\Sigma_h$ for the diffeomorphism $h : \Sigma \to \Sigma$ induced by an automorphism $f : P \to P$, $P$ is a nontrivial $SO(3)$ bundle over a Riemann surface $\Sigma$, as suggested by Floer. In this situation, $P_f \to \Sigma_h$ is the $SO(3)$ bundle with $w_2(P_f) \neq 0$. The moduli space of flat connections over this bundle is a compact **smooth** symplectic manifold with dimension $6g - 6$, where $g \geq 2$ is the genus of $\Sigma$. The flat connections over $P_f$ correspond naturally to the fixed points of the symplectomorphism $\phi_f : \mathcal{M}(\Sigma) \to \mathcal{M}(\Sigma)$ induced by $f$ on the moduli space. Using $\phi_f$ there are two well-defined Floer homology theories $HF_*^{\text{sym}}(\mathcal{M}(\Sigma), \phi_f)$ (see [16]) and $HF_*^{\text{ins}}(\Sigma_h, P_f)$ (see [18]). In [11], Dostoglou and Salamon proved that there is an natural isomorphism

$$HF_*^{\text{sym}}(\mathcal{M}(\Sigma), \phi_f) = HF_*^{\text{ins}}(\Sigma_h, P_f). \tag{1.2}$$

For homology 3-sphere case, one needs to take care of the singularity of the symplectic space $\mathcal{M}(\Sigma)$. In this paper, we combine a technique of Oh [29] and minimal surface theory to get a well-defined Floer symplectic homology $HF_*(\mathcal{M}(\Sigma); L_1, L_2)$ $(= HF_*^{\text{sym}}(\mathcal{R}(Y_1), \mathcal{R}(Y_2); \mathcal{R}(Y_0)))$ its definition in §2.2.3). In [29], Oh suggested a possible definition for the representation variety $\mathcal{R}(Y_0)$ and also indicated some of the technical difficulties along the singular strata. With this definition of $HF_*^{\text{sym}}(\mathcal{R}(Y_1), \mathcal{R}(Y_2); \mathcal{R}(Y_0))$ at hand, Atiyah's conjecture can be stated as whether it agrees with the instanton theory. The purpose of this paper is to offer an affirmative answer to this conjecture.

**Theorem:** *For a Heegaard decomposition $(Y; Y_1, Y_2, , Y_0)$ of a homology 3-sphere $Y$, there is a well-defined Floer symplectic homology $HF_*^{sym}$ for the Lagrangian intersections of representation spaces with genus $g \geq 3$. There is also a natural isomorphism*

$$HF_*(Y) \cong HF_*^{sym}(\mathcal{R}(Y_1), \mathcal{R}(Y_2); \mathcal{R}(Y_0)). \tag{1.3}$$

The paper is organized as follows. In §2, we review the Floer instanton homology theory briefly, and discuss how to deal with the possible singularities in symplectic setting, and the Oh's extension for Floer symplectic homology theory. §3 is a discussion of perturbations for both theories to identify the generators of the chain groups, also an application of the work by S. Cappell, R. Lee and E. Miller [6] and Yoshida [36] to identify the spectral flow with Maslov index. The uniformly lowest eigenvalues for self-duality operator and Cauchy-Riemann operator are studied in §4. We give some estimates essentially due to [11] in §5. In §6, the main effort of this paper shows that the nonlinear comparisions between instanton and pseudoholomorphic curve (we will omit the adjective "pseudo" henceforth) can be deformed into each other for one dimensional moduli space. The difficulty is to deform an ASD connection into a holomorphic curve (since the idea in [11] to show the surjectivity of the deformation from holmorphic curve to ASD connection does not work in this case). We further exploit the gauge theory of Heegaard decomposition from [34]. Using Kuranishi technique and deformation of metrics on handlebodies, we can first deform an ASD connection into the right Sobolev space of symplectic setting and then deform again to a holomorphic curve via inverse function theorem. Hence following the same orientation which identified in [11], we identify the two Floer boundary maps.

**Acknowledgements:** Both authors would like to thank Y. Oh and D. Salamon for many helpful discussions and useful comments. R.L acknowledges partially support from NSF. W.L acknowledges MSRI and I.Newton institute for their hospitality and support.

## 2 Floer homologies

### 2.1 Floer instanton homology for homology 3-spheres

In this subsection, we will give a brief description of the definition of Floer instanton homology. For details see [9], and [13].

Let $Y$ be a homology 3-sphere, i.e. an oriented closed 3-dimensional smooth manifold with $H_1(Y, Z) = 0$, and let $P \to Y$ be a smooth principal $SU(2)$-bundle (this bundle is trivial). Fix a trivialization $Y \times SU(2)$ of $P$ and let $\theta$ be the associated trivial connection. Denote the Sobolev $L_k^p$-space of connections on $P$ by $\mathcal{A}(P)$. This space has a natural affine structure with underlying vector space $\Omega^1(Y, adP)$, where $adP$ is the adjoint bundle. $\mathcal{A}(P)$ is acted upon by the gauge group $\mathcal{G}$ of $L_{k+1}^p$-automorphisms of $P$, and the orbit space $\mathcal{B}(P) = \mathcal{A}(P)/\mathcal{G}$ is well-defined when $k + 1 > \frac{3}{p}$. The irreducible connections form an open dense subspace $\mathcal{B}^*(P)$ of $\mathcal{B}(P)$ which is a Banach manifold with

$$T_a \mathcal{B}^*(P) \equiv \{\alpha \in L_k^p(\Omega^1(Y, adP)) | \ d_a^* \alpha = 0\},$$

where $d_a^*$ is the $L^2$-adjoint of $d_a$ (covariant derivative on sections of $adP$) with respect to some metric on $Y$.

The Chern-Simons functional $cs : \mathcal{A}(P) \to \mathbf{R}$ is defined as

$$cs(a) = \frac{1}{2} \int_Y tr(a \wedge da + \frac{2}{3} a \wedge a \wedge a),$$

and satisfies $cs(g \cdot a) = cs(a) + 2\pi deg(g)$ for gauge transformations $g : Y \to SU(2)$. Thus $cs$ is well-defined on $\tilde{\mathcal{B}}(P) = \mathcal{A}(P)/\{g \in \mathcal{G} : deg(g) = 0\}$ and it descends to a functional $cs$ which plays the role of a Morse function in defining Floer homology. Its differential is given by

$$dcs(a)(\alpha) = \int_Y tr(F_a \wedge \alpha),$$

hence its critical set consists of the flat connections $\mathcal{R}(\mathcal{B}(P)) = \{a \in \mathcal{B}(P) | \ F_a = 0\}$. (Here $F_a$ is the curvature 2-form on $Y$) It is well-known that elements of $\mathcal{R}(\mathcal{B}(P))$ are in 1-1 correspondence with those of

$$\mathcal{R}(Y) = Hom(\pi_1(Y), SU(2))/adSU(2),$$

the $SU(2)$-representations of $\pi_1(Y)$ modulo conjugacy. Given any metric on $Y$, the Hodge star operator applied to the curvature $F_a$ gives a vector field

$$f(a) = \star F_a \in L_{k-1}^p(\Omega^1(Y, adP)).$$

Comparing with $T_a \mathcal{B}^*(P)$, we note the different Sobolev norm and denote the latter by $\mathcal{L}_a$. Hence $f$ is a section of the bundle with fiber $\mathcal{L}_a$. A representation $\alpha \in \mathcal{R}(Y)$ is called *nondegenerate* if the twisted cohomology $H^1(Y; ad\alpha)$ vanishes. Note that this is the same as requiring that $\ker df(a) = \ker \star d_a = 0$, where $\star d_a$ is the Hessian of the Chern-Simons functional.

A 1-parameter family $\{a(t)| \ t \in \mathbf{R}\}$ of connections on $P$ gives rise to a connection $A$ with vanishing $t$-component on the trivial $SU(2)$ bundle over $Y \times \mathbf{R}$. Floer's crucial observation is that trajectories of the vector field $f$, i.e. the flow lines of

$$\frac{\partial a}{\partial t} + f(a(t)) = 0 \quad \text{or} \quad \frac{\partial a}{\partial t} = \star F(a(t)), \tag{2.1}$$

can be identified with instantons $A$ on $Y \times \mathbf{R}$, and $A|_{Y \times \{t\}} = a(t)$. A trajectory flow "connects" two flat connections on $Y$ if and only if the Yang-Mills energy of the trajectory is finite. One needs to show that all zeros of $f$ are nondegenerate and that their stable and unstable manifolds intersect transversally in smooth finite dimensional manifolds. Floer has shown that one can perturb the Chern-Simons functional to achieve this (see [13]). In fact Taubes has shown that the Wilson loop perturbations have those desired properties. (c.f. [34] Proposition 1.5 and Proposition 5.1) For the rest of this paper, we assume that the Chern-Simon functional has been so perturbed, so that all irreducible representations are isolated and nondegenerate. Since $\mathcal{R}(Y)$ is compact, they are also finite.

For the analysis of ASD connections, it is convenient to work with the weighted Sobolev space $L^p_{k,\delta}$ that we will introduce in §3. For each connection $A$ the anti-self-duality operator induces a Fredholm operator

$$d_A^* \oplus d_A^+ : L^p_{k+1,\delta}(\Omega^1(Y \times \mathbf{R}, adP)) \to L^p_{k,\delta}((\Omega^0 \oplus \Omega^2_+)(Y \times \mathbf{R}, adP)). \tag{2.2}$$

We say that $A$ is *regular* if $d_A^* \oplus d_A^+$ is surjective, i.e. $H_A^0 = 0$ (irreducible) and $H_A^2 = 0$ (generic). For a nondegenerate critical point $\alpha$ of $cs$, the spectral flow is $SF(\alpha, \theta) = Index(d_A^* \oplus d_A^+)(\alpha, \theta)$, the Atiyah-Patodi-Singer index of the anti-self-duality operator over $Y \times \mathbf{R}$. So

$$\mu(\alpha) \equiv Index(d_A^* \oplus d_A^+)(\alpha, \theta) \pmod{8}$$

where $A$ is any family of connections $\{a(t)\} \in \mathcal{B}(P)$ over $Y$ with $a(+\infty) = \theta, a(-\infty) = a_\alpha$ (see [13]). Floer's chain group $C_j(Y)$ is defined to be the free module generated by irreducible flat connections $\alpha$ with $\mu(\alpha) = j \pmod{8}$.

Define $\mathcal{M}_{Y \times \mathbf{R}}$ to be the moduli space of finite-energy ASD connections on $Y \times \mathbf{R}$ and let $\mathcal{M}(\alpha, \beta)$ be the subspace of those $A$ such that $\lim_{t \to -\infty} A = \alpha$, $\lim_{t \to +\infty} A = \beta$ for fixed flat connections $\alpha$ and $\beta$. It is a smooth, canonically oriented manifold which has dimension congruent to $\mu(\alpha) - \mu(\beta) \pmod 8$. The moduli space $\mathcal{M}(\alpha, \beta)$ has infinitely many connected components each of which admits a proper, free $\mathbf{R}$-action arising from translations in $Y \times \mathbf{R}$. If $\mu(\alpha) - \mu(\beta) = 1 \pmod 8$, let $\mathcal{M}^1(\alpha, \beta)$ be the union of 1-dimensional components of $\mathcal{M}(\alpha, \beta)$. Further perturbations make all the $\mathcal{M}^1(\alpha, \beta)$ regular. Then $\mathcal{M}^1(\alpha, \beta)/\mathbf{R}$ will be a compact oriented 0-manifold, i.e. it is a finite set of signed points. The differential $\partial : C_j \to C_{j-1}$ of Floer's chain complex is defined

by
$$\partial \alpha = \sum_{\beta \in C_{j-1}} \#\hat{\mathcal{M}}(\alpha, \beta) \beta \tag{2.3}$$

where $\hat{\mathcal{M}}(\alpha, \beta) = \mathcal{M}^1(\alpha, \beta)/\mathbf{R}$, and $\#\hat{\mathcal{M}}(\alpha, \beta)$ is the algebraic number of points. The sign in this formula is given by the spectral flow. Floer has shown that $\partial^2 = 0$, hence $\{C_j, \partial\}_{j \in Z_8}$ is a chain complex graded by $\mathbf{Z}_8$. The homology of this complex is called the Floer homology, denoted by $HF_j$. Floer has shown that it is independent of the choice of metric on $Y$ and of perturbations (see [9], [13]).

## 2.2 Floer homology for Lagrangian intersections

In this subsection, we briefly recall the definition of Floer symplectic homology of Lagrangian intersections and Oh's extension. Then we will discuss the well-definedness of Floer homology for the representation spaces. For details on the background see [15], [29], and [16].

In [14] and [15], Floer studied the symplectic homology theory for Lagrangian intersections with the assumption $\pi_2(P, L) = 0$, where $(P, \omega)$ is a symplectic manifold of dimension $2n$ and $L$ is a Lagrangian submanifold. For any exact diffeomorphism $\phi$ of $P$ with the property that $\phi(L)$ intersects $L$ transversally, Floer obtained a symplectic homology for Lagrangian intersections to prove the Arnold conjecture. But the condition $\pi_2(P, L) = 0$ is too restrictive for most situations.

**Example of $\pi_2(P, L) \neq 0$:** Let $Y$ be a homology 3-sphere and $Y = Y_1 \cup_\Sigma Y_2$ be a Heegaard decomposition where $\Sigma$ is a Riemann surface. The representation space $\mathcal{R}(\Sigma)$ of $\pi_1(\Sigma)$ in $SU(2)$ carries a natural symplectic structure, but it is not a smooth manifold. $\mathcal{R}(Y_i), i = 1, 2$ is a Lagrangian subspace, not a smooth one. ¿From simple calculation, we have $\pi_2(\mathcal{R}(\Sigma), \mathcal{R}(Y_i)) = Z$, for $i = 1, 2$.

### 2.2.1 Oh's extension of Floer homology for Lagrangian intersections

Thus one needs to extend the Floer homology of Lagrangian intersections which drops the unpleasant condition on $\pi_2$ and careful anaylsis on the stratified space such as the representation space. Y. Oh [29] has generalized to the situation which is more suitable to our case. We will review his construction for the Floer symplectic homology of Largrangian intersections without assuming $\pi_2(P, L) = 0$, and anaylize the singularities in the representation space which is just a matter of technicalities for homology 3-spheres (see [29] §6).

Let $(P, \omega)$ be a compact symplectic manifold with symplectic structure $[\omega] \in H^2(P, \mathbf{R})$ a nontrivial second cohomology class. By choosing an almost complex structure $J$ on $(P, \omega)$ such that $\omega(\bullet, J\bullet)$ defines a Riemannian metric, we have an integer valued second cohomology class $c_1(P) \in H^2(P, Z)$ the first Chern class. These classes define two homomorphisms

$$I_\omega : \pi_2(P) \to \mathbf{R}; \qquad I_{c_1} : \pi_2(P) \to Z.$$

If $u : (D^2, \partial D^2) \to (P, L)$ is a smooth map of pairs, there is a unique trivialization up to homotopy of the pull-back bundle $u^*TP \cong D^2 \times C^n$ as a symplectic vector bundle. This trivialization defines a map from $S^1 = \partial D^2$ to $\Lambda(C^n)$ the set of Lagrangians in $C^n$. Let $\mu \in H^1(\Lambda(C^n), Z)$ be the well-known Maslov class. Then we define a map

$$I_{\mu, L} : \pi_2(P, L) \to Z,$$

by $I_{\mu, L}(u) = \mu(\partial D^2)$, this Maslov index is invariant under any symplectic isotopy of $P$.

**Definition 2.2.1** *(i) $(P, \omega)$ is a monotone symplectic manifold if*

$$I_{c_1} = \alpha I_\omega, \quad \text{for some } \alpha > 0$$

*(ii) A Lagrangian submanifold $L$ on $P$ is monotone if*

$$I_{\mu, L} = \lambda I_\omega, \quad \text{for some } \lambda > 0.$$

**Remark:** The monotonicity is preserved under the exact deformations of $L$. By the canonical homomorphism $f : \pi_2(P) \to \pi_2(P, L)$, one has

$$I_\omega(x) = I_\omega(f(x)), \qquad I_{\mu, L}(f(x)) = 2I_{c_1}(x),$$

where $x \neq 0 \in \pi_2(P)$. Therefore if $L$ is a monotone Lagrangian submanifold, then $P$ must be a monotone symplectic manifold and $\lambda = 2\alpha$. In fact the constant $\lambda$ does not depend on the Lagrangian $L$, but on $(P, \omega)$ if $I_\omega|_{\pi_2(P)} \neq 0$.

The following proposition (c.f. Proposition 2.7 in [29]) is a crucial ingredient to prove compactness property for the holomorphic trajectories connecting two intersection points of Lagrangian submanifolds.

**Proposition 2.2.2** *Suppose that (1) $L_i, i = 1, 2$ are monotone Lagrangian submanifolds, (2) Image $Im\pi_1(L_i) \subset \pi_1(P)$ under the canonical homomorphism is a torsion subgroup for at least one of $L_i$, (3) $u_i : [0, 1] \times [0, 1] \to P, i = 1, 2$ be maps such that*

$$u_i(\bullet, j-1) \in L_j, \quad u_i(0, \bullet) = x, \quad u_i(1, \bullet) = y; \quad j = 1, 2; x, y \in L_1 \cap L_2.$$

*Then (i) $[\omega](u_0) = [\omega](u_1)$ if and only if $\mu_{u_0}(x, y) = \mu_{u_1}(x, y)$, where $\mu_u$ is the Maslov-Viterbo index.*

*(ii) If $u_i$ are J-holomorphic with respect to an almost complex structure J compatible with $\omega$, then*

$$\int \|\nabla u_0\|_J^2 = \int \|\nabla u_1\|_J^2 \quad \text{if and only if } \mu_{u_0}(x, y) = \mu_{u_1}(x, y).$$

Then define a trajectory $u$ as the solution of the Cauchy-Riemann equation

$$\overline{\partial}_J u = \frac{\partial u}{\partial \tau} + J_t \frac{\partial u}{\partial t} = 0, \tag{2.4}$$

where $(\tau, t)$ is a coordinate in $[0, 1] \times [0, 1]$. For a monotone Lagrangian $L$, let $\sigma(L)$ be the positive generator for the subgroup $[\mu|_{\pi_2(P,L)}]$ of $Z$. Under the following assumptions, Oh in [29] extends the Floer symplectic homology for Lagrangian intersections **without** $\pi_2(P, L_i) = 0$:

**Assumptions:**

1. $P, L_1, L_2$ are monotone in the sense of Definition 2.2.1.

2. $\sigma(L_i) \geq 3, i = 1, 2$ (for compactness reason).

3. $Im(\pi_1(L_i)) \subset \pi_1(P)$ is a torsion subgroup for at least one of $L_i$ (for regularity reason).

4. Lagrangian submanifolds $L_1, L_2$ intersect transversally.

A trajectory flow connects two points of Lagrangian intersections if and only if the symplectic action of the trajectory is finite (see [13]). ¿From the assumption (4), the moduli space of holomorphic curves $\mathcal{M}_J(x, y)$ is a smooth manifold with dimension $\mu_u(x, y)$, which admits a proper free **R**-action arising from translation in the $t$-direction. Then the one dimensional components $\hat{\mathcal{M}}_J(x, y) = \mathcal{M}_J(x, y)/\mathbf{R}$ will be a compact oriented 0-dimensional manifold, i.e., it is a finite set of signed points. The symplectic chain group $C_*^{\text{sym}}(L_1, L_2; P)$ is defined by the free module generated by the transversal intersection points in $L_1 \cap L_2$. The Floer symplectic boundary map

$$\partial^{\text{sym}} x = \sum_{\mu_u(x,y)=1} \#\hat{\mathcal{M}}_J(x, y) \cdot y, \tag{2.5}$$

makes $(C_*^{\text{sym}}(L_0, L_1; P), \partial^{\text{sym}})$ chain complex via $\partial^{\text{sym}} \circ \partial^{\text{sym}} = 0$ and its homology is denoted by $HF_*^{\text{sym}}(L_0, L_1; P)_J$. Note that $HF_*^{\text{sym}}(L_0, L_1; P)_J$ has a graded abelian group structure with grading given by $g.c.d(\sigma(L_1), \sigma(L_2))$. Floer and Oh have shown that this symplectic Floer homology is independent of the choice of $J$ and hamiltonian perturbations (see [15] and [29] Theorem 1.2).

### 2.2.2 Singularities on the representation spaces

Let $(Y, Y_1, Y_2, Y_0)$ denote a Heegaard decomposition of $Y$, where $Y_1 = Y_+ \cup \Sigma \times [0, 1]$ and $Y_2 = Y_- \cup \Sigma \times [-1, 0]$, so that they have overlap

$$Y_0 = Y_+ \cap Y_- = \Sigma \times [-1, 1].$$

Since both $Y_1$ and $Y_2$ are handle bodies, their fundamental groups are free groups in $g$ generators ($g$ = genus of $\Sigma$) and $\pi_1(Y_0) \cong \pi_1(\Sigma)$. (c.f. [1], [6], [24], [34] and [35]) There is a pull-back diagram of representation spaces

$$\begin{array}{ccc} \mathcal{R}(Y) & \to & \mathcal{R}(Y_1) \\ \downarrow & & \downarrow \\ \mathcal{R}(Y_2) & \to & \mathcal{R}(Y_0) \end{array}$$

from Seifert-Van Kampen theorem. In other words, the natual mappings $\mathcal{R}(Y_i) \to \mathcal{R}(Y_0)$, $i = 1, 2$, can be regarded as inclusions of subspaces $\mathcal{R}(Y_1), \mathcal{R}(Y_2)$ into $\mathcal{R}(Y_0)$ and their intersection is $\mathcal{R}(Y)$. In order to get a well-defined Floer symplectic homology of the representation spaces from the Heegaard decomposition of $Y$, one has to understand those four assumptions for our particular situation and to understand the effect of the reducible representation stratum. Note that $\pi_1(\mathcal{R}(Y_0)) = 0$ for genus $g \geq 3$, so the assumption (3) that $Im(\pi_1(\mathcal{R}(Y_i))) \subset \pi_1(\mathcal{R}(Y_0))$ is a torsion subgroup of at least one of $i$ is automatically satisfied. Huebschmann in [22] has verified the following stratified symplectic structure is indeed in the sense of Sjamaar and Lerman in [33].

The following theorem gives the stratified structure for our case. (see [1] and [35])

**Theorem 2.2.3** *(i) $\mathcal{R}(Y_j)$ is a stratified space, $j = 0, 1, 2, \emptyset$. The top strata, denoted by $\mathcal{R}^*(Y_j)$, consists of all the irreducible $SU(2)$-representations of $\pi_1(Y_j)$ which have dimension $6g - 6, 3g - 3, 3g - 3$ and $0$ respectively.*

*(ii) The other two singular strata $\mathcal{R}(Y_j) \supset \mathcal{S}(Y_j) \supset \mathcal{P}(Y_j)$,*

$$\mathcal{S}(Y_j) = Hom(\pi_1(Y_j), U(1))/Z_2, \quad j = 0, 1, 2,$$

*consists of representations with image in the $U(1)$-subgroup and $Z_2$-action sends an $U(1)$ representation $\rho : \pi_1(Y_j) \to U(1)$ to its complex conjugate $\overline{\rho}$, which have dimension $2g, g$, and $g$ respectively.*

$$\mathcal{P}(Y_j) \cong Hom(\pi_1(Y_j), Z_2) \cong H^1(Y_j, Z_2),$$

*consists of representations into the center $Z(SU(2)) \cong \{\pm I\}$ of $SU(2)$. (see [35])*

¿From deformation theory, the Zariski tangent space of $\mathcal{R}^*(Y_j)$ at $\rho \in \mathcal{R}^*(Y_j)$ is naturally isomorphic to the cohomology $H^1(Y_j, Ad\rho)$, but it may not be the actual tangent space because there are obstructions for the existence of deformation in $H^2(Y_j, Ad\rho)$. An element $\alpha \in H^1(Y_j, ad\rho)$ is tangent to a curve in $\mathcal{R}(Y_j)$ if and only if $[\alpha \wedge \alpha] = 0$. From the viewpoint of symplectic geometry, this corresponds to the fact that the intersection of $\mathcal{R}^*(Y_1)$ and $\mathcal{R}^*(Y_2)$ are not necessarily clean intersection. The symplectic structure on $\mathcal{R}^*(Y_0)$ comes from a nondegenerate skew symmetric pairing $\omega$ defined by

$$\omega : H^1(Y_0, Ad\rho) \times H^1(Y_0, Ad\rho) \to H^2(Y_0, Ad\rho \otimes Ad\rho) \to H^2(Y_0, \mathbf{R}) \cong \mathbf{R},$$

where the first arrow is cup product and the second is induced by the Killing form $Ad\rho \otimes Ad\rho \to \mathbf{R}$ on the Lie algebra $Ad\rho \cong \mathbf{su}(2)$. Since $H^2(Y_j, \mathbf{R}) = 0$ for $j = 1, 2$, the restriction of $\omega$ to the subspace $H^1(Y_j, Ad\rho)$ is trivial, and since $H^1(Y_j, Ad\rho)$ has half of the dimension of $H^1(Y_0, Ad\rho)$ it is a Lagrangian subspace in $H^1(Y_0, Ad\rho)$. In other words, $\mathcal{R}^*(Y_1)$ and $\mathcal{R}^*(Y_2)$ are Lagrangian submanifolds. In general, the dimension $H^1(Y, Ad\rho)$ may jump from point to point depending on the intersection of the two Lagrangian subspaces $H^1(Y_1, Ad\rho), H^1(Y_2, Ad\rho)$ as they fit into a Mayer-Vietoris sequence:

$$0 \to H^1(Y, Ad\rho) \to H^1(Y_1, Ad\rho) \oplus H^1(Y_2, Ad\rho) \to H^1(Y_0, Ad\rho) \to 0.$$

Let $\mathcal{S}^*(Y_j) = \mathcal{S}(Y_j) \setminus \mathcal{P}(Y_j), j = 0, 1, 2, \emptyset$, denote the complement of $\mathcal{P}(Y_j)$ in $\mathcal{S}(Y_j)$. Then, similar to the situation of irreducible representations, $\mathcal{S}^*(Y_0)$ is a nonsingular symplectic manifold of dimension $2g, g = $ genus of $\Sigma$, and $\mathcal{S}^*(Y_1), \mathcal{S}^*(Y_2)$ are lagrangian submanifolds of $\mathcal{S}^*(Y_0)$ with $\mathcal{S}^*(Y)$ as their common intersection. ¿From Mayer-Vietoris exact sequence

$$0 \to H^1(Y, U(1))(= 0) \to H^1(Y_1, U(1)) \oplus H^1(Y_2, U(1)) \to H^1(Y_0, U(1)) \to 0,$$

the intersection of $\mathcal{S}^*(Y_1)$ and $\mathcal{S}^*(Y_2)$ are always clean intersection.

**Definition 2.2.4** *The moment map $\mu : H^1(Y_0, Ad\rho) \to Hom(\mathbf{h}, \mathbf{R})$ is defined by sending $\xi \in H^1(Y_0, Ad\rho)$ to the homomorphism $\eta \mapsto <\xi, \frac{d}{dt}exp(t\eta) \cdot \xi > |_{t=0}$, where $\eta \in \mathbf{h} = R, exp(t\eta) \cdot \xi$ is the action of $exp(t\eta)$ on $\xi$ and $<,>$ is the symplectic pairing on $H^1(Y_0, Ad\rho)$.*

Taking the zero set $\mu^{-1}(0)$ of $\mu$ and factoring out the action of the isotropy subgroup $Z(\rho)$, the quotient space $\mu^{-1}(0)/Z(\rho)$ is isomorphic to a neighborhood of $\rho$ in $\mathcal{S}^*(Y_0)$. In the present situation, we can work out $\mu^{-1}(0)/Z(\rho)$ explicitly by presenting $\rho$ in the form of diagonal matrices

$$\rho(x) = \begin{pmatrix} \sigma(x) & 0 \\ 0 & \sigma^{-1}(x) \end{pmatrix} \quad \sigma(x) : \pi_1(Y_0) \to U(1).$$

**Proposition 2.2.5** *(i) $\mathcal{S}^*(Y_0)$ has the tangent bundle fibre $H^1(Y_0, \mathbf{R}) \times \{0\}$ and normal bundle fibre $c(S^{2g-3} \times S^{2g-3})/U(1)$ in $\mathcal{R}^*(Y_0)$, since the moment map $\overline{\mu} : H^1_{\overline{\partial}}(\Sigma, \sigma^{\otimes 2}) \times H^1_{\overline{\partial}}(\Sigma, \overline{\sigma}^{\otimes 2}) \to \mathbf{R}$ is given by $\overline{\mu}(x,y) = -\|x\|^2 + \|y\|^2$.*

*(ii) $\mathcal{P}(Y_0)$ has normal bundle fibre in $\mathcal{R}^*(Y_0)$ equal to $((\mathbf{R}^{2g} \otimes SU(2))_0/SU(2))$, and normal bundle fibre in $\mathcal{S}^*(Y_0)$ equal to $((\mathbf{R}^{2g} \otimes SU(2))_d)/SU(2) = \mathbf{R}^{2g}/Z_2$, where $(\mathbf{R}^{2g} \otimes SU(2))_d$ denote the set of decomposable elements and $(\mathbf{R}^{2g} \otimes SU(2))_0 = (\mathbf{R}^{2g} \otimes SU(2)) \setminus (\mathbf{R}^{2g} \otimes SU(2))_d$.*

The above proposition gives the local model for singular points in $\mathcal{R}(Y_0)$. In particular, the angle of the normal cone in $\mathcal{R}^*(Y_0)$ is $\pi/2$. Now we define the strata monotonicity for stratified symplectic and Lagrangian spaces.

**Definition 2.2.6** *(i) The stratified symplectic space $(P = \cup_{i \in I} S_i, \omega)$ is monotone at the strata $S_i$ if*

$$I_{c_1}|_{S_i} = \alpha_i I_\omega|_{S_i}$$

*is true for the strata $S_i$ for some $\alpha_i > 0$.*

*(ii) The stratified Lagrangian subspace $L = \cup_{i \in I} L_i$ is monotone at strata $L_i$ if*

$$I_{\mu, L_i} = \lambda_i I_\omega|_{L_i}$$

*holds for the strata $L_i$ of $L$ for some $\lambda_i > 0$.*

Now we show that the representation spaces are monotone in the stratumwise sense.

**Lemma 2.2.7** *The top-strata $\mathcal{R}^*(Y_0)$ in $\mathcal{R}(Y_0)$ is a monotone symplectic manifold. Also the top-strata $\mathcal{R}^*(Y_j)$ in $\mathcal{R}(Y_j), j = 1, 2$ is a monotone Lagrangian submanifold in $\mathcal{R}^*(Y_0)$.*

Proof: Atiyah and Bott in [4] observed that the first Chern class of tangent bundle of $\mathcal{R}^*(Y_0)$ determines an isomorphism of $\pi_2(\mathcal{R}^*(Y_0))$ with the even integer. They also determine an integrable complex structure on $\mathcal{R}^*(Y_0)$ by the Hodge star operator on $H_A^1(Y_0, adSU(2))$ which is compatible with $\omega$ the symplectic structure. For the top-strata $\mathcal{R}^*(Y_0)$, $\alpha_{\mathcal{R}^*(Y_0)} = \frac{1}{4\pi^2} > 0$ which is verified in [4] and [12].

Both $\mathcal{R}(Y_i)$ are symmtric in $\mathcal{R}^*(Y_0)$, since there exist anti-symmetric involutions $\tau_i : \mathcal{R}(Y_0) \to \mathcal{R}(Y_0)$ with $\text{Fix}(\tau_i) = \mathcal{R}(Y_i)$. In a suitable basis, we have that $\pi_1(\Sigma)$ is generated by $a_1, \cdots, a_g$, and $b_1, \cdots, b_g$ with the single relation

$$\coprod_{i=1}^{g}[a_i, b_i] = 1,$$

$\pi_1(\Sigma) \to \pi_1(Y_1)$ sends each $b_i \to 1$ and the images of the $a_i$ freely generate $\pi_1(Y_1)$. For $\pi_1(Y_2)$ a similar description holds relative to a different basis, i.e. after applying an automorphism of $\pi_1(\Sigma)$ (see Chapter II.1.(d) in [1]). So it is monotone by a result of Oh in [29]. ∎

**Remark:** The above Lemma is known to Oh [29]. One can use the induced symplectic form at each strata to define the monotonicity stratumwise. This will be a future study. For this paper, we only need to work on the top-strata monotone symplectic and lagrangian manifolds.

### 2.2.3 Floer symplectic homology for the representation spaces

Because all the representations of $Y$ are either trivial or irreducible $SU(2)$-representations (for integral homology 3-spheres), we can concentrate on the Lagrangian intersection points in $L_1 \cap L_2$ which are irreducible. Any $J$-holomorphic disk with $\{Id\}$ as its corner point will have dimension $> 2$ (see §3.2). We consider only the Floer homology boundary map which involves one and two dimensional moduli space. So the ad-trivial representations do not come into our discussion.

As the space in question is a stratified space, the question of $J$-holomorphic curve and its regularity may seem to be ambigous as in the case of minimal surfaces [8]. For the representation spaces, we make use of the identification of geometric quotient and symplectic quotient for the reducible strata. Then the definition of $J$-holomorphic curve at reducibles is given in the principal of $U(1)$-lifting of $J$-holomorphic curve on reducible strata [33]. Since the geometric quotient gives the precise complex structures of ambient space and its quotient, so we define the $J$-holomorphic curve through the geometric quotient which is transverse to the $\mathbf{C}^*$-orbits. The local structure of $J$-holomorphic curve at reducible representation is studied in order to give a well-defined symplectic Floer homology of representation spaces. This kind of analysis has been done in the gauge theory for instanton connections by Fukaya in [19].

Let us clarify the almost complex structure on the normal direction first. The vector spaces $\{H^1(Y_0, \mathbf{h}_{Ad\rho}^\perp) \mid p \in \mathcal{S}^*(Y_0)\}$ form a symplectic vector bundle $\nu$ over $\mathcal{S}^*(Y_0)$. There is a Hermitian structure on $\nu$ compatible with its symplectic structure $\omega$. After picking a complex structure on $\Sigma$ and identifying $H^1(Y_0, \mathbf{h}_{Ad\rho}^\perp)$ with $H_{\bar{\partial}}^1(\Sigma, \sigma^{\otimes 2} \otimes C)$, we have a complex structure on each of these vector spaces with $U(1)$-action. Note that $H^1(Y_j, Ad\rho)$ is decomposed into the sum $H^1(Y_j, \mathbf{R}) \oplus H^1(Y_j, \sigma^{\otimes 2})$, where $H^1(Y_j, \mathbf{R})$ and $H^1(Y_j, \sigma^{\otimes 2})$ are respectively the real and complex Lagrangian subspaces in $H^1(Y_0, \mathbf{R})$ and $H^1(Y_0, \sigma^{\otimes 2})$. By **complex Lagrangian**, we mean a totally real subspace in the complex symplectic space $H^1(Y_0, \sigma^{\otimes 2})$ which is a Lagrangian and is invariant under the $U(1)$-action. From the definition of moment map $\mu$, the quotient

$$H^1(Y_j, Ad\rho)/U(1) = H^1(Y_j, \mathbf{R}) \times \{H^1(Y_j, \sigma^{\otimes 2})/U(1)\},$$

is isomorphic to a neighborhood of $\rho$ in $\mathcal{R}(Y_j)$. The first factor $H^1(Y_j, \mathbf{R}) \times \{0\}$ is mapped to an Euclidean neighborhood in $\mathcal{S}^*(Y_j)$ and the second factor $\{0\} \times (H^1(Y_j, \sigma^{\otimes 2})/U(1))$ which is a cone over the complex projective space $CP^{2g-2}$ is mapped into the intersection of $\mathcal{R}(Y_j)$ with the normal cone $\bar{\mu}^{-1}(0)/U(1)$.

The moment map in definition 2.2.4 defines a Lie algebra homomorphism $\eta \to H_\eta$, where $H_\eta(\xi) = \mu(\xi)(\eta)$ is Hamiltonian with $U(1)$-action. Since the zero element $0 \in \mathbf{h}$ is a fixed point of the coadjoint action, its inverse image $\mu^{-1}(0)$ is invariant under $U(1)$ and is a submanifold of $\mathcal{R}(Y_0)$. So $\mu^{-1}(0)$ is coisotropic and the corresponding isotropic foliation is given by the orbits of $U(1)$, i.e. the leaves of the null foliation of $\omega|_{\mathcal{S}^*(Y_0)}$ are the $U(1)$-orbits. Thus $\mathcal{S}^*(Y_0) = \mu^{-1}(0)/U(1)$ is a symplectic manifold called the symplectic quotient with $\pi^*\omega|_{\mathcal{S}^*(Y_0)} = \omega|_{\mu^{-1}(0)}$, where $\pi : \mu^{-1}(0) \to \mathcal{S}^*(Y_0)$ is the orbit map. (see [33])

Based on the stratified structure, Sjamaar and Lerman in [33] showed that Hamiltonian flows are strata-preserving and gave a recipe for lifting a reduced Hamiltonian flow to the level set $\mu^{-1}(0)$. This provides that the Wilson loop perturbations can be lifted since they are hamiltonian. Because of $U(1)$-equivariant the lift of Hamiltonian flow preserves the orbit-type stratification of symplectic space. We may assume that the domain of $u$ is a closed disc $D = D_r$ centered at $0$ and of radius $r$ in $\mathbf{C}$, and the neighborhood of $u(0) \in \mathcal{S}^*(Y_0)$ is $\mathbf{C}^{3g-3}$ equipped with an induced complex structure $J$.

**Definition 2.2.8** $u$ is $J$-holomorphic at $u(0) \in \mathcal{S}^*(Y_0)$ if there exists a $J$-holomorphic curve $\tilde{u} : D \to \mathbf{C}^* \cdot \mu^{-1}(0)$ which is transverse to the $\mathbf{C}^*$-orbits and $Im(\tilde{u}) = \pi^*(Imu)$, where $\pi : \mathbf{C}^* \cdot \mu^{-1}(0) \to \mathbf{C}^* \cdot \mu^{-1}(0)//\mathbf{C}^*(= \mu^{-1}(0)/U(1))$ is the orbit map.

¿From Geometric invariant theory point of view, there is a correspondence of complex structures between the ambient space $\mathbf{C}^* \cdot \mu^{-1}(0)$ and its geometric quotient $\mathbf{C}^* \cdot \mu^{-1}(0)//\mathbf{C}^*$. The identification between geometric and symplectic quotients ([23] 7.5) make the definition of $J$-holomorphic curve natural through the geometric quotient viewpoint, the infinite dimensional version has been studied in [4] for Yang-Mills connections over Riemann surface. The advantage of geometric quotient is the globally holomorphic structure defined over $\mathbf{C}^* \cdot \mu^{-1}(0)//\mathbf{C}^*$.

Denote $\partial$ and $\overline{\partial}$ for the usual del and delbar operators on $\mathbf{C}^{3(g-1)}$. Let $u_i^l, \overline{u}_i^l, i = 1, \cdots, 3(g-1)$ for the components of $\tilde{u}$ and its conjugates. We only consider the local neighborhood of $u(0) \in \mathcal{S}^*(Y_j), j = 0, 1, 2$ since otherwise it has been explained in [28].

**Lemma 2.2.9** $u$ is $J$-holomorphic at $u(0) \in \mathcal{S}^*(Y_j)$ if and only if

$$\overline{\partial} u_i^l + A_{im}^l(u^l(z))\partial u_m^l = 0,$$

where for each $w \in \mathbf{C}^{3(g-1)}$, $A_{im}^l(u^l(z))$ is a certain $3(g-1) \times 3(g-1)$ complex valued matrix which is $\mathbf{C}^*$-invariant from the entries of $J(w)$ on $\mathbf{C}^* \cdot \mu^{-1}(0)$ and vanishes when $J(w)$ is standard complex structure on $\mathbf{C}^{3(g-1)}$. Thus $A_{im}^l(0) = 0$ for all $i, m$.

Proof: This is Lemma 2.1 in [28] with $\mathbf{C}^{3(g-1)}$ and the (almost) complex structure $J$ is induced from the geometric quotient. Then the argument applied to the lifted holomorphic curve $\tilde{u}$ with property of transversing the $\mathbf{C}^*$-orbits. ■

Based on Lemma 2.2.9, one can have the local behavior for $J$-holomorphic curve at singular point $u(0) \in \mathcal{S}^*(Y_j), j = 0, 1, 2$. In particular one can use Aronszajn's unique continuation theorem (see [28]) on $\tilde{u}$. With the above definition and Lemma 2.2.9, a $J$-holomorphic curve is regular if it is smooth with respect to the local coordinates and belongs to a Sobolev class in $L_k^p$ as in Lemma 2.2 [14] and Lemma 3.2 [29]. Because regularity is a local problem, the lifting holomorphic curve $\tilde{u}$ which is transverse to $\mathbf{C}^*$-orbits can be shown to be smooth by the standard elliptic estimates and bootstrapping argument (see [14], [30] and [33]).

There is a problem whether the 1- or 2-dimenional components of $\mathcal{M}_J(a,b)$ with $\mu(a,b) = 1, 2$ are manifolds. When the linearized operator $E_u$ of $J$-holomorphic curve equation is onto, the space

$$\mathcal{M}_{1,\delta}^p(a,b) = \{u + \xi | \xi \in L_{1,\delta}^p(u^*TP), u + \xi \text{ flat on the handle bodies } Y_1, Y_2 \text{ and}$$

$$\text{flat on } \Sigma \times \{s\} \quad \overline{\partial}_J(u+\xi) = 0, \quad \lim_{t \to -\infty}(u+\xi) = a, \lim_{t \to +\infty}(u+\xi) = b\}, \tag{2.6}$$

is a smooth manifold [11] p34. Also there exists a Baire second category of perturbation data which make $E_u$ onto. Note that the ontoness condition does not mean $u(,s)$ irreducible for all $s, 0 \leq s \leq 1$. Denote $\mathcal{M}_J(a,b)$ the moduli space of $J$-holomorphic curves on the top-strata.

In the following, we are going to characterize the property of $J$-holomorphic curves passing through $U(1)$-strata (Lemma 2.2.12). Using the grafting technique, we describe the local model of $M_J^{\text{sing}}(a,b)$ (see (2.10)). Then we show that the 1-dimensional moduli space of $J$-holomorphic curves does not contain any point in $U(1)$-strata, i.e. $\mathcal{M}_{1,\delta}^p(a,b) = \mathcal{M}_J(a,b)$ (Proposition 2.2.17). For 2-dimensional moduli space $\mathcal{M}_{1,\delta}^p(a,b)$, $\mu(a) - \mu(b) = 2$, we show that the compactification of $\mathcal{M}_J(a,b) \subset \mathcal{M}_{1,\delta}^p(a,b)$ does not intersect with the compact closed piece $M_J^{\text{sing}}(a,b)$ (Lemma 2.2.20).

**Lemma 2.2.10** *Suppose that $u \in L_1^2([0,1] \times \mathbf{R}, \mathcal{R}(Y_0))$ satisfy $\int_{[0,1] \times \mathbf{R}} |\nabla u|^2 \leq M$ for some $M$. Then for $0 < r < 1$*

$$L(u|_{\partial D_r}) = \int_{\partial D_r} |\frac{\partial u}{\partial \theta}| d\theta,$$

*is defined almost everywhere in $[0,1]$.*

Proof: Let us first assume that $(s,t)$ is in the interior of $[0,1] \times \mathbf{R}$. Hence the following estimates holds for $0 \leq \varepsilon \leq 1$.

$$\int_{\varepsilon^2}^{\varepsilon} \frac{L(u|_{\partial D_r})^2}{r} dr \leq \int_{\varepsilon^2}^{\varepsilon} \frac{1}{r} \int_{\partial D_r} |\frac{\partial u}{\partial \theta}|^2 2\pi d\theta dr$$

$$\leq 2\pi \int_{D_\varepsilon - D_{\varepsilon^2}} |\nabla u|^2 \leq 2\pi M$$

¿From this inequality, we conclude that the function $L(u|_{\partial D_r})$ is defined almost everywhere. ∎

**Lemma 2.2.11** *Let $u \in L_1^2(D_r, \mathcal{R}(Y_0)) \cap \mathcal{M}_{1,\delta}^p(a,b)$. Then there exists a constant $C_0$ such that*

$$\int_{D_r} |\nabla u|^2 \leq C_0 r L(u|_{\partial D_r}).$$

Proof: Since $u \in \mathcal{M}_{1,\delta}^p(a,b)$, we have $\int |\nabla u|^2 < \infty, \overline{\partial}_J u = 0$, so $u$ is continuous by the Sobolev embedding theorem for $L_1^p, p \geq 2$. Hence for some constant $C_1 > 0$

$$|u(s,t) - u(p)| \leq C_1 r. \tag{2.7}$$

Note that $\overline{\partial}_J u = 0$ gives $|\frac{\partial u}{\partial s}| = |\frac{\partial u}{\partial t}|$, also $|\frac{\partial u}{\partial r}| = |\frac{\partial u}{\partial \theta}|$ for polar coordinate $(r, \theta)$.

$$\begin{aligned}
\int_{D_r} |\nabla u|^2 &= \int_{D_r} <\nabla u, \nabla(u - u(p))> \\
&= \int_{\partial D_r} <\frac{\partial u}{\partial r}, (u - u(p))> d\theta \\
&\leq \int_{\partial D_r} |\frac{\partial u}{\partial r}||u - u(p)| d\theta \\
&\leq C_0 r \int_{\partial D_r} |\frac{\partial u}{\partial \theta}| d\theta = C_0 r L(u|_{\partial D_r}).
\end{aligned}$$

The second equality follows from integration by parts (see Thereom 2 in [8] p266). ∎

**Remarks:** (1) In minimal surface theory, this type of estiamte is called a linear isoperimetric inequlity (c.f. [8] §6.3). Both isoperimetric inequality and Courant-Lebesgue lemma (used in [30]) play important roles in minimal surface theory ([8] Chapter 4).

(2) Any $J$-holomorphic curve is energy minimizing in its homological class relative to the free Lagrangian boundary condition (see Proposition 2.3 [15]). Thus the image of the curve is a minimal surface (see [29]). In [2], Almgren studied the general regularity problem and described the singular set in terms of Hausdorff measure. For our purpose, one does not need to follow the geometric measure theory treatment, even though $\mathcal{R}(Y_j)$ is a rectifiable set (see §2.3 [2]).

**Lemma 2.2.12** *Suppose that $u$ is a $J$-holomorphic curve in $L_1^2([0,1] \times \mathbf{R}; \mathcal{R}(Y_0))$ with $u(i-1, t) \subset \mathcal{R}(Y_i), i = 1, 2$ which can be reducible at $(s,t)$, i.e. $u(s,t) \in \mathcal{S}^*(Y_0)$. Then $\nabla u(s,t) = 0$.*

Proof: From symplectic quotient point of view, there is a normal cone of $u(s,t)$ by Proposition 2.2.5. Then there exist $C_2 > 0$ such that

$$\frac{L(u|_{\partial D_r})}{|u|_{\partial D_r} - u(s,t)|} \leq C_2 2\pi \sin(\frac{\theta}{2}) < \infty, \qquad (2.8)$$

where $\theta$ is the normal cone angle. By Lemma 2.2.11 and (2.7), we have

$$\frac{L(u|_{\partial D_r})}{|u|_{\partial D_r} - u(s,t)|} \geq \frac{\int_{D_r} |\nabla u|^2}{C_0 C_1 r^2}. \qquad (2.9)$$

Note that $u(D_r)$ is inside the normal cone at $u(s,t)$ for $r$ sufficiently small. Since $\int_{D_r} |\nabla u|^2 \leq M$, so that

$$\lim_{r \to 0} \int_{D_r} |\nabla u|^2 = 0.$$

$$\lim_{r \to 0} \frac{\int_{D_r} |\nabla u|^2}{r^2} = \lim_{r \to 0} \frac{\int_0^r \int_0^{2\pi} |\frac{1}{r} \frac{\partial u}{\partial \theta}|^2 r d\theta dr}{r^2}$$
$$= \lim_{r \to 0} \frac{\int_0^{2\pi} |\frac{\partial u}{\partial \theta}|^2 d\theta}{2r^2} + \frac{\int_0^r \int_0^{2\pi} (\frac{2}{r} < \frac{\partial^2 u}{\partial r \partial \theta}, \frac{\partial u}{\partial \theta} > -\frac{1}{r^2}|\frac{\partial u}{\partial \theta}|^2)}{2r}$$

Since $u$ is smooth in $D_r$, if $p$ is not a critical point of $u$ (i.e. $\lim_{r \to 0} \frac{\partial u}{\partial \theta} \neq 0$), then we have

$$\lim_{r \to 0} \frac{L(u|_{\partial D_r})}{|u|_{\partial D_r} - u(p)|} = \infty$$

which contradicts with (2.8). ∎

To get a well-defined Floer homology, we need to consider the following space with $\mu(a,b) \leq 2$,

$$M_J^{\text{sing}}(a,b) = \{u : [0,1] \times \mathbf{R} \to \mathcal{R}(Y_0) | \ \overline{\partial}_J u = 0, u(s, \pm\infty) \in \mathcal{R}^*(Y) \quad \text{there exists}$$

$$\text{a point } (s,t) \in [0,1] \times \mathbf{R} \text{ such that } \nabla u(s,t) = 0, u(s,t) \in \mathcal{S}^*(Y_0)\}. \qquad (2.10)$$

The definition of $M_J^{\text{sing}}(a,b)$ also includes the case for $s = 0,1$ with $u(0,t) \in \mathcal{S}^*(Y_1), u(1,t) \in \mathcal{S}^*(Y_2)$. The special space $M_J^{\text{sing}}(a,b)$ with $\mu(a,b) = 2$ may be empty, isolated points. We are going to study the space of balanced $J$-holomorphic curves $\hat{M}_J^{\text{sing}}(a,b)$.

The local model for $M_J^{\text{sing}}(a,b)$ is very similar to the one for instanton with $U(1)$-reducible flat connection in the middle (see §7 in [19]). The local $U(1)$- group action around the $U(1)$-reducible connection can not be extended to be a global $U(1)$-group action since the element $u \in M_J^{\text{sing}}(a,b)$ has the isotropy group $\{\pm 1\}$ for $a, b \in \mathcal{R}^*(Y)$. We are going to glue this local $U(1)$-parameter for $J$-holomorphic curve like instanton case in [19] §10.

Let $u \in M_J^{\text{sing}}(a,b)$ be a $J$-holomorphic curve. Note that $u|_{D_r}$ is the $J$-holomorphic curve defined in Definition 2.2.8. The different lifting $g\tilde{u}, g \in U(1)$, of $J$-holomorphic curve in $\mathbf{C}^* \cdot \mu^{-1}(0)//\mathbf{C}^* = \mu^{-1}(0)/U(1)$ gives different element in $\mathbf{C}^* \cdot \mu^{-1}(0)$. Let $g_r$ be a path in $SU(2)$ such that $g_r = Id, r \geq r_1; g_r = g, r \leq r_2$. Let

$$\exp : D_r \times T(\mathbf{C}^* \cdot \mu^{-1}(0)) \to \mathbf{C}^* \cdot \mu^{-1}(0)$$

be a smooth family of charts of the ambient space such that

$$\exp((s,t); \xi) = g\tilde{u}, \quad \exp((s,t); \zeta) = g_r^* u,$$

where $\xi(s,t) \in T(\mathbf{C}^* \cdot \mu^{-1}(0))$ is transverse vector field to the $\mathbf{C}^*$-orbits. We define an almost $J$-holomorphic curve using a cutoff function. There is a function $\varepsilon : (0, r_1] \to \mathbf{R}_+$ with $\lim_{r \to 0} \varepsilon(r) = 0$ and

$$\|\xi|_{D_r}\|_{L_1^p} \leq \varepsilon(r), \quad \|\zeta|_{D_r}\|_{L_1^p} \leq \varepsilon(r).$$

Let $\chi_r$ be the smooth cutoff function with $\chi_r = 1, r \geq r_2; \chi_r = 0, r \leq r_3$. Note that $r_3 \leq r_2 \leq r_1$.

**Definition 2.2.13** *For $u \in \hat{M}_J^{sing}(a,b)$, we define the map from $\hat{M}_J^{sing}(a,b) \times (0, r_1] \times U(1)$ to $\mathcal{P}(a,b)$ which is given by*

$$\Phi(u, r, g) = u_{\chi_r, g} = \begin{cases} u & \text{for } [0,1] \times \mathbf{R} \setminus D_{r_1} \\ g_r^* u & \text{for } D_{r_1} \setminus D_{r_2} \\ \exp(\chi_r \xi + (1 - \chi_r) \zeta) & \text{for } D_{r_2} \setminus D_{r_3} \\ g\tilde{u} & \text{for } D_{r_3} \end{cases} \quad (2.11)$$

*and $\xi, \zeta$ are defined as above (see [15] §4 and [19] §8).*

If the singular point is at $s = 0, 1$, then the disk is understood as half-disk. ¿From the definition of $\Phi(u, r, g)$, it is easy to see from [15] §4 that $\Phi(u, r, g)$ is a continuous map and its image is almost holomorphic.

**Lemma 2.2.14** *For any $u \in M_J^{sing}(a,b)$ with $\mu(a) - \mu(b) = 1$, there exists $\varepsilon(r)$, $\lim_{r \to 0} \varepsilon(r) = 0$ such that for $r \leq r_1$,*

$$\|\bar{\partial}_J \Phi(u, r, g)\|_{L_0^p} < \varepsilon(r); \quad \|\Phi(u, r, g) - u\|_{L_0^p} < \varepsilon(r). \quad (2.12)$$

**Lemma 2.2.15** *The balanced space $\hat{M}_J^{sing}(a,b)$ with $\mu(a) - \mu(b) = 1$ is compact, and there is a constant $C$ independent of $u \in \hat{M}_J^{sing}(a,b)$ such that for all $u \in \hat{M}_J^{sing}(a,b)$, $\xi \in L_0^p(u^*(T\mathcal{R}(Y_0)))$ and $p \geq 2$, we have*

$$C\|\xi\|_{L_0^p}^p \leq \|E_u^*\xi\|_{L_0^p}^p, \tag{2.13}$$

*where $E_u^*$ is the $L^2$-adjoint operator of $E_u$.*

Note that the uniformly constant holds for any compact space of regular $J$-holomorphic curves.

**Proposition 2.2.16** *For any $u \in \hat{M}_J^{sing}(a,b), \mu(a) - \mu(b) = 1$, there is a local $U(1)$-action on $\hat{M}_J^{sing}(a,b)$ and an injective map*

$$\tilde{\Phi} : \hat{M}_J^{sing}(a,b) \times (0, r_1] \times U(1) \to M_J^{sing}(a,b).$$

Proof: Let us solve $\overline{\partial}_J \tilde{\Phi}(u,r,g) = 0$ for $\tilde{\Phi}(u,r,g) = \Phi(u,r,g) + \tilde{\xi}$. The uniformly bounded right inverse of $E_u$ from Lemma 2.2.15 also gives the bound for the right inverse of $E_{\Phi(u,r,g)}$, since we have that $E_{\Phi(u,r,g)}^* - E_u^*$ is zero order compact operator and from (2.12)

$$\begin{aligned}\|(E_{\Phi(u,r,g)}^* - E_u^*)\xi\|_{L_0^p} &\leq C\|\Phi(u,r,g) - u\|_{L_0^p}\|\xi\|_{L_1^p} \\ &\leq C\varepsilon(r)\|\xi\|_{L_1^p}.\end{aligned}$$

Then we obtain the injective map $\tilde{\Phi}$ from Lemma 4.2 [15] or Lemma 6.1.2 in §6.1. ∎

We used the global property of $E_{\Phi(u,r,g)}$, rather than local information at ends as Lemma 4.3 and Lemma 5.3 in [15]. This reflects the $U(1)$-group only acting locally, not globally well-defined.

**Proposition 2.2.17** *For holomorphic curve $u \in \mathcal{M}_{1,\delta}^p(a,b)$, with $dim\mathcal{M}_{1,\delta}^p(a,b) = 1$, one has that $u(0,R) \subset \mathcal{R}^*(Y_1), u(1,R) \subset \mathcal{R}^*(Y_2)$ and $u(s,t) \in \mathcal{R}^*(Y_0)$ for all the interior points in $[0,1] \times \mathbf{R}$. I.e. the $J$-holomorphic curve $u$ does not slide inside the sigular strata.*

Proof: Suppose that there is a singular point $u(s,t) \in \mathcal{S}^*(Y_0)$. Then from Proposition 2.2.16, we take different lifting and regluing together. So the local action of $U(1)$ gives extra solution piece in $\mathcal{M}_J^s(a,b)$ which contradicts with $dim\mathcal{M}_J^s(a,b) = 1$. For singularity occured at $\mathcal{S}^*(Y_j), j = 1, 2$, the same proof applied to half disk. ∎

**Lemma 2.2.18** *The minimal Maslov number $\sigma(\mathcal{R}^*(Y_j)) = 8$ for $j = 1, 2$.*

Proof: By definition $\sigma(\mathcal{R}^*(Y_j))$ is the minimal number of $\mu(u)$ for

$$u : (D^2, \partial D^2) \to (\mathcal{R}^*(Y_0), \mathcal{R}^*(Y_j)).$$

The Maslov index $\mu(\partial D^2) = I_{(\mathcal{R}^*(Y_0), \mathcal{R}^*(Y_j))}(u)$. Let $a \in u(\partial D^2) \cap \mathcal{R}^*(Y)$. Hence $\sigma(\mathcal{R}^*(Y_j)) = 8 > 3$ for $j = 1, 2$ follows from Proposition 3.2.2 and the fact that $SF(D(a, a)) \equiv 0 \pmod{8}$ for self-duality operator $D$. ∎

**Lemma 2.2.19** *For $\mu(a) - \mu(b) = 2$, the space $\hat{M}_J^{sing}(a, b)$ is a compact, closed manifold.*

Proof: Let $u_n \in \hat{M}_J^{sing}(a, b)$ be a sequence with finite energy. Then again from $\sigma = 8 > 2$ and $\mu = 2 < 8$ there is no bubbling occured. Hence there exists a subsequence $\{u_n\}$ that has weak limit to a holomorphic curve $u$ with $\lim_{t \to -\infty} u = a, \lim_{t \to \infty} u = b$ in the $L_1^2$ sense (strong convergence on the complement of bubbling points). Such a weak limit can not be split along $c \in \mathcal{R}^*(Y)$ because we are going to show that there is a limit $p$ of singular points with $\nabla u(p) = 0$. On the other hand the balanced moduli space with $\mu(a, c) = 1$ or $\mu(c, b) = 1$ can not contain a singular point from Proposition 2.2.17. So any sequence in $\hat{M}_J^{sing}(a, b)$ has a convergence subsequence with limit in $\hat{M}_J^{sing}(a, b)$. There are $\{p_n\} \in \mathcal{S}^*(Y_j)$ such that $\nabla u_n(p_n) = 0$, then from compactness of $\mathcal{S}^*(Y_j)$ and diagonal argument, we have that $\|p_n - p\|_{C^0} \le \varepsilon$. Then around the neighborhood of $p$, the normal directional $U(1)$-invariant holomorphic curve has the following estimates from Corollary 3.4 [30].

For any $r < 1, \max_{|x| < r} |\nabla u_n(x) - \nabla u(x)| \le C_3(r) \|\nabla u_n - \nabla u\|_{L^2(D_1)}$, where $C_3(r)$ depends on $r, \varepsilon$ the smaller number from $\|p_n - p\|_{C^0} \le \varepsilon$ and $\|\nabla u_n - \nabla u\|_{L_1^2(D_1)}$, but independent of $u_n - u$. Hence

$$
\begin{aligned}
|\nabla u(p)| &\le |\nabla u_n(p) - \nabla u(p)| + |\nabla u_n(p)| \\
&\le C_3(r) \|\nabla u_n - \nabla u\|_{L^2(D_1)} + |\nabla u_n(p) - \nabla u_n(p_n)| \\
&\le C_3(r)\varepsilon + \|u_n\|_{C^1} \|p_n - p\|_{C^0}
\end{aligned}
$$

The $C^1$-norm is bounded from regularity of the normal directional $u_n$. Hence $\nabla u(p) = 0$. The result follows. ∎

Once the balanced space $\hat{M}_J^{\text{sing}}(a,b)$ is compact for $\mu(a) - \mu(b) \leq 2$, then any element $u \in \hat{M}_J^{\text{sing}}(a,b)$ has a local $U(1)$-action around the reducible representations, i.e. locally there is a extra dimension piece of the $J$-holomorphic curves.

**Lemma 2.2.20** *The compactification of $\hat{\mathcal{M}}_J(a,b)$ for $\mu(a,b) = 2$ does not intersect with $\hat{M}_J^{sing}$, i.e. $\overline{\hat{\mathcal{M}}_J(a,b)} \cap \hat{M}_J^{sing} = \emptyset$.*

Proof: Suppose that a sequence $\{u_n\} \in \hat{\mathcal{M}}_J(a,b)$ has a subsequence converging to $u \in \hat{M}_J^{\text{sing}}(a,b)$. Note that there are no bubbling occured due to the dimension and index reasons, so the subsequence (still denoted by $u_n$) converges strongly to $u$. By definition of $\hat{M}_J^{\text{sing}}(a,b)$, there exists $(s,t) \in [0,1] \times \mathbf{R}$ such that $u(s,t) \in \tilde{S}^*(Y_j)$. Thus there is extra $U(1)$-parameter in $\mathcal{M}_J^{\text{sing}}(a,b)$ which contradicts with $\dim \mathcal{M}_J^{\text{sing}}(a,b) = 2$. If $u \in M_J^{\text{sing}}(\alpha, \beta)$ for $\alpha \neq a$ or $\beta \neq b$, then we have $\mu(\alpha) - \mu(\beta) = 1$ (for having 2-dimensional piece of $J$-holomorphic curves). The spectral flow is preserved under the strong convergence, i.e. $2 = \mu_{u_n}(a,b) = \mu_u(\alpha,b)$ which contradicts with $\mu(\alpha) - \mu(b) = 1$. So we obtain the result. ∎

**Proposition 2.2.21** *There is a dense set $\mathcal{J}_d$ in $\mathcal{J}_{reg}(\mathcal{R}^*(Y_1), \mathcal{R}^*(Y_2))$ such that the one dimensional component of $\hat{\mathcal{M}}_J(a,b)$ is compact up to the splitting of two isolated trajectories for $J \in \mathcal{J}_d$,*

$$\hat{\mathcal{M}}_J(a,c) \times \hat{\mathcal{M}}_J(c,b), \quad \text{for} \quad \mu(a,c) = \mu(c,b) = 1 \quad \text{and} \quad c \in \mathcal{R}^*(Y).$$

*For $\partial^{sym} a = \sum_{\mu(a,\alpha)=1} \#\hat{\mathcal{M}}_J(a,\alpha)\alpha$, then we have*

$$\partial^{sym} \circ \partial^{sym} = 0. \tag{2.14}$$

Proof: Note that the proof is same as in Proposition 4.3 [29]. Since for a sequence $u_n \in \mathcal{M}_J(a,b) \cap \mathcal{R}^*(Y_0)$ the limit $u$ can not have a sphere bubbling off ($\mu(u_n) = 2 < 8$) or a disk bubbling off ($\sigma = 8 > 2$). Both are ruled out by dimension reason and $\sigma(\mathcal{R}^*(Y_j)) = 8 \geq 3$. There is a possibility for $u$ to split along reducible flat connection $\alpha_0 \in \mathcal{S}^*(Y)$, but $\mathcal{S}^*(Y) = \emptyset$ for integral homology 3-sphere. Hence the 1-dimensional $J$-holomorphic moduli space $\hat{\mathcal{M}}_J(a,b)$ can only split along $c \in \mathcal{R}^*(Y)$ with $\mu(a,c) = \mu(c,b) = 1$. For the 1-dimensional holomorphic curves with asymptotic values in $\mathcal{R}^*(Y)$, Proposition 2.2.17 avoid the singularities occured. Hence it gives the spaces $\mathcal{M}_J(a,c), \mathcal{M}_J(c,b)$. The sequence $\{u_n\}$ can not have a subsequence converging to

$u \in M_J^{\text{sing}}$ by Lemma 2.2.20, then the proof is completed. ∎

**Theorem 2.2.22** *With this boundary map $\partial^{sym}$ and the chain groups*

$$C_*(\mathcal{R}^*(Y_1), \mathcal{R}^*(Y_2); \mathcal{R}^*(Y_0)) = Z\{x \in \mathcal{R}^*(Y_1) \cap \mathcal{R}^*(Y_2)\},$$

*we have a well-defined Floer symplectic homology of Lagrangian intersections of representation spaces*

$$H_*(C_*(\mathcal{R}^*(Y_1), \mathcal{R}^*(Y_2); \mathcal{R}^*(Y_0)), \partial^{sym}) = HF^{sym}(\mathcal{R}^*(Y_1), \mathcal{R}^*(Y_2); \mathcal{R}^*(Y_0)),$$

*which is also denoted by $HF^{sym}(\mathcal{R}(Y_1), \mathcal{R}(Y_2); \mathcal{R}(Y_0))$ for the Heegaard decomposition of $Y$.*

In the rest of the paper, we are going to identify this Floer symplectic homology with the Floer instanton homology for homology 3-spheres. Hence the symplectic Floer homology is independent of perturbations, almost complex structures and also Heegaard decompositions.

# 3 Comparing the chain complexes

## 3.1 Perturbations

In this subsection we are going to compare the generators for Floer instanton chain groups and Floer symplectic chain groups. As we indicate in §2, the Floer instanton chain group $C_*(Y)$ is generated by irreducible flat connections over the bundle $Y \times SU(2)$, which is the intersection between $\mathcal{R}(Y_1)$ and $\mathcal{R}(Y_2)$ in $\mathcal{R}^*(Y_0)$. The Floer symplectic chain group $C_*^{sym}(\mathcal{R}(Y_1), \mathcal{R}(Y_2); \mathcal{R}^*(Y_0))$ is defined to be the free module generated by the transversal intersection points $\mathcal{R}(Y_1) \cap \mathcal{R}(Y_2)$. Hence if all the irreducible flat connections $\mathcal{R}^*(Y)$ are nondegenerate, we have actually the identification between these two groups. But one has to make sure that all the critical points of $cs$ are nondegenerate. A nondegenerate zero of $f$ defined in §2.1 is isolated (see [34] Lemma 1.2). Note that $\nabla f$ has index zero so that a suitably generic perturbation of $f$ will have isolated nondegenerated zeros in $\mathcal{B}_Y^*$.

Let $Y_j, j = 0, 1, 2,$ denote as in §2.2 the various submanifolds associated to a Heegaard decomposition of $Y (j = \emptyset)$. There is the principal $SU(2)$-bundle $P_j = Y_j \times SU(2)$ over $Y_j$ and hence the space of $C^\infty$-connections on $P_j$. Denote by $\mathcal{A}_j = \mathcal{A}(P_j)$ the completion of this last space with respect to the Sobolev $L_k^p$-norm. Then, analogous to the pull-back digram of (3.1), there is an exact sequence of Banach manifolds:

$$\mathcal{A} \xrightarrow{i_1^* \times i_2^*} \mathcal{A}_1 \times \mathcal{A}_2 \xrightarrow{j_1^* \times j_2^*} \mathcal{A}_0 \times \mathcal{A}_0, \tag{3.1}$$

where $i_1^*, i_2^*, j_1^*, j_2^*$ are the natural mappings induced by inclusions $i_1 : Y_1 \to Y, i_2 : Y_2 \to Y; j_1 : Y_0 \to Y_1, j_2 : Y_0 \to Y_2$. The term exact refers to the fact that $i_1^* \times i_2^*$ is an imbedding, $j_1^* \times j_2^*$ is a submersion and the image of $i_1^* \times i_2^*$ equals to the preimage of the diagonal in $\mathcal{A}_0 \times \mathcal{A}_0$ under $j_1^* \times j_2^*$.

Let $\mathcal{G}_j = Aut(P_j)$ be the gauge group of bundle automorphisms of $P_j$, completed with repsect to the appropriate Sobolev norm $L_{k+1}^p$ with $k + 1 > \frac{3}{p}$. Then $\mathcal{G}_j$ operates on $\mathcal{A}_j$ and induces an diagram of maps

$$\mathcal{B} \xrightarrow{i_1^* \times i_2^*} \mathcal{B}_1 \times \mathcal{B}_2 \xrightarrow{j_1^* \times j_2^*} \mathcal{B}_0 \times \mathcal{B}_0, \tag{3.2}$$

on the quotient spaces $\mathcal{B}_j = \mathcal{B}(P_j) = \mathcal{A}_j/\mathcal{G}_j$. However these quotient spaces $\mathcal{B}_0, \mathcal{B}_1, \mathcal{B}_2$ and $\mathcal{B}$ are no longer manifolds but infinite dimensional stratified spaces [10] §4.2.2.

The top strata $\mathcal{B}_j^*$ of $\mathcal{B}_j, j = 0, 1, 2, \emptyset$, are respectively the subspaces of gauge equivalent classes of irreducible connections on $P_j, \mathcal{B}_j^* = \mathcal{A}_j^*/\mathcal{G}_j$. They are Banach manifolds with the tangent spaces $T_{[A]}\mathcal{B}_j^*$ at a point $[A]$ given by

$$\mathcal{T}_{jA} = \{a \in L_k^p(\Omega^1(Y_j, AdP_j)) \mid d_A^* a = 0 \text{ and } i_j^*(*a) = 0\}. \tag{3.3}$$

(c.f. [34] Proposition 4.1) Here $i_j : \partial Y_j \to Y_j$ is the inclusion, $i_j^*(*a)$ is the pullback of the 2-form $*a$ via $i_j$ to the boundary. When $j = \emptyset$ the boundary condition $i_j^*(*a)$ becomes vacuous and so $\mathcal{T}_A$ is defined in the same way as in (2.1).

There are the vector bundles $\mathcal{L}_j$ over $\mathcal{A}_j^*$ (which is not the tangent bundle) whose fiber $\mathcal{L}_{jA_j}$ at a point $A_j$ is isomorphic to the Sobolev space

$$\mathcal{L}_{jA_j} = \{a \in L_{k-1}^p(\Omega^1(Y_j, AdP_j)) \mid d_{A_j}^* a = 0, \quad i_j^*(*a) = 0\}. \tag{3.4}$$

By assigning to $A_j \in \mathcal{A}_j^*$ the 1-form $*F_{A_j} \in \mathcal{L}_{jA_j}$, we obtain a section $f_j : \mathcal{A}_j^* \to \mathcal{L}_j$ of this bundle. Since $\mathcal{L}_j$ and $f_j$ are equivariant with respect to the gauge group action, they descend to bundle $\mathcal{L}_j \to \mathcal{B}_j^*$ and section $f_j : \mathcal{B}_j^* \to \mathcal{L}_j$ over $\mathcal{B}_j^*$. (Following [34], we use the same notation for these bundles and sections on $\mathcal{A}_j^*$ and $\mathcal{B}_j^*$.) The zero set $f_j^{-1}(0)$ is the space of flat irreducible connections on $P_j$. Let $\nabla f_j : \mathcal{T}_j \to \mathcal{L}_j$ denote the covariant derivative of $f_j$. Then, the kernels and cokernels of $\nabla f_{jA_j}$ are cohomologies of $Y_j$ with twisted coefficients

$$\ker \nabla f_{jA_j} = H^1(Y_j, Ad\rho_j), \quad \operatorname{coker} \nabla f_{jA_j} = H^1(Y_j, \partial Y_j; Ad\rho_j), \tag{3.5}$$

where $\rho_j$ is the holonomy representation of $A_j$.

Thus the flat connection $A$ is nondegenerate if and only if these subspaces $H^1(Y_1, Ad\rho_1)$ and $H^1(Y_2, Ad\rho_2)$ intersect each other transversely. We consider only those perturbations $\pi = (\gamma, \phi) \in \Pi$ for which the framed loops $\gamma_j(S^1 \times D^2), j = 1, \cdots, m$, lie inside the submanifold $\Sigma \times (-2, -1) \subset Y_1$ ([13], [34]). We denote by

$$f_j^\pi = f_j + \nabla h : \quad \mathcal{B}_j^* \to \mathcal{L}_j^*, \quad j = 0, 1, 2, \emptyset, \tag{3.6}$$

the new sections after the perturbation. Since $\Sigma \times (-2, -1)$ is outside of $Y_0$ and $Y_2$, the sections $f_0^\pi$ and $f_2^\pi$ are the same as before while $f_1^\pi$ and $f^\pi$ are changed from $f_1$ and $f$ by adding the term $\operatorname{grad} h_\pi$.

**Proposition 3.1.1** *(Taubes) Let $\mathcal{R}_j^\pi = (f_j^\pi)^{-1}(0) \cap \mathcal{B}_j^*$, $j = 0, 1, 2, \emptyset$. For a Baire first category of perturbations in $\Pi$, we have*

1. *The sets $\mathcal{R}_j^\pi$ and $\mathcal{R}^*(Y_j)$ are identical on an open neighborhood of their intersections with $\mathcal{B}_j \setminus \mathcal{B}_j^*$, that is on the complement of compact sets in $\mathcal{R}_j^\pi$ and $\mathcal{R}^*(Y_j)$.*

2. *$\mathcal{R}_j^\pi$ is a smoothly embedded submanifold of $\mathcal{B}_j^*$, which is smoothly isotopic to $\mathcal{R}^*(Y_j)$ by an ambient isotopy which is the identity on an open neighborhood of the intersection of $\mathcal{R}_j^\pi$ with $\mathcal{B}_j \setminus \mathcal{B}_j^*$.*

3. *For given $\varepsilon > 0$, the $L_1^2$-distance moved by a point in $\mathcal{R}_j^\pi$ by the isotopy in (2) is less than $\varepsilon$.*

4. *The intersection of $\mathcal{R}^*(Y_1)$ with $\mathcal{R}^*(Y_2)$ in $\mathcal{R}^*(Y_0)$ is transverse. In fact the transversal intersection is equivalent to that the intersection point is a nondegenerated zero of $f^\pi$.*

The proof is in [34] (see Proposition 1.5 and Proposition 5.1).

If we verify that the perturbations in gauge theory also provid Hamiltonian perturbations in symplectic theory, then once again we have the same generators for two Floer chain groups. Now we define an appropriate set of perturbations of Chern-Simons functional (see [13], [34]).

Let $\gamma = \{\gamma_i\}_1^m$ be a finite collection of disjoint embeddings of solid tori $\gamma_i : S^1 \times D^2 \to Y$. Choose a corresponding collection of functions $h_i' \in C^2([-2,2], \mathbf{R})$ and $h_i = h_i' \circ tr : SU(2) \to \mathbf{R}$. Given the collections $\{\gamma_i\}_1^m$ and $\{h_i\}_1^m$ we define a function $h : \mathcal{A}(Y) \to \mathbf{R}$ by

$$h(A) = \sum_{i=1}^m \int_{D^2} h_i(hol_{\gamma_i}(x, A))\mu, \tag{3.7}$$

where a smoothly compact supported volume form $\mu$ on the interior of $D^2$ satisfies $\int_{D^2} \mu = 1$. The Ad-invariance of the $\{h_i\}$ insures that the resulting function $h$ does not depend on the choice of basepoint in defining the holonomy around a loop $\gamma_i(S^1 \times \{x\})$, and in particular, descends to a map $h_i^j : \mathcal{B}(Y_j) \to \mathbf{R}$. The components of the space of such functions are in bijective correspondence with isotopy classes of links in $Y$.

**Lemma 3.1.2** *(Taubes [34], Fleor [13]) Let $h$ be the above function, referred to as admissible function. The $L^2$-gradient of $h$ is a 1-form $\nabla h$ such that for any tangent vector $a \in T_A\mathcal{A}(Y)$,*

$$Dh(A)(a) = < \nabla h(A), a >_{L^2}, \tag{3.8}$$

*where the zero of $\nabla h$ is outside of $\cup_i \gamma_i(S^1 \times D^2)$, and $\nabla h(A) = \sum_i \nabla_i h * d^2\theta$.*

Let $\omega$ be the symplectic structure on $\mathcal{R}^*(Y_0)$ and $\overline{\omega} : T\mathcal{R}^*(Y_0) \to T^*\mathcal{R}^*(Y_0)$ be the isomorphism of the tangent and cotangent bundle determined by $\omega$. If $f \in C^\infty(\mathcal{R}^*(Y_0))$ is a smooth function, its derivative $df$ is a 1-form and $\overline{\omega}^{-1}(df)$ is a vector field, denote by $Hf$. We call $Hf$ a Hamiltonian vector field. Our perturbations also give arise to functions

$$h_i^j : Hom(\pi_j, SU(2)) \to \mathbf{R}, \quad i = 1, 2, \cdots, m; \quad j = 0, 1, 2, \emptyset,$$

by $h_i^j(hol_{\gamma_i}(\cdot, A))$, where $\gamma_i \in \pi_1(Y_j) = \pi_j, tr : SU(2) \to \mathbf{R}$. Since $tr$ is an invariant function, $h_i^j$ is $SU(2)$-invariant and defines a function on $Hom(\pi_j, SU(2))/SU(2)$, also denoted $h_i^j$.

For $\phi \in \mathcal{R}(Y_j), u \in Z^1(\pi_j, ad\phi)$ is a cocycle representing $[u] \in H^1(\pi_j, ad\phi)$, we have

$$dh_i^j(\phi)([u]) = \frac{d}{dt}|_{t=0} h_i^j(\phi_t(u)) = B(F_i^j(\phi(\alpha)), u(\alpha)), \tag{3.9}$$

where $\phi_t(u(\alpha)) = \exp(tu(\alpha) + o(t^2))\phi(\alpha)$, $B$ is the Killing form in which the symplectic form $\omega(a,b) = B_*(a,b) \cap [\mathcal{R}^*(Y_0)]$ ($[\mathcal{R}^*(Y_0)]$ fundamental class) and $F_i^j$ is the differential of $h_i^j$ with repect to $B$. Since $B$ is nondegenerate, it defines an isomorphism $\overline{B} : \mathbf{su(2)} \to \mathbf{su(2)}^*$ ($\mathbf{su(2)}$ is the Lie algebra of $SU(2)$), and induces isomorphism on the first group cohomology of $\pi_j$. Let $\overline{B}^t : H^1(\pi_j, \mathbf{su(2)}^*)^* \to H^1(\pi_j, \mathbf{su(2)})$ be the transpose of $\overline{B}$. So we have

$$(\overline{B}^t)^{-1}(dh_i^j(\phi)) = \eta_j(\alpha \otimes F_i^j(\phi(\alpha))),$$

where $\eta_j : H_1(\pi_j, \mathbf{su(2)}) \to H_1(\pi_j, \mathbf{su(2)}^*)^*$ is arising from the cap product. Let $\theta_j$ be the map arising from $\omega$:

$$\theta_j : H^1(\pi_j, \mathbf{su(2)}) \to H^1(\pi_j, \mathbf{su(2)}^*)^*.$$

Then we have the following commutative diagram from Lemma 3.8 in [21]:

$$\begin{array}{ccc} H^1(\pi_j, \mathbf{su(2)}) & \stackrel{\cap [\pi_j]}{\to} & H_1(\pi_j, \mathbf{su(2)}) \\ \overline{\omega} \downarrow & \searrow \theta_j & \downarrow \eta_j \\ H^1(\pi_j, \mathbf{su(2)})^* & \stackrel{\overline{B}^t}{\leftarrow} & H^1(\pi_j, \mathbf{su(2)}^*)^*. \end{array}$$

Here $[\pi_j]$ denotes the fundamental class in $H^2(\pi_j, \mathbf{R})$. Therefore

$$\begin{aligned} \alpha \otimes F_i^j(\phi(\alpha)) &= \eta_j^{-1}(\overline{B})^{-1}(dh_i^j(\phi)) \\ &= (\overline{\omega})^{-1}(dh_i^j(\phi)) \cap [\pi_j] \\ &= H(h_i^j(\phi)) \cap [\pi_j]. \end{aligned}$$

Thus the Poincaré duality isomorphism takes the homology class of the cycle $\alpha \otimes F_i^j(\phi(\alpha))$ to the Hamiltonian vector field $H(h_i^j(\phi))$ ( see Proposition 3.7 in [21]). Now we take that $B$ is the canonical bilinear form on $\mathbf{su(2)}$, thus $F_i^j = \nabla_i h^j$ the partial derivative of the lifting of $h$ to $SU(2)^m$ in the direction of $i$th factor. Hence we use the theorems in [21] to show that $h$ also provides an Hamiltonian vector field perturbation.

**Proposition 3.1.3** *(1) Any admissible perturbation in Floer instanton theory also gives a Hamiltonian vector field perturbation in Floer symplectic theory for Heegaard splitting of a 3-manifold.*

*(2) As set, $C_*(Y) = C_*^{sym}(\mathcal{R}^*(Y_1), \mathcal{R}^*(Y_2); \mathcal{R}(Y_0)^*)$, both generators for Floer instanton chain groups and Floer symplectic chain groups are exactly same.*

### 3.2  Maslov index and spectral flow

The symplectic action $\alpha : \Omega(L_1, L_2) \to \mathbf{R}$ is given by

$$< d\alpha(u), \xi > = \int_0^1 \omega(\frac{du}{ds}, \xi) ds, \qquad (3.10)$$

where $\xi$ is a smooth section of the induced bundle $u^*TP$ on $[0,1]$. The Hessian of the symplectic functional $\alpha$:

$$E_u : T_u \mathcal{P}^p_{k,\delta}(a,b) \to L^p_{k-1,\delta}(u^*TP),$$

with the following properties:

$$u(0,t) \in L_1, \quad u(1,t) \in L_2, \quad \text{for any } t \in \mathbf{R},$$

$$u(s,-\infty) = a, \quad u(s,+\infty) = b, \quad \text{for all } s \in [0,1],$$

where $a, b \in L_1 \cap L_2$. We would like to associate two indices, so called spectral flow and Maslov index to the family of operators $\{E_{u(\cdot,t)}\}$ both of which depend only on $a,b \in L_1 \cap L_2$ and the homotopy class of the map $u$.

There are some problems in defining spectral flow and Maslov index. (1) If $u \in \Omega(L_1,L_2)$ is a constant path, then $E_u$ is formally symmetric. If $P$ is Kähler and $g$ is a Kähler metric, then $E_u$ is symmetric for any path $u \in \Omega(L_1,L_2)$. But the operator $E_u$ is not even formally symmetric in general; (2) Each operator $E_{u(\cdot,t)}$ acts on a different space of sections. One must make the domains of the operators to be identical one, and give a definite meaning of the continuity of the family $\{E_{u(\cdot,t)}\}$ with respect to the parameter $t$; (3) $\{E_{u(\cdot,t)}\}$ does not form a loop, if $a \neq b$. One may make it a loop in anatural way. Otsuki and Furutani have made this clear and identify the spectral flow with the Maslov index in [31]. $SF(E_{u(\cdot,t)}) = Mas\{(L_1(t), L_2(t))\}$ has been proved in [17], [31]. Thus we will use Cappell, Lee and Miller's theorem [7] and Yoshida's [36] to identify the spectral flow from instanton with the spectral flow from holomorphic curve.

Let $(Y, Y_1, Y_2, Y_0)$ be a Heegaard decomposition of $Y$, $E$ a vector bundle over $Y$, and $\{D(t) | 0 \leq t \leq 1\}$ a family of first order, self-adjoint, elliptic operator $D(t) : \Gamma(E) \to \Gamma(E)$ on the space of smooth sections of $E$. Furthermore on $Y_0$, the operator has the following form

$$D(t)|_{Y_0} = \sigma \circ (\frac{\partial}{\partial s} + \hat{D}(t)), \tag{3.11}$$

with $\hat{D}(t)$ self-adjoint, elliptic operator on $\Gamma(E|_{Y_0})$ and $\sigma$ a bundle automorphism. In an obvious way, we replace $Y$ by $Y(s) = Y_1 \cup (\Sigma \times [-s,s]) \cup Y_2$, the bundle $E$ by $E(s)$ which is the same as $E$ over $Y_1, Y_2$ and $E(s)|_{\Sigma \times [-s,s]} = E|_{\Sigma} \times [-s,s]$, and the operator $D(t)$ by $D(t,s) : \Gamma(E(s)) \to \Gamma(E(s))$ which is given by (3.11) on $\Sigma \times [-s,s]$ and extends to $Y_1, Y_2$ as $D(t)$. Note that $s$ is our notation for the variable along the surface tube.

Assume that there exist $k \geq 0$ and $\delta > 0$ such that $\hat{D}(t)$ has no eigenvalues in the intervals $(k, k+\delta)$, $(-k-\delta, -k)$. Let $H(t,k)$ be the symplectic vector space spanned by eigensections $\phi_j$ of $\hat{D}(t), \hat{D}(t)\phi_j = \lambda_j \phi_j$ with $|\lambda_j| \leq k$, and let $P_\pm(t,k)$ be the $L^2$-closure of the spanned eigensections $\phi_j$ with eigenvalues $|\lambda_j| > k$. By the spectral decomposition theorem $L^2(E|_\Sigma) = P_+(t,k) \oplus H(t,k) \oplus P_-(t,k)$, and as $t$ varies these spaces $P_\pm(t,k), H(t,k)$ vary smoothly. In the situation $k = 0$, we denote the spaces by $P_\pm(t)$ and $H(t)$.

Let $\mathcal{L}_1(t)$ and $\mathcal{L}_2(t), a \leq t \leq b$ be a choice of smoothly varying Lagrangian pairs in $H(t,k)$ such that they satisfy the endpoint condition:

$$\mathcal{L}_1(t) = L_1(t) \oplus [P_+(t) \cap H(t,k)], \quad \text{if} \quad t = a, b,$$

$$\mathcal{L}_2(t) = L_2(t) \oplus [P_-(t) \cap H(t,k)], \quad \text{if} \quad t = a, b,$$

where $L_1(t), L_2(t)$ are the subspaces in $H(t)$ consisting of extended $L^2$-solutions of $D(t)|_{Y_1}$ and $D(t)|_{Y_2}$ respectively. With respect to these choices of Lagrangians $\mathcal{L}_1(t), \mathcal{L}_2(t)$ there are the following self-adjoint, elliptic operators with global boundary conditions:

$$D_1(t, \mathcal{L}_1(t)) : L^2_1(E(s)|_{Y_1(s)}; \mathcal{L}_1(t) \oplus P_+(t,k)) \to L^2(E(s)|_{Y_1(s)}),$$

$$D_2(t, \mathcal{L}_2(t)) : L^2_1(E(s)|_{Y_2(s)}; \mathcal{L}_2(t) \oplus P_-(t,k)) \to L^2(E(s)|_{Y_2(s)}),$$

and hence well-defined $(1/s^2)$-spectral flow of $[D_j(t, \mathcal{L}_j(t)); a \leq t \leq b]$ ( for $\delta$-spectral flow definition see [7], [9] and [13]) on $Y_1(s), Y_2(s)$.

**Theorem 3.2.1** *(The spectral flow decomposition theorem) With the notation as above, we have for $s$ sufficiently large*

$$(\frac{1}{s^2}) - SF\{D(s,t)|_{Y(s)} | a \leq t \leq b\} = \sum_{j=1}^{2} (\frac{1}{s^2}) - SF\{D_j(t, \mathcal{L}_j(t))|_{Y_j(s)} | a \leq t \leq b\}$$

$$+ Mas\{(\mathcal{L}_1(t), \mathcal{L}_2(t)) | a \leq t \leq b\} + \frac{1}{2}[dim \ker \hat{D}(b) - dim \ker \hat{D}(a)], \quad (3.12)$$

*where the second term stands for the Maslov index of the Lagrangian path $(\mathcal{L}_1(t), \mathcal{L}_2(t))$. (see [7])*

Under certain hypothesis, Yoshida first gave such a decomposition formula for spectral flow in [36]. Motived by the computation of Casson invarinat in [6], Cappell, Lee and Miller generalized the decomposition formula of spectral flow for general 3-manifolds. There are active recent studies on this kind of problem. We refer the reader to [7].

**Proposition 3.2.2** *Let $D(s,t)$ be the self-duality operator defined in §2.1 (2.2). Then*

$$SF(D(s,t) : 0 \leq t \leq 1) = Mas\{(L_1(t), L_2(t)) : 0 \leq t \leq 1\}. \qquad (3.13)$$

Proof: The identiy follows by applying Theorem 3.2.1. From the transverse intersections at the end $t = 0, 1$,

$$\dim \ker \hat{D}(1) - \dim \ker \hat{D}(0) = 0.$$

Also the $L^2$-solution space from $Y_j$ is unchanged and $\ker D(t, L_j(t)) = L_j(t)$ gives of constant dimension, i.e

$$(\frac{1}{s^2}) - SF\{D_j(t, \mathcal{L}_j(t))|_{Y_j(s)} | a \leq t \leq b\} = 0.$$

Hence from Theorem 3.2.1, we have identified the $(1/s^2)$-spectral flow of $D$ with the Maslov index. Because of irreduciblity there is no zero mode, so there is no descrepency between the $(1/s^2)$-spectral flow and the usual spectral flow. Therefore, the identity holds. ∎

**Remark:** For $J$-holomorphic curve $u$ over $\Sigma \times [-s, s] \times \mathbf{R}$, if $u$ connecting two transversal Lagrangian intersections restricts to a connection $A(t, s)$ on $\Sigma \times [-s, s]$ without $ds, dt$ forms for $A(t, s) \in \mathcal{R}^*(Y_0)$, then Proposition 3.2.2 just follows from Theorem 1.1 in [36]. That is why we use the more general spectral flow decomposition theorem from [7].

¿From [15] and [31], the spectral flow of $E_u$ operator is the same as the Maslov index from Lagrangian intersections. Hence combining the above Proposition 3.2.2, we have that

$$SF(D(t, s) : 0 \leq t \leq 1) = SF((E_{u(t, \circ)} : 0 \leq t \leq 1)).$$

Note that this identity gives us a correspondence between the instanton and symplectic chain complex. I.e.

$$C_j(Y) \equiv C_j(\mathcal{R}_1^*, \mathcal{R}_2^*; \mathcal{R}^*(Y_0)). \qquad (3.14)$$

For a fixed element $x$ of $\mathcal{R}^*(Y)$, we have

$$C_j(Y) = \{a \in \mathcal{R}^*(Y) : SF(D(a, x)) = j \pmod{8}\},$$

$$C_j(\mathcal{R}_1^*, \mathcal{R}_2^*; \mathcal{R}^*(Y_0)) = \{a \in \mathcal{R}_1^* \cap \mathcal{R}_2^* : SF(E_u(a, x)) = SF(D(a, x))\}.$$

Now the chain complexes for instantons and Lagrangian intersections are identically same.

# 4 Smallest eigenvalues

## 4.1 Smallest eigenvalue for self-duality operator

Let $Y$ be a homology 3-sphere. For $\delta \geq 0$ (to be determined), let $e_\delta : Y \times \mathbf{R} \to \mathbf{R}$ be a smooth positive function with $e_\delta(y,t) = e^{\delta|t|}$ for $|t| \geq 1$. Let $E$ be an $SU(2)$-vector bundle over $Y \times \mathbf{R}$ with a translationally invariant metric and a metric-preserving connection. To define the Banach manifolds $\mathcal{B}(a,b)$ of paths connecting two flat connections $a$ and $b$ in $\mathcal{B}_Y$, we choose representatives $a, b \in \mathcal{A}_Y$ and a smooth connection $C$ on $Y \times \mathbf{R}$ which coincides with $a$ for $t \leq -1$ and with $b$ for $t \geq 1$. Define the $L^p_{k,\delta}$-norm on sections $u$ of $E$ by

$$\|u\|_{L^p_{k,\delta}} = \|e_\delta \cdot u\|_{L^p_k}. \tag{4.1}$$

Then

$$\mathcal{A}_\delta(a,b) = C + L^p_{1,\delta}(\Omega^1_{ad}(Y \times \mathbf{R}))$$

is an affine space and is independent of the choice of $C$. The corresponding gauge group is the group $\mathcal{G}_\delta$ obtained by completing the compactly supported gauge transformations in the $L^p_{2,\delta}(\Omega^0_{ad}(Y \times \mathbf{R}))$. We need $p > 2$ to construct the orbit space $\mathcal{B}_{Y \times \mathbf{R}} = \mathcal{A}_\delta / \mathcal{G}_\delta$.

**Proposition 4.1.1** *1. Let $D_a : L^p_{1,\delta}(\Omega^1 \oplus \Omega^0)(Y, adSU(2)) \to L^p_{0,\delta}(\Omega^1 \oplus \Omega^0)(Y, adSU(2))$ be the operator $D_a(\alpha, \beta) = (\star d_a \alpha - d_a \beta, -d^\star_a \alpha)$. Then there exists a positive $\lambda_0$ such that for all $a \in \mathcal{R}^*(Y)$ the eigenvalues of $D_a$ satisfy $|\lambda(D_a)| \geq \lambda_0$.*

*2. If $F(A)$ is in $L^p$ for $p \geq 2$, then there is a constant $C_A$ such that*

$$\sup |F_A|_{y,t} \leq C_A e^{-\gamma|t|}, \tag{4.2}$$

*where $\gamma = \gamma(\lambda_0) > 0$, and $C_A$ is continuous in $A \in \mathcal{A}_\delta(a,b)$.*

*3. For an anti-self-dual connection $A$ over $Y \times \mathbf{R}$, there is a gauge transformation $g$ on the bundle $Y \times \mathbf{R} \times SU(2)$ such that $g^*(A) = a_\pm + a$ with $F(a_\pm) = 0$ and*

$$\sup |a|_{y,t} \leq C_A e^{-\gamma|t|} \tag{4.3}$$

*Moreover, we can choose $a$ so that all derivatives decay exponentially:*

$$\sup |\nabla^{(l)} a|_{(y,t)} \leq C_A e^{-\gamma|t|}.$$

**Proof:** The first statement follows from [13], and the second and third are in [9] (see 4.1). ∎

The above estimates remain valid for perturbed ASD connections and general 3-manifolds as in [19]. Fix a positive $\delta_{ins} < \min\{\lambda_0, \frac{\gamma}{2}\}$. We will henceforth use the norm (4.1) with $\delta = \delta_{ins}$. Let us denote $\|u\|_{L^p_{1,\delta}(A)} = \|\nabla_A u\|_{L^p_{0,\delta}} + \|u\|_{L^p_{0,\delta}}$ (and $\|u\|_{L^p_{1,\delta}} = \|u\|_{L^p_{1,\delta}(C)}$).

**Definition 4.1.2** : *(1) The balancing function $b : \mathcal{B}_{Y \times \mathbf{R}} \to \mathbf{R}$ is given by the equation:*

$$\int_{-\infty}^{b(A)} \|F(A)\|^2_{L^2(Y)} = \int_{b(A)}^{\infty} \|F(A)\|^2_{L^2(Y)}.$$

*(2) Set the balanced moduli space $\hat{\mathcal{M}}_{Y \times \mathbf{R}} = \{A \in \mathcal{M}_{Y \times \mathbf{R}} \subset \mathcal{B}_{Y \times \mathbf{R}} | \ b(A) = 0\}$.*

(So the value $b(A)$ is the time which splits the action of $A$ in half.) We have the following lemma which is proved in [13] Proposition 3b.2 and [25] Lemma 3.1.5.

**Lemma 4.1.3** *If $dim\mathcal{M}_{Y \times \mathbf{R}} \leq 1$, then $\hat{\mathcal{M}}_{Y \times \mathbf{R}}$ is compact.*

Let $d_A$ denote the covariant derivative corresponding to the connection $A$ and $d_A^{*\delta} = e_\delta^{-1} d_A^* e_\delta$ be the adjoint of $d_A$ with respect to the $L^2_{0,\delta}$-norm. Floer has proved the following in [13].

**Proposition 4.1.4 (Floer)** (i) *For positive $\delta$, $\mathcal{G}_\delta$ is a Banach Lie group with Lie algebra (which can be identified with) $L^p_{2,\delta}(\Omega^0_{ad}(Y \times \mathbf{R}))$.*

(ii) *The quotient space $\mathcal{B}_\delta(a,b) = \mathcal{A}^*_\delta(a,b)/\mathcal{G}_\delta$ is a smooth Banach manifold with tangent spaces*

$$T_{[A]}\mathcal{B}_\delta(a,b) = \{\alpha \in L^p_{1,\delta}(\Omega^1_{ad}(Y \times \mathbf{R})) \mid d_A^{*\delta}\alpha = 0\}.$$

(iii) *The 2-form $F_A^-$ representing the anti-self-dual part of the curvature of $A$ is smooth and $\mathcal{G}_\delta$-equivariant.*

(iv) *If $\delta > 0$ is smaller than the smallest nonzero absolute value of an eigenvalue of $D_a$ or $D_b$, then for any anti-self-dual connection $A \in \mathcal{B}_\delta(a,b)$ the anti-self-duality operator*

$$D_A^\delta = d_A^{*\delta} \oplus d_A^+ : L^p_{1,\delta}\Omega^1_{ad}(Y \times \mathbf{R}) \to L^p_{0,\delta}(\Omega^0_{ad} \oplus \Omega^2_{ad,+})(Y \times \mathbf{R})$$

*is Fredholm. Furthermore, $D_A^\delta = \frac{\partial}{\partial t} + D_{a_t}^\delta$, where*

$$D_{a_t}^\delta = \begin{pmatrix} \star d_a & -d_a \\ -d_a^\star & \delta \end{pmatrix}.$$

$D^\delta_{a_t}$ is self-adjoint on $\Omega^1_{ad}(Y) \oplus \Omega^0_{ad}(Y)$ and $\star$ is the Hodge operator on the 3-manifold $Y$. If $a$ and $b$ are irreducible nondegenerate flat connections, then one can take $\delta = 0$.

**(v)** Let $\mathcal{M}$ be the moduli space of all equivalence classes of nonflat anti-self-dual connections $A$ on $Y \times \mathbf{R}$ with finite action $||\frac{\partial A}{\partial t}||^2_2$. There is a first category set of metrics on $Y$ such that the anti-self-duality operator $D^\delta_A$ is surjective for all $A \in \mathcal{M} \cap \mathcal{B}_\delta$.

**Theorem 4.1.5** *Suppose $\dim \mathcal{M}_{Y \times \mathbf{R}} = 1$. Then there exists a positive constant $C_p$ such that for all $A \in \hat{\mathcal{M}}_{Y \times \mathbf{R}}$, and for all $p \geq 2, u \in L^p_{0,\delta}(\Omega^2_{ad,+}(Y \times \mathbf{R}))$ we have*

$$C_p \cdot \int_{Y \times \mathbf{R}} e^{p\delta \cdot |t|} |u|^p \leq \int_{Y \times \mathbf{R}} e^{p\delta \cdot |t|} |(d^+_A)^{*\delta} u|^p.$$

This is Theorem 3.2.6 in [25].

## 4.2 Smallest eigenvalue for Cauchy-Riemann operator

Let $P = \mathcal{R}^*(Y_0)$ be the symplectic manifold with symplectic form $\omega$. Given two Lagrangian submanifolds $L_j = \mathcal{R}^*(Y_j), j = 1, 2$, of $P$, we define

$$\mathcal{P}^p_{k,\text{loc}} = \{u \in L^p_{k,\text{loc}}([0,1] \times \mathbf{R}, P) |\ u(0, \mathbf{R}) \subset L_1, u(1, \mathbf{R}) \subset L_2\}, \tag{4.4}$$

for $k \geq 2/p$. We will identify $[0,1] \times \mathbf{R}$ with $iI \times \mathbf{R} \subset \mathbf{C}$. Let $S_\omega$ be the bundle over $P$ whose fibre is given by

$$S_x = \{J \in End(T_x P) |\ J^2 = -id \text{ and } \omega(\cdot, J\cdot) \text{ is a metric}\}.$$

Denote by $\mathcal{J} = C^\infty([0,1] \times S_\omega)$ the set of almost complex structures on $[0,1]$. We will always use $J$ to be an element of $\mathcal{J}$ unless otherwise stated. Define

$$\overline{\partial}_J u(s,t) = \frac{\partial u(s,t)}{\partial t} + J_s \frac{\partial u(s,t)}{\partial s}$$

on $\mathcal{P}^p_{k,\text{loc}}$ and then the solution of $\overline{\partial}_J u(s,t) = 0$ is translationally invariant in $t$.

For $\delta \geq 0$, define the Banach manifold of paths connecting two transversal intersection points $a, b \in L_1 \cap L_2$,

$$\mathcal{P}^p_{k,\delta}(a,b) = \{u \in \mathcal{P}^p_{k,\text{loc}}\ |\ \|e_\delta \cdot u\|_{L^p_k([0,1] \times \mathbf{R})} < \infty\},$$

where $e_\delta$ is defined in the same way as in (4.1). If $E$ is a smooth vector bundle over $P$ and $u \in \mathcal{P}^p_{k,\text{loc}}$, then the pullback bundle $u^*E$ has the structure of a $L^p_k$-bundle, i.e. of a locally trivial bundle with transition maps in $L^p_k(U, GL(k, \mathbf{R}))$ for $U \subset [0,1] \times \mathbf{R}$. So

$$L^p_{k,\delta}(u^*E) = \{\xi \in \mathcal{P}^p_{k,\text{loc}}(u^*E) \mid \|\xi\|_{L^p_{k,\delta}} < \infty\}. \tag{4.5}$$

We will take $E$ to be the tangent bundle of $P$. Then for $k = 1, p > 2$, the set $\mathcal{P}^p_{1,\delta}(a,b)$ for any $a, b \in L_1 \cap L_2$ is a smooth Banach manifold with tangent spaces

$$T_u \mathcal{P}^p_{k,\delta} = \{\xi \in L^p_{k,\delta}(u^*TP) \mid \xi(0,t) \in T_{u(0,t)}L_1, \; \xi(1,t) \in T_{u(1,t)}L_2 \quad \text{for all } t \in \mathbf{R}\}$$

(This is Theorem 3 in [14]).

**Proposition 4.2.1** *(1) Let $A_a = J\frac{d}{ds} : L^p_{1,\delta}(u^*T_aP) \to L^p_{0,\delta}(u^*T_aP)$ be the first order selfadjoint Dirac operator. Then there exists a constant $\lambda_{sym} > 0$ such that for all $a \in L_1 \cap L_2$, the eigenvalues of $A_a$ satisfy $|\lambda(A_a)| \geq \lambda_{sym}$.*

*(2) If $\nabla u$ is in $L^p$ for $p \geq 2$, then there exists a constant $C_u$ such that*

$$\int_{[0,1]} |u(s,t)|^2 ds \leq C_u e^{-\gamma_1 |t|}, \tag{4.6}$$

*where $\gamma_1 = \gamma_1(\lambda_{sym}) > 0$ and $C_u$ is continuous in $u$.*

Proof: (1) follows from Lemma 4.3 in [14], and the second follows from the proof of Theorem 4 in [14]. ∎

Fix a positive $\delta_{sym} < \min\{\lambda_{\text{sym}}, \frac{\gamma_1}{2}\}$. ¿From now on, we will simply use the $\delta$-weight for both instanton and symplectic theories with $\delta \leq \min\{\delta_{ins}, \delta_{sym}\}$. Similarly we define the balancing function in symplectic setting.

**Definition 4.2.2** *(1) The function $b_{sym} : \mathcal{P}^p_{k,\delta}(a,b) \to \mathbf{R}$ is given by*

$$\int_{-\infty}^{b_{sym}(u)} \|\nabla u\|^2_{L^2([0,1])} = \int_{b_{sym}(u)}^{+\infty} \|\nabla u\|^2_{L^2([0,1])}.$$

*(2) Set the balanced moduli space*

$$\hat{\mathcal{M}}_J(a,b) = \{u \in \mathcal{M}^s_J(a,b) \subset \mathcal{P}^p_{k,\delta}(a,b) \mid b_{sym}(u) = 0\}. \tag{4.7}$$

Note that our definition is same as Floer in [14] (3.2) and [15] (2.14).

**Lemma 4.2.3** *Suppose that $L_1$ intersects with $L_2$ transversally. Then there is a subset $\mathcal{J}_1(L_1, L_2)$ of $\mathcal{J}$ which is dense and open in $\mathcal{J}$, for which the balanced moduli space $\hat{\mathcal{M}}_J(a,b)$ is compact if $dim \mathcal{M}_J(a,b) \leq 1$. (c.f. [29] for construction of $\mathcal{J}_1(L_1, L_2)$)*

Proof: This is Oh's Proposition 4.1 in [29]. The bubbling off a sphere can be easily ruled out by the dimension counting (or index computation); the bubbling off a disk can be ruled out by the analysis in [30]. By Proposition 2.2.14, a holomorphic curve in an one dimensional moduli space can not touch any point in $\mathcal{S}^*(Y_j)$. So Oh's argument can be applied for $L_j = \mathcal{R}_j^*$. ∎

For every $u \in \mathcal{P}_{k,\delta}^p(a,b)$, we can define the linearized operator of $\overline{\partial}_J$:

$$E_u = D\overline{\partial}_J(u) : T_u \mathcal{P}_{k,\delta}^p(a,b) \to L_{k-1,\delta}^p(u^*TP),$$

which is given by

$$E_u(\xi)(s,t) = (\nabla_t + J\nabla_s)(\xi)(s,t) + B(s,t)(\xi)(s,t) \tag{4.8}$$

with matrix operator $B$ only depending on the choice of the connection $\nabla$. As $t \to \pm\infty$, the operator $E_u$ approaches translationally invariant operators of the form

$$E_{\pm\infty}\xi = \frac{\partial}{\partial t} + A_{u(\pm\infty)},$$

where $A_{u(\pm\infty)}$ is independent of $t$. Therefore the Fredholm property of $E_u$ for $u \in \mathcal{M}_J^s(a,b)$, $a,b$ transversal intersections in $L_1 \cap L_2$, will follow from the asymptotically constant elliptic operators $A_{u(\pm\infty)}$.

**Definition 4.2.4** $u \in \mathcal{M}_{k,\delta}^p(a,b)$ *is regular if $Coker E_u = 0$; $(J, L_1, L_2)$ is regular if $u$ is regular for all $u \in \mathcal{M}_{k,\delta}^p(a,b)$ and $a,b \in L_1 \cap L_2$.*

**Proposition 4.2.5 (Floer)** *(1) For any path $u \in \mathcal{P}_{k,\delta}^p(a,b), \delta \geq 0$, the operator $E_u$ extends to bounded linear operators*

$$E_u : T_u \mathcal{P}_{l,\delta_1}^q(a,b) \to L_{l-1,\delta_1}^q(u^*TP),$$

*for all $1 \leq l \leq k, 1 \leq q, l - 2/q \leq k - 2/p, \delta_1 \leq \delta$.*

*(2) $E_u$ is Fredholm if and only if $\delta$ does not lie in the set $\sigma(A_a) \cup \sigma(A_b)$ where $\sigma(A_a)$ denotes the spectrum ( eigenvalues) of $A_a$. I.e. for transversal interesections $a, b \in L_1 \cap L_2$, the section $\overline{\partial}_J$ is a Fredholm section of the bundle $T\mathcal{P}^p_{k-1,\delta}(a,b)$ over $\mathcal{P}^p_{k,\delta}(a,b)$.*

*(3) The zero set of $\overline{\partial}_J$ with topology induced by the Banach manifold topology of $\mathcal{P}^p_{1,\delta}(a,b)$ is homeomorphic to $\mathcal{M}^p_{1,\delta}(a,b)$ with the topology of local convergence. The index of the operator $E_u$ $Ind E_u = dim \mathcal{M}^p_{1,\delta}(a,b)$ is the Maslov-Viterbo index $\mu_u(a,b)$.*

*(4) There exists a dense and open set $\mathcal{J}_{reg}(L_1, L_2) \subset \mathcal{J}$ such that $(J, L_1, L_2)$ is regular for $J \in \mathcal{J}_{reg}(L_1, L_2)$, i.e. $E_u$ is surjective for all $u \in \mathcal{M}^p_{1,\delta}(a,b)$.*

Proof: (1) is Lemma 4.1 in [14] with weighted Sobolev embedding theorem in [26]. (2) and (3) come from Theorem 4 in [14]. (4) combines Theorem 5 [14] and Proposition 3.2 in [29]. The proof is contained in Appendix of [29]. ∎

**Theorem 4.2.6** *Suppose $dim \mathcal{M}_J(a,b) = 1$, then there exists a positive constant $C_p$ independent of $u \in \mathcal{M}_J(a,b)$ such that for all $u \in \mathcal{M}_J(a,b)$ and $p \geq 2$, we have*

$$C_p \|\xi\|^p_{L^p_{0,\delta}} \leq \|E_u^* \xi\|^p_{L^p_{0,\delta}},$$

*for $\xi \in L^p_{k-1,\delta}(u^*TP)$ and $E_u^*$ is the $L^2_{0,\delta}$-adjoint of $E_u$.*

Proof: Proposition 4.2.5 (4) implies that $E_u^*$ has trivial kernel, thus we have the following inequality

$$C_{p,u} \|\xi\|^p_{L^p_{0,\delta}} \leq \|E_u^* \xi\|^p_{L^p_{0,\delta}},$$

for $\xi \in L^p_{k-1,\delta}(u^*TP)$. The constant $C_{p,u}$ is continuous in $u$. Hence the result follows from Lemma 4.2.3. ∎

# 5 Comparing the linearized operators

## 5.1 Curvature estimates

Let $P = \mathcal{R}^*(Y_0)$ be the space of irreducible flat $SU(2)$-connections over surface $\Sigma$. Given $\mathcal{R}^*(Y_j), j = 1, 2$ two Lagrangian submanifolds, we have

$$\mathcal{P}^p_{k,\text{loc}} = \{u \in L^p_{k,\text{loc}}([-1,1] \times \mathbf{R}, \mathcal{R}^*(Y_0)) |\ u(-3+2j, \mathbf{R}) \subset L_j, j = 1, 2\},$$

where $u$ is a 1-parameter family of flat connections on $Y_0$ which smoothly extends to two handlebodies via the boundary condition. So for $u \in \mathcal{P}^p_{k,\delta}(a,b)$ we have a corresponding element $u$ in

$$\nabla + L^p_{k,\delta}(\Omega^1(\Sigma \times [-1,1] \times \mathbf{R})),$$

where $\nabla$ is a smooth connection which coincides with $a$ for $t \leq -1$ and for $b$ for $t \geq 1$. Now the holomorphic curve can be viewed as

$$u = A + \phi ds + \psi dt, \quad s \in [-1,1], t \in \mathbf{R}, \tag{5.1}$$

with $A \in \Omega^1(\Sigma), \phi \in \Omega^0(\Sigma), \psi \in \Omega^0(\Sigma)$ and satisfies the following.

1. Flatness: The curvature of $u$ can be computed as

$$F(u) = F_A + (-\frac{\partial A}{\partial s} + d_A\phi) \wedge ds + (-\frac{\partial A}{\partial t} + d_A\psi) \wedge dt - (\frac{\partial \phi}{\partial t} - \frac{\partial \psi}{\partial s} - [\phi,\psi]) \wedge ds \wedge dt. \tag{5.2}$$

   The curvature of $A$ is flat at each $(s,t) \in [-1,1] \times \mathbf{R}$, i.e. $F_A = 0$.

2. Asymptotic flat: $u(s, \pm\infty) = u^\pm (u^+ = b, u^- = a)$ are flat connections on $Y$ since $u^\pm \in L_1 \cap L_2 = \mathcal{R}^*(Y)$. Hence $u^\pm = \alpha^\pm + \Phi^\pm ds$ with properties:

$$F(\alpha^\pm) = 0, \quad \frac{\partial \alpha^\pm}{\partial s} - d_{\alpha^\pm}\Phi^\pm = 0. \tag{5.3}$$

3. Holomorphicity: Let $J$ be an almost complex structure in $End(H^1(Y_0, adSU(2)))$ consisting of elements $J_s \in End(H^1(\Sigma, ad\rho))$ which takes the form $J_s = *_s$ the Hodge star operator of the surface $\Sigma \times \{s\}$, $J_s^2 = -id, \omega(\cdot, J_s\cdot) > 0$. With respect to the almost complex structure $*_s$ and $u \in \mathcal{M}^s_J(a,b)$, the holomorphic curve equation is

$$(\frac{\partial A}{\partial t} - d_A\psi) + *_s(\frac{\partial A}{\partial s} - d_A\phi) = 0, \tag{5.4}$$

in $H^1(Y_0, adSU(2))$ where $\phi, \psi$ are uniquely determined by

$$d^{*_s}_A(\frac{\partial A}{\partial t} - d_A\psi) = 0, \quad d^{*_s}_A(\frac{\partial A}{\partial s} - d_A\phi) = 0. \tag{5.5}$$

4. Symplectic action is given by

$$E(u) = \frac{1}{2} \int_R \int_{[-1,1]} (\|\frac{\partial A}{\partial t} - d_A\psi\|^2_{L^2(\Sigma \times \{s\})} + \|\frac{\partial A}{\partial s} - d_A\phi\|^2_{L^2(\Sigma \times \{s\})}) ds dt. \tag{5.6}$$

By (5.4) above, we have $E(u) = \int_R \int_{[-\frac{1}{2},\frac{1}{2}]} \|\frac{\partial A}{\partial t} - d_A\psi\|^2_{L^2(\Sigma \times \{s\})} ds dt$.

**Lemma 5.1.1** *For $u \in \mathcal{M}_{1,\delta}^p(a,b)$ the self-dual curvature of $u$ is*

$$F^+(u) = -\frac{1+*}{2}(\frac{\partial \phi}{\partial t} - \frac{\partial \psi}{\partial s} - [\phi, \psi]) \wedge ds \wedge dt. \qquad (5.7)$$

*There exists a constant $C > 0$ such that*

$$sup|F(u)|_{y,t} \leq Ce^{-\gamma|t|},$$

*with the constant independent of $u$ if $dim\mathcal{M}_{1,\delta}^p(a,b) = 1$.*

Proof: The (5.7) follows from (5.2) and (5.4) above. Once we view $u$ as a connection on $Y \times \mathbf{R}$ which approaches to flat connections $a, b$, then the Yang-Mills functional gives the finite energy bound:

$$\frac{1}{8\pi^2} \int_{Y \times \mathbf{R}} |F(u)|^2 = cs(b) - cs(a) < \infty.$$

Hence the result follows from Proposition 4.1.1 and Lemma 4.2.3. ∎

¿From Lemma 5.1.1, we have the curvature $F(u) \in L_{0,\delta}^p(\Omega^2(Y \times \mathbf{R}))$ (see the definition of $\delta$ in §4.2). In order to deform the holomorphic curve $u$ into an anti-self-dual connection, one needs to change the metric on $\Sigma$ by the familiar process of the shrinking neck $\Sigma \times [-1, 1]$. The key property is that the new metric should agree with old metric on a large open set of $Y$. Over the " neck " we define $g_\varepsilon$ as follows:

$g_\varepsilon = g$ on the complement of $\Sigma \times [-1 - \varepsilon, 1 + \varepsilon]$;

$g_\varepsilon = ds^2 + \chi_\pm g_\Sigma$ on $\Sigma \times ([-1-\varepsilon, -1] \cup [1, 1+\varepsilon], |d\chi_\pm| \leq 2/\varepsilon$;

$g_\varepsilon = ds^2 + \varepsilon^2 g_\Sigma$ on $\Sigma \times [-1, 1]$.

The shrinking process is to let $\varepsilon \to 0$ which will cause the curvature of the Levi-Civita connection of the manifold $\Sigma \times \{s\}$ to blow up. Our concern is with the estimate for the curvature of $u$ which lies on the $SU(2)$-principle bundle, not the tangent bundle of $Y \times \mathbf{R}$. The volume $Vol(\Sigma \times [-1, 1], g_\varepsilon) = C\varepsilon^2$, where $C$ is a constant.

**Lemma 5.1.2** *For any $\delta_7' > 0$, there exist $\varepsilon_0$ and $C_3$ (independent of $u$) such that for $0 < \varepsilon < \varepsilon_0$, any $u \in \mathcal{M}_J(a,b)$ with $dim\mathcal{M}_J(a,b) \leq 1$ and $p \geq 2$,*

$$\|F(u)^+\|_{L_{0,\delta}^p(Y \times \mathbf{R}, dt^2 + g_\varepsilon)} \leq \delta_7'.$$

Proof: Note that the support of $F(u)^+$ lies on $(\Sigma \times [-1,1]) \times \mathbf{R}$. The metric $g_\varepsilon$ changes uniformly along $s$-direction and so

$$\|F(u)^+\|_{L^p_{0,\delta}(Y \times \mathbf{R}, dt^2 + g^\varepsilon)} \leq (\int_R \int_{\Sigma \times [-1,1]} |e^{\delta|t|} F(u)|^p vol_{g^\varepsilon})^{1/p} \leq C_{\delta,\gamma} \varepsilon^{2/p},$$

by Lemma 5.1.1. Thus we obtain the result. ■

## 5.2 Comparing the anti-self-dual operator and Cauchy-Riemann operator

As in §4.2, the linear operator

$$E_u = D\overline{\partial}_J : T_u \mathcal{P}^p_{k,\delta}(a,b) \to L^p_{k-1,\delta}(u^* TP),$$

is given by

$$E_u \alpha = \pi_A (\nabla_t \alpha + *_s \nabla_s \alpha - *_s dX_s(A)\alpha).$$

Here $\nabla_s = \frac{\partial}{\partial s} + \phi, \nabla_t = \frac{\partial}{\partial t} + \psi$, and $\pi_A : \Omega^1(\Sigma \times \{s\}, ad) \to H^1_A(\Sigma \times \{s\}, ad)$ is the orthogonal projection.

For an anti-self-dual connection $A + \Phi ds + \Psi dt \in \Omega^1 \oplus \Omega^0 \oplus \Omega^0(\Sigma, ad)$, its linearized operator $d^+_{A+\Phi ds + \Psi dt}$ is

$$d^+_{A+\Phi ds+\Psi dt}(a + \phi ds + \psi dt) = \begin{pmatrix} \nabla_t + *_s \nabla_s - *_s dX_s(A) & *_s d_A & -d_A \\ {}_s d_A & \nabla_t & -\nabla_s \\ 0 & 0 & 0 \end{pmatrix} \begin{pmatrix} a \\ \phi \\ \psi \end{pmatrix}, \quad (5.8)$$

with the gauge fixing condition

$$d^*_{A+\Phi ds+\Psi dt}(a + \phi ds + \psi dt) = d^*_A a + *_s \nabla_s *_s \phi + \nabla_t \psi = 0. \quad (5.9)$$

In this subsection, the estimates are essentially due to Dostoglou and Salamon in [11]. The only difference is that we keep these estimates in the weighted Sobolev norm rather than the usual Sobolev norm with conformal factors. This is because they consider nontrivial $SO(3)$ bundle with $w_2 \neq 0$ and in this case $\delta$ can be taken to be zero. Now denote the self-duality operator $D = d^+_{A+\Phi ds+\Psi dt} \oplus d^*_{A+\Phi ds+\Psi dt}$, we need to get the smallest eigenvalue for $D$ from the operator $E_u$.

**Lemma 5.2.1** *There exist $\varepsilon_0, c > 0$ such that for $\xi \in \Omega^1 \oplus \Omega^0 \oplus \Omega^0(\Sigma, ad)$ and $0 < \varepsilon < \varepsilon_0$, we have the following*

$$\|\xi\|_{L^p_{1,\delta}((\Sigma\times[-1,1],g_\varepsilon)\times \mathbf{R})} \leq c(\varepsilon \|D^*_\delta \xi\|_{L^p_{0,\delta}((\Sigma\times[-1,1],g_\varepsilon)\times \mathbf{R})} + \|\pi_A \xi\|_{L^p_{0,\delta}((\Sigma\times[-1,1],g_\varepsilon)\times \mathbf{R})});$$

$$\|\pi_A^\perp \xi\|_{L^p_{1,\delta}((\Sigma\times[-1,1],g_\varepsilon)\times \mathbf{R})} \leq c\varepsilon(\|D^*_\delta \xi\|_{L^p_{0,\delta}((\Sigma\times[-1,1],g_\varepsilon)\times \mathbf{R})} + \|\pi_A \xi\|_{L^p_{0,\delta}((\Sigma\times[-1,1],g_\varepsilon)\times \mathbf{R})}).$$

Proof: Essentially, this follows from Lemma 5.2 in [11]. ∎

**Lemma 5.2.2** *There exists a constant $C_4 = C(p,\delta) > 0$ such that*

$$\|\pi_A(D^*_\delta \xi) - E_u(\pi_A \xi)\|_{L^p_{0,\delta}} \leq C_4 \|\pi_A^\perp \xi\|_{L^p_{0,\delta}},$$

*over the domain $((\Sigma \times [-1,1], g_\varepsilon) \times \mathbf{R})$.*

Proof: The calculation for the difference is same as in [11]. The estimate follows by replacing the weighted Sobolev norm in their Lemma 5.3 [11]. ∎

**Lemma 5.2.3** *For $u \in \mathcal{M}_J(a,b)$ with $\dim \mathcal{M}_J(a,b) = 1$, there exist constant $\varepsilon_1, C_5 > 0$ such that for $0 < \varepsilon < \varepsilon_1$, and all $\xi$,*

$$\|\xi\|_{L^p_{1,\delta}} \leq C_5(\varepsilon \|D^*_\delta \xi\|_{L^p_{0,\delta}} + \|\pi_A(D^*_\delta \xi)\|_{L^p_{0,\delta}});$$

$$\|\pi_A^\perp \xi\|_{L^p_{1,\delta}} \leq C_5 \varepsilon \|D^*_\delta \xi\|_{L^p_{0,\delta}},$$

*on the domain $((\Sigma \times [-1,1], g^\varepsilon) \times \mathbf{R})$.*

Proof: For $u \in \mathcal{M}_J(a,b)$ with $\dim \mathcal{M}_J(a,b) = 1$, we have

$$C_p^s \|\xi\|_{L^p_{0,\delta}} \leq \|E_u^* \xi\|_{L^p_{0,\delta}},$$

from Thereom 4.2.6. Therefore applying to $\pi_A \xi$, one obtains

$$\begin{aligned}
C_p^s \|\pi_A \xi\|_{L^p_{0,\delta}} &\leq \|E_u^* \pi_A \xi\|_{L^p_{0,\delta}} \\
&\leq \|\pi_A(D^*_\delta \xi)\|_{L^p_{0,\delta}} + \|E_u^* \pi_A \xi - \pi_A(D^*_\delta \xi)\|_{L^p_{0,\delta}} \\
&\leq \|\pi_A(D^*_\delta \xi)\|_{L^p_{0,\delta}} + c\|\pi_A^\perp \xi\|_{L^p_{0,\delta}} \\
&\leq \|\pi_A(D^*_\delta \xi)\|_{L^p_{0,\delta}} + c'\varepsilon(\|D^*_\delta \xi\|_{L^p_{0,\delta}} + \|\pi_A \xi\|_{L^p_{0,\delta}}).
\end{aligned}$$

The third inequality follows from Lemma 5.2.2 and the last one follows from Lemma 5.2.1. Hence we choose $\varepsilon_1 < \varepsilon_0$ such that $c'\varepsilon < 1/2$ and the inequality becomes

$$\|\pi_A \xi\|_{L^p_{0,\delta}} \leq c''(\varepsilon \|D_\delta^* \xi\|_{L^p_{0,\delta}} + \|\pi_A(D_\delta^* \xi)\|_{L^p_{0,\delta}}). \tag{5.10}$$

Therefore the estimates follows by applying Lemma 5.2.1 again. ∎

**Proposition 5.2.4** *For $u \in \mathcal{M}_J(a,b)$ with $\dim \mathcal{M}_J(a,b) = 1$, there exist $\varepsilon_1, C_6 > 0$ independent of $u$ and $\varepsilon$ such that for $0 < \varepsilon < \varepsilon_1$, and $\xi \in \Omega^1 \oplus \Omega^0 \oplus \Omega^0(\Sigma, ad)$,*

$$C_6 \|\xi\|_{L^p_{0,\delta}} \leq \|D_\delta^* \xi\|_{L^p_{0,\delta}},$$

*where the inequality is over $((\Sigma \times [-1,1], g^\varepsilon) \times \mathbf{R})$.*

Proof: Theorem 4.2.6 gives us a uniform bound for the first eigenvalue of operator $E_u$ and all holomorphic curve $\mathcal{M}_J(a,b)$. By Lemma 5.2.3

$$\|\xi\|_{L^p_{0,\delta}} \leq \|\xi\|_{L^p_{1,\delta}} \leq C_5(\varepsilon \|D_\delta^* \xi\|_{L^p_{0,\delta}} + \|\pi_A(D_\delta^* \xi)\|_{L^p_{0,\delta}}).$$

Hence the result follows. ∎

Note that the above Proposition is what we need to apply the inverse function theorem and Dostoglou and Soloman in [11] p 19 also pointed out this fact.

# 6 Comparing the Floer boundary maps

## 6.1 Deforming holomorphic curves into anti-self-dual connections

We will deform every holomorphic curves $u \in \mathcal{P}^p_{k,\delta} \cap \mathcal{M}_J(a,b)$ with $\dim \mathcal{M}_J(a,b) = 1$ into anti-self-dual connections. The method is to apply the inverse function theorem, equivalent to the Newton iteration method in [11].

**Proposition 6.1.1** *For $0 < \varepsilon < \varepsilon_1$, there exists a constant $C_7$ independent of $\varepsilon, u \in \mathcal{M}_J(a,b)$ such that the self-duality operator $D$ has a bounded right inverse $G$ with*

$$\|G\xi\|_{L^p_{1,\delta}(Y \times \mathbf{R}, g_\varepsilon)} \leq C_7 \|\xi\|_{L^p_{0,\delta}(Y \times \mathbf{R}, g_\varepsilon)},$$

$$\|G\xi\|_{L^q_{0,\delta}} \leq C_7 \|\xi\|_{L^p_{0,\delta}}, \quad 1/4 + 1/q \geq 1/p.$$

Proof: The result follows from Proposition 5.2.4. The last inequality follows from weighted Sobolev embedding theorem [26]. ∎

Our goal is to deform the holomorphic curve, which by Lemma 5.1.2 is an almost anti-self-dual connection, to a nearby anti-self-dual connection. This is solving the nonlinear anti-self-duality equation

$$F^+(u) + d_u^+ a + (a \wedge a)_+ = 0. \tag{6.1}$$

Proposition 6.1.1 solves the linearized anti-self-duality equation for regular $u$. We shall use the inverse function theorem to deform the almost anti-self-dual connection.

**Lemma 6.1.2** *(c.f. [13]) Let $f : E \to F$ be a $C^1$ map between Banach spaces. Assume that in the first order Taylor expansion $f(\xi) = f(0) + Df(0)\xi + N(\xi)$, $Df(0)$ has a finite dimensional kernel and a right inverse $G$ such that for $\xi, \zeta \in E$*

$$\|GN(\xi) - GN(\zeta)\|_E \leq C(\|\xi\|_E + \|\zeta\|_E)\|\xi - \zeta\|_E$$

*for some constant $C$. Let $\delta_1 = (8C)^{-1}$. Then if $\|Gf(0)\|_E \leq \frac{\delta_1}{3}$, there exists a $C^1$-function*

$$\phi : K_{\delta_1} \to ImG$$

*with $f(\xi + \phi(\xi)) = 0$ for all $\xi \in K_{\delta_1}$ and furthermore we have the estimate*

$$\|\phi(\xi)\|_E \leq \frac{4}{3}\|Gf(0)\|_E + \frac{1}{3}\|\xi\|_E$$

*where $K_{\delta_1} = KerDf(0) \cap \{\xi \in E : \|\xi\|_E < \delta_1\}$.*

Applying Lemma 6.1.2 to $f(a) = F^+(u) + D_u a + (a \wedge a)_+$ with $f(0) = F^+(u), u \in \mathcal{M}_J(a,b)$ and $\dim \mathcal{M}_J(a,b) = 1$, $N(a) = (a \wedge a)_+$, $Df(0) = d_u^+ \oplus d_u^*$ with the bounded right inverse $G$ from Proposition 6.1.1, $E = L^p_{1,\delta} \cap L^q_{0,\delta}(\Omega^1 \oplus \Omega^0 \oplus \Omega^0((\Sigma, g_\varepsilon), ad))$ and $F = L^p_{0,\delta}(\Omega^1 \oplus \Omega^0 \oplus \Omega^0((\Sigma, g_\varepsilon), ad))$, we have the following

**Theorem 6.1.3** *Let $u \in \hat{\mathcal{M}}_J(a,b)$ with $\dim \mathcal{M}_J(a,b) = 1$, and let $\varepsilon_1$ be the constant of Proposition 5.2.4. Then if $0 < \varepsilon < \varepsilon_1$, we can deform $u = A + \Phi ds + \Psi dt$ to a smooth anti-self-dual connection over $Y \times \mathbf{R}$.*

Proof: Using Proposition 6.1.1 and Lemma 5.1.2, we have

$$\|GF^+(u)\|_{L^q_{0,\delta}} \leq C_7\|F^+(u)\|_{L^p_{0,\delta}} \leq C_7 C_3 \varepsilon^{2/p},$$

and $N(a) - N(b) = ((a-b) \wedge a)_+ + (b \wedge (a-b))_+$. We use weighted Hölder inequality and Lemma 7.2 in [26]

$$\|((a-b) \wedge a)_+\|_{L^p_{0,\delta}} \leq \|a-b\|_{L^q_{0,\frac{\delta}{2}}} \|a\|_{L^4_{0,\frac{\delta}{2}}} \leq C_\delta \|a-b\|_{L^q_{0,\delta}} \|a\|_{L^q_{0,\delta}}.$$

So

$$\|GN(a) - GN(b)\|_{L^q_{0,\delta}} \leq C_7 C_\delta \|a-b\|_{L^q_{0,\delta}} (\|a\|_{L^q_{0,\delta}} + \|b\|_{L^q_{0,\delta}}).$$

Thus by Lemma 6.1.2 with $\delta_1 = (8C_7 C_\delta)^{-1}$, there exists $\phi : T_u \mathcal{P}^p_{1,\delta}(a,b) \to \text{Im}G$ with $f(\xi + \phi(\xi)) = 0$; here $\phi(u) = a_u$. So $u + a_u$ is an ASD connection over $(Y \times \mathbf{R}, g_\varepsilon + dt^2)$ with $\|a_u\|_{L^q_{0,\delta}}$ small and is smooth by standard elliptic regularity. ∎

Theorem 6.1.3 provides an injective map

$$T_\varepsilon : \hat{\mathcal{M}}_J(a,b) \to \hat{\mathcal{M}}(a,b),$$

for $0 < \varepsilon < \varepsilon_1$ and $\dim \mathcal{M}(a,b) = dim\mathcal{M}_J(a,b) = 1$. The injectivity follows from the inverse function theorem (see also [10] §7.2). In particular, for $\dim \mathcal{M}_J(a,b) = 1$, the cardinality of holomorphic curve moduli space $\hat{\mathcal{M}}_J(a,b)$ is less than or equal to one of instanton moduli space $\hat{\mathcal{M}}(a,b)$.

## 6.2 Deforming anti-self-dual connection into holomorphic curve

In [11], Dostoglou and Salamon showed that the map $T_\varepsilon$, in fact, is onto for a mapping cylinder with nontrivial $SO(3)$ bundle. Their arguments will not go through for $SU(2)$ Floer homology of a homology 3-sphere, since the analysis has to deal with reducible representations on each handlebody. In addition, the anti-self-dual connection may not lie in the right space for the map $T_\varepsilon$ which consists of flat connections on each handlebody. The curvature of the $SU(2)$ anti-self-dual connection does not neccessarily has the exponential decay property.

Instead of stretching the tube, we first try to deform the anti-self-dual connection into an almost flat connection at each slice, then to put such a deformed anti-slef-dual connection in the holomorphic curve setting. The extra deforming terms make the new connection not to satisfy the holomorphic curve equation. Using the analysis of the uniformly bounded inverse operator of Cauchy-Rieman operator from anti-self-duality operator, we deform this new connection again into a holomorphic curve by inverse function theorem. This produces another injective map from $\hat{\mathcal{M}}(a,b)$ to $\hat{\mathcal{M}}^p_{1,\delta}(a,b)$, thus we get the one to one and onto map between the balanced moduli spaces $\hat{\mathcal{M}}(a,b)$ and $\hat{\mathcal{M}}_J(a,b)$ from Lemma 4.1.3 and Lemma 4.2.3 since the 1-dimensional moduli space $\hat{\mathcal{M}}^p_{1,\delta}(a,b)$ is $\hat{\mathcal{M}}_J(a,b)$ by Proposition 2.2.17. A prior we do not have the way to rule out the singularity, only after all these deformations, we know that at the end the elements in $\hat{\mathcal{M}}^p_{1,\delta}(a,b)$ do not contain any singular point.

### 6.2.1 Deforming ASD into a path in $\mathcal{R}_j$

Let $A \in \mathcal{M}(a,b)$, with $\dim \mathcal{M}(a,b) = 1$, be an anti-self-dual connection on $Y \times \mathbf{R}$. Let $A_j(t)$ be its restriction on section $Y_j \times \{t\}$ for $j = 0, 1, 2, \emptyset$. In order to deform $A$ into a pseudoholomorphic curve, one needs to put $A$ in the right space, namely $A|_{Y_j} \in \mathcal{R}_j, j = 0, 1, 2$ for every $t$. This subsection is devoted to discuss this issue.

Recall the notation from §3.1 or in [34]. We have a section $f_j : \mathcal{B}^*_j \to \mathcal{L}_{jA_j(t)}$ (see §3 (3.4)) over $\mathcal{B}^*_j$. The problem is to deform $A_j(t)$ into the zero set of $f_j^\pi = f_j + grad h_\pi$.

**Proposition 6.2.1** *Let $\nabla f_{jA}^\pi : \mathcal{T}_{jA} \to \mathcal{L}_{jA}$ denote the covariant derivative of $f_{jA}^\pi$ at $A$. Then the followings hold:*

1. $\nabla f_{jA}^\pi(a) = *d_A a + Hess h_\pi[A](a) - d_A u_j(a)$ *where* $Hess h_\pi[A](a) = \frac{d}{ds}(grad h_\pi[A + sa])|_{s=0}$ *and* $u_j(a) \in \Omega^0(Y_j, AdP_j)$ *obeys*

   $$d_A^* d_A u_j(a) - *(F_A \wedge a - a \wedge F_A) - d_A^* Hess h_\pi[A] \cdot a = 0 \quad and \quad i_j^*(u_j(a)) = 0.$$

2. $\nabla f_{jA}^\pi$ *is bounded Fredholm operator with index 0 for $j = \emptyset$, index $3g - 3$ for $j = 1, 2$, and index $6g - 6$ for $j = 0$.*

3. *The difference $\nabla f_{jA}^\pi - \nabla f_{jA}$ is a compact operator.*

4. *The assignment of $\nabla f_{jA}^\pi$ to $[A] \in \mathcal{B}^*_j$ defines a smooth section of $Fred_d(\mathcal{T}_j, \mathcal{L}_j)$ over $\mathcal{B}^*_j$, where $d = Index(\nabla f_{jA}^\pi)$.*

The proof of Proposition 6.2.1 can be found in Proposition 5.1, 6.1, 6.2 of [34].

Now we extend the operator $\nabla f^\pi_{jA_j(t)}$ to

$$D^\pi_{A_j(t)} : L^p_k(\Omega^1 \oplus \Omega^0)(Y_j, adSU(2)) \to L^p_{k-1}(\Omega^1 \oplus \Omega^0)(Y_j, adSU(2)), \qquad (6.2)$$

where the operator $D^\pi_{A_j(t)}$ is given by $D^\pi_{A_j(t)}(\alpha, \beta) = (\nabla f^\pi_{jA_j(t)}\alpha - d_{A_j(t)}\beta, -d^*_{A_j(t)}\alpha)$.

For an irreducible connection $A_j(t) \in \mathcal{B}^*_j$, the covariant Laplacian

$$d^*_{A_j(t)} d_{A_j(t)} : L^2_k \Omega^0(Y_j, adSU(2)) \to L^2_{k-2} \Omega^0(Y_j, adSU(2))$$

is invertible with Neumann condition on the boundary $\partial Y_j$. The inverse of $d^*_{A_j(t)} d_{A_j(t)}$ defines a bounded linear map from $L^2_{k-2} \Omega^0(Y_j, adSU(2)) \to L^2_k \Omega^0(Y_j, adSU(2))$.

The operator $D^\pi_{A_j(t)}$ make sense for any connection $A_j(t) \in \mathcal{A}_j$. The gauge equivariant property of $D^\pi_{A_j(t)}$ defines an assignment from $\mathcal{A}_j$ to the space of bounded, real Fredholm operators from $\mathcal{T}^1_j \oplus \mathcal{T}^0_j$ (see §3 (3.3)) to $\mathcal{L}^1_j \oplus \mathcal{L}^0_j$ (see Lemma A.1 in [34]):

$$\mathcal{T}^1_j = \mathcal{T}_j, \quad \mathcal{T}^0_j = \{b \in L^p_k(\Omega^0(Y_j, adP_j)) | i^*_j(b) = 0\};$$

$$\mathcal{L}^1_j = \mathcal{L}_j, \quad \mathcal{L}^0_j = \{b \in L^p_{k-1}(\Omega^0(Y_j, adP_j)) | i^*_j(b) = 0\}.$$

¿From simple calculation, one has

$$\ker D^\pi_{A_j(t)} = \ker \nabla f^\pi_{jA_j(t)}, \quad Co \ker D^\pi_{A_j(t)} = Co \ker \nabla f^\pi_{jA_j(t)}.$$

For $A_j(t)$ flat connection, we know from (3.5) that

$$\ker \nabla f_{jA_j(t)} = H^1(Y_j, ad\rho_j), \quad Co \ker \nabla f_{jA_j(t)} = H^1(Y_j, \partial Y_j, ad\rho_j),$$

where $\rho_j$ is the holonomy representation of $A_j(t)$.

**Lemma 6.2.2** *The Hilbert space $L^2_k(\Omega^1 \oplus \Omega^0)(Y_j, adSU(2))$ can be decomposed into the following orthogonal decomposition (Hodge decomposition):*

$$L^2_k(\Omega^1 \oplus \Omega^0)(Y_j, adSU(2)) = Im D^\pi_{A_j(t)} \oplus \ker D^\pi_{A_j(t)} \oplus Co \ker D^\pi_{A_j(t)}. \qquad (6.3)$$

∎

### 6.2.1.a Kuranishi picture for $\mathcal{R}_j$

The main concern is the ellipticity for $f_j^\pi = 0$. Therefore we replace the equation $f_j^\pi = 0$ by $p_{A_j} f_j^\pi = 0$, where

$$p_{A_j} : \mathcal{L}_{jA_j} \to Co\ker \nabla f_{jA_j}^\pi$$

is the orthogonal projection onto $Co\ker \nabla f_{jA_j}^\pi$ where $p = 2$ in $\mathcal{L}_{jA_j}$.

**Lemma 6.2.3** *Let $h$ be an admissible perturbation and $A_j$ be a smooth connection with $f_j^\pi(A_j) = 0$. There exists an open $\varepsilon_2$-neighborhood of $A_j$*

$$U_{A_j,\varepsilon_2} = \{A \in \mathcal{A}_j | \quad \|A - A_j\|_{L_k^2} < \varepsilon_2\}$$

*such that if $A \in U_{A_j,\varepsilon_2}$ then*

$$p_A : Co\ker \nabla f_{jA_j}^\pi \to Co\ker \nabla f_{jA}^\pi$$

*is injective.*

Proof: Let $A = A_j + a_j$ and $\alpha \in Co\ker \nabla f_{jA_j}^\pi$ such that $p_A(\alpha) = 0$. Hence $\alpha$ has the following property:

$$(\nabla f_{jA_j}^\pi)^*\alpha = 0, \quad \alpha = \nabla f_{jA}^\pi u, \quad \text{for} \quad u \in \mathcal{T}_{jA}, i_j^*(*\alpha) = 0. \tag{6.4}$$

Thus the composition

$$(\nabla f_{jA}^\pi)^*(\nabla f_{jA_j}^\pi u) = [(\nabla f_{jA}^\pi)^* - (\nabla f_{jA_j}^\pi)^*]\alpha,$$

where the term $[(\nabla f_{jA}^\pi)^* - (\nabla f_{jA_j}^\pi)^*]$ is zero order compact operator. Now the standard bootstrap arguments show that $\alpha \in \mathcal{T}_{jA_j}$ (gain one more derivative). Also the zeroth order compact operator is bounded in the norm in $\mathcal{T}_{jA_j}$, i.e.

$$\|[(\nabla f_{jA}^\pi)^* - (\nabla f_{jA_j}^\pi)^*]\alpha\|_{L_{k-1}^2} \leq C(A_j)\|a_j\|_{L_k^2}\|\alpha\|_{L_k^2}. \tag{6.5}$$

(see also Proposition 4.9 (3) in [34]) Furthermor there is a constant $C(A_j)'$ such that

$$\|\alpha\|_{L_k^2} \leq C(A_j)'\|(\nabla f_{jA}^\pi)^*\alpha\|_{L_{k-1}^2}, \tag{6.6}$$

because $\alpha = \nabla f_{jA}^\pi u$ is perpendicular to $\ker(\nabla f_{jA}^\pi)^*$. From (6.4) and (6.5),

$$\begin{aligned}
\|\alpha\|_{L_k^2} &\leq C(A_j)'\|(\nabla f_{jA}^\pi)^*\alpha\|_{L_{k-1}^2} \\
&= C(A_j)'\|[(\nabla f_{jA}^\pi)^* - (\nabla f_{jA_j}^\pi)^*]\alpha\|_{L_{k-1}^2} \\
&\leq C(A_j)'C(A_j)\|a_j\|_{L_k^2}\|\alpha\|_{L_k^2}
\end{aligned}$$

If $\varepsilon_2$ satisfies $C(A_j)'C(A_j)\varepsilon_2 < 1$, then $\alpha = 0$. ∎

Thus we have showed that the equations $f_j^\pi = 0$ and $p_{A_j} f_j^\pi = 0$ are equivalent for connections in $U_{A_j, \varepsilon_2}$ for $f_j^\pi(A_j) = 0$, i.e. they have the same zero set near $A_j$. Using the Kuranishi deformation technique, one can describe a finite dimensional local model for $(f_j^\pi)^{-1}(0)$. This has been done in [27].

**Theorem 6.2.4** *Fix an admissible perturbation $h$ and $A_j \in (f_j^\pi)^{-1}(0) = \mathcal{R}_j^\pi$. Let $stab(A_j)$ denote the subgroup of gauge group which keeps $A_j$ invariant. There are*

1. *a $stab(A_j)$-equivariant neighborhood $V_{A_j}$ of $0$ in $\ker \nabla f_{jA_j}^\pi$,*

2. *a $\mathcal{G}_j$-equivariant neighborhood $U_{A_j}$ of $A_j$ in $\mathcal{A}_j$,*

3. *a $stab(A_j)$-equivariant real analytic embedding*

$$\phi_{A_j} : V_{A_j} \to U_{A_j} \cap \mathcal{L}_{jA_j},$$

*whose differential at $0$ is the inclusion of $\ker \nabla f_{jA_j}^\pi$ into $\ker d_{A_j}^* \cap \mathcal{L}_{jA_j}$,*

4. *and a $stab(A_j)$-equivariant map $\Phi_{A_j} : V_{A_j} \to Co\ker \nabla f_{jA_j}^\pi$ such that $\phi_{A_j}$ maps $\Phi_{A_j}^{-1}(0)$ homeomorphically onto $(f_j^\pi|_{U_{A_j} \cap \mathcal{L}_{jA_j}})^{-1}(0)$.*

Proof: We define a map

$$G : \ker \nabla f_{jA_j}^\pi \oplus Im(\nabla f_{jA_j}^\pi)^* \to Im(\nabla f_{jA_j}^\pi),$$

by the formula $G(\alpha, \beta) = (Id - p_{A_j}) f_j^\pi(A_j + \alpha + \beta)$. Its differential at $(0,0)$

$$\frac{\partial G}{\partial \beta}(0,0) = (Id - p_{A_j}) \nabla f_{jA_j}^\pi(\beta), \tag{6.7}$$

which is surjective by the orthogonal decomposition. By the implicit function theorem, there is an open neighborhood $V_{A_j}$ of $0 \in \ker \nabla f_{jA_j}^\pi$ and a $stab(A_j)$ equivariant map $\phi'_{A_j} : V_{A_j} \to U_{A_j} \cap \mathcal{L}_{jA_j}$ with

$$p_{A_j} f_j^\pi(A_j + \alpha + \phi'_{A_j}(\alpha)) = 0.$$

For $\alpha \in V_{A_j}$, define $\phi(A_j)(\alpha) = A_j + \alpha + \phi'_{A_j}(\alpha)$ and

$$\Phi(A_j)(\alpha) = p_{A_j} f_j^\pi(\phi(A_j)(\alpha)) \in Co\ker \nabla f_{jA_j}^\pi.$$

(4) follows from the above Lemma and implicit function theorem again. ■

Note that in [24] we have extend the Taubes construction to the reducible connections. For a reducible flat connection $A_j \in \mathcal{R}_j$ satisfying $f_j^\pi(A_j) = 0$, we also have the space $\mathcal{L}_{jA_j}$ and $\mathcal{T}_{jA_j}$ with similar orthogonal decomposition (c.f [24] §3.2).

6.2.1.b **Neighborhood of $\mathcal{R}_j$**

We are going to define a neighborhood of $\mathcal{R}_j$ and show that any connection in such a neighborhood can be deformed into $\mathcal{R}_j$.

**Definition 6.2.5** *Set*

$$U_{\mathcal{R}_j,\varepsilon_2} = \{A \in \mathcal{B}_j \mid \text{ there exists } A_j \in \mathcal{R}_j \text{ such that}$$

$$\|A - A_j\|_{L_k^2} < \varepsilon_2, \quad \|f_j^\pi(A)\|_{L_{k-1}^2} < \varepsilon_2\}, \quad j = 1, 2.$$

**Lemma 6.2.6** *There is a constant $C_8$ independent of $\varepsilon_2$ such that*

$$\|u\|_{L_k^2} \leq C_8 \|(\nabla f_{jA}^\pi)^* u\|_{L_{k-1}^2}$$

*for all $A \in U_{\mathcal{R}_j,\varepsilon_2}$ and $u$ is perpendicular to $\ker((\nabla f_{jA_j}^\pi)^*) \cap \mathcal{L}_{jA_j}$, and $A_j$ is an element in $\mathcal{R}_j$ which is $\varepsilon_2$ close to $A$.*

Proof: Since for $A_j \in \mathcal{R}_j$ the connected component of $A_j$ is compact, there is a constant $C_9$ independent of $A_j$ such that

$$\|u\|_{L_k^2} \leq C_9 \|(\nabla f_{jA_j}^\pi)^* u\|_{L_{k-1}^2} \tag{6.8}$$

for all $u$ perpendicular to $\ker((\nabla f_{jA_j}^\pi)^*) \cap \mathcal{L}_{jA_j}$. By (6.5), we have

$$\begin{aligned}
\|u\|_{L_k^2} &\leq C_9 \|(\nabla f_{jA_j}^\pi)^* u\|_{L_{k-1}^2} \\
&\leq C_9 \|[(\nabla f_{jA_j}^\pi)^* - (\nabla f_{jA}^\pi)^*] u\|_{L_{k-1}^2} + C_9 \|(\nabla f_{jA}^\pi)^* u\|_{L_{k-1}^2} \\
&\leq C_9 C(A_j) \|A - A_j\|_{L_k^2} \|u\|_{L_k^2} + C_9 \|(\nabla f_{jA}^\pi)^* u\|_{L_{k-1}^2}.
\end{aligned}$$

The result follows from choosing $C_9 C(A_j) \varepsilon_2 < \frac{1}{2}$ and $C_8 = 2C_9$. ■

¿From the definition of $f_j^\pi$, we have that $f_j^\pi(A_j + a) = \nabla f_{jA_j}^\pi a + N(a)$ with

$$
\begin{aligned}
N(a) &= *F(A_j + a) + grad_\pi h(A_j + a) - (*d_{A_j}a + Hess(h_\pi)a - d_{A_j}u_j(a)) \\
&= *(a \wedge a) + (grad_\pi h(A_j + a) - grad_\pi h(A_j) - Hess(h_\pi)a) + d_{A_j}u_j(a). \quad (*)
\end{aligned}
$$

Note that $d_{A_j}u_j(a)$ is the term which is determined by projection of $\nabla f_{jA_j}^\pi(a)$ onto the space $\mathcal{L}_{jA_j}$ (c.f. Proposition 6.2.1 (1)).

**Lemma 6.2.7** *Fix a smoothly embedded loop in $Y_j$ and consider the perturbation $h_\pi$. Then there exists a constant $C_{10}$ independent of the connection $A_j$ such that*

$$\|N(a) - N(b)\|_{L^2_{k-1}} \leq C_{10}\|a - b\|_{L^2_k}(\|a\|_{L^2_k} + \|b\|_{L^2_k}),$$

*for any $a, b$.*

Proof: We denote the term

$$H(a) = grad_\pi h(A_j + a) - grad_\pi h(A_j) - Hess(h_\pi)a, \tag{6.9}$$

Then $H(a) - H(b)$ is the only term needed to check in $(*)$. In [34] formula (8.6), Taubes obtains the uniformly bound

$$|\nabla^{(n)}grad_\pi h(A_j)(\{a_i\}_{i=1}^{n-1})| \leq C_{11}\Pi_{i=1}^{n-1}(\|\nabla_{A_j}a_i\|_{L^2(\gamma)} + \|a_i\|_{L^2(\gamma)}),$$

where $\gamma$ is the collection of the loops in $Y_j$. By the mean value theorem, Hölder inequlity and Sobolev embedding theorem,

$$
\begin{aligned}
\|H(a) - H(b)\|_{L^2_{k-1}} &= \|Hessh_\pi(A_j)(a - b) + \nabla^2 grad_\pi h(A_j + \xi)(\{a, a\} - \{b, b\})\|_{L^2_{k-1}} \\
&\leq C_{11}\|a - b\|_{L^2_k}(\|a\|_{L^2_k} + \|b\|_{L^2_k}).
\end{aligned}
$$

∎

**Proposition 6.2.8** *Let $\varepsilon_2$ be the one in Lemma 6.2.3. Then if $0 < \varepsilon \leq \varepsilon_2$, then any connection $A \in U_{\mathcal{R}_j, \varepsilon}$ can be deformed into a smooth connection $A + a_j$ on $Y_j$, $j = 1, 2$ with*

$$f_j^\pi(A + a_j) = 0, \quad \|a_j\|_{L^2_k} < \varepsilon.$$

Proof: This is a consequence of Lemma 6.1.2. Applying to $f(a_j) = f_j^\pi(A + a_j)$ for $A = A_j + a \in U_{\mathcal{R}_j,\varepsilon}$ with $f(0) = f_j^\pi(A)$, $Df(0) = \nabla f_{jA}^\pi$ has bounded right inverse from Lemma 6.2.6 and Lemma 6.2.7. Taking $E = L_k^2 \cap \mathcal{T}_{jA}$ and $F = L_{k-1}^2 \cap \mathcal{L}_{jA}$ in Lemma 6.1.2, the result follows. ∎

Now the task is to push an anti-self-dual connection $A_j(t)$ on $Y_j \times \{t\}$ for $t \in \mathbf{R}$ into the neighborhood $U_{\mathcal{R}_j,\varepsilon}$ for $0 < \varepsilon \leq \varepsilon_2$.

**Deformation of $A_j(t)$ into $\mathcal{R}_j$ on $Y_j \times \mathbf{R}$**

Recall that an anti-self-dual connection $A(t)$ decays exponentially to the irreducible deformed-flat connection $A(\pm\infty) \in \mathcal{R}^\pi(Y)$. A connection $a_j$ on $Y_j$ is called deformed-flat if $f_j^\pi(a_j) = 0$ for $j = 0, 1, 2, \emptyset$. $A(\pm\infty) \in \mathcal{R}^\pi(Y) = \mathcal{R}_\pi^*(Y_1) \cap \mathcal{R}_\pi^*(Y_2)$, hence we denote $A_j(\pm\infty) \in \mathcal{R}_\pi^*(Y_j), j = 0, 1, 2, \emptyset$. Morgan, Mrowka and Ruberman in [27] have extensively studied the ASD connection around the flat connections at the ends. Our approach is to get a similar $\delta$-decay property.

### 6.2.1.c Deforming around infinity

**Lemma 6.2.9** *For $\alpha_j \in \mathcal{R}_\pi^*(Y_j), j = 1, 2$, there exists a constant $C_j > 0$ which is independent of $\alpha_j$ such that*

$$\|a_j\|_{L_k^2} \leq C_j \|D_{\alpha_j}^\pi a_j\|_{L_{k-1}^2}, \tag{6.10}$$

*for $a_j \in (\mathcal{T}^1 \oplus \mathcal{T}^0)_{j\alpha_j} \cap (\ker D_{\alpha_j}^\pi)^\perp$.*

Proof: Note that $\text{Coker} D_{\alpha_j}^\pi = H^1(Y_j, \partial Y_j, ad\alpha_j) = 0$ for $\alpha_j \in \mathcal{R}_\pi^*(Y_j), j = 1, 2$. The (6.11) is true for $\alpha_j \in \mathcal{R}_\pi^*(Y_j), j = 1, 2$, with $C_j$ depending continuously on $\alpha_j$. Then the result follows from the compactness of $\mathcal{R}_\pi^*(Y_j), j = 1, 2$. ∎

Now the operator $(D_{\alpha_j}^\pi)^* D_{\alpha_j}^\pi$ is elliptic with the corresponding boundary condition. Its index equals to $3g-3$ for $j = 1, 2$. It has pure point spectrum, all real. The multiplicity of any eigenvalue is finite and there are no accumulation points. Let $\delta_j'$ be a positive number which is smaller than the smallest absolute value of the eigenvalues of $(D_{\alpha_j}^\pi)^* D_{\alpha_j}^\pi$. We define a new norm around the slice of $A_j(\pm\infty)$ to get the exponential decay property in $t$-direction.

Let $\mathbf{B}(k)$ be a Banach space of continuous path $a_j : [T_0, \infty) \to (\mathcal{T}^1 \oplus \mathcal{T}^0)_{jA_j(\infty)}$ (similarly for $(-\infty, -T_0]$) such that

$$\|a_j\|_{B(k)} = \sup_{|t| \geq T_0} e^{\delta_j |t|} \|a_j(t)\|_{L^2_k(Y_j)} < \infty, \tag{6.11}$$

$\delta_j = \min\{\delta'_j/2, \gamma/2\}$ with $\gamma$ in Proposition 4.1.1. (4.2) (4.3). Here $T_0$ is determined by the following.

**Lemma 6.2.10** *For any $\varepsilon_3 > 0$, there exist $T_j > 0$ and a gauge transformation $g_j$ such that for $T_0 = \max\{T_1, T_2\}$,*

$$\|g_j^* A_j - A_j(\pm\infty)\|_{B(k)} < \varepsilon_3, \quad \|f_j^\pi(A_j)\|_{B(k-1)} < \varepsilon_3.$$

Proof: By Proposition 4.1.1 (3), we have that

$$|\nabla^{(l)}(g_j^* A_j - A_j(\pm\infty))|_{y,t} \leq C_A e^{-\gamma|t|}.$$

Hence by the choice of $\delta_j$,

$$e^{\delta_j |t|} \|(g_j^* A_j - A_j(\pm\infty))\|_{L^2_k(Y_j)} \leq C'_A e^{(\delta_j - \gamma)|t|},$$

so that there exists a $T_j$ such that for $|t| \geq T_j$, we have $C'_A e^{(\delta_j - \gamma)|t|} < \varepsilon_3$. The term $f_j^\pi(A_j)$ can also be estimated in a similar way by using Proposition 4.1.1 (2) for the perturbed ASD connections (see §9 in [19]). ∎

**Lemma 6.2.11** *Let $a'_j(t) = g_j^* A_j(t) - A_j(\pm\infty)$ for $g_j$ in Lemma 6.2.10. Then for $\alpha_j \in \mathcal{R}^*_\pi(Y_j)$ and $|t| \geq T_0$, we have*

1. *$a'_j(t)$ is perpendicular to $\ker D^\pi_{\alpha_j}$,*

2. *$\delta_j^{1/2} \|a'_j\|_{B(k)} \leq \|D^\pi_{\alpha_j} a'_j\|_{B(k-1)}$.*

Proof: The first claim that $a'_j(t)$ can be putted into the slice of $\alpha_j$ for all $|t| \geq T_0$ follows from Theorem 2.6.3 in [27]. The finiteness of $\|D^\pi_{\alpha_j} a'_j\|_{B(k-1)}$ is from Proposition 4.1.1. Hence the result follows by multiplied $e^{\delta_j |t|}$ on

$$\delta'_j \|a'_j(t)\|_{L^2_k(Y_j)} \leq \|D^\pi_{\alpha_j} a'_j(t)\|_{L^2_{k-1}(Y_j)},$$

and taking sup-norm for $|t| \geq T_0$. ∎

We are going to find a path $a_j(t) \in \mathbf{B}(k)$ such that $A_j(t) + a_j(t) \in \mathcal{R}^*_\pi(Y_j)$. Define an open neighborhood of $A_j(\pm\infty)$ on the slice $\mathcal{T}^1_{jA_j(\pm\infty)}$ by

$$U_{\delta_j\varepsilon_j,\pm} = \{A_j \in \mathcal{A}_j | \ \|A_j - A_j(\pm\infty)\|_{B(k)} < \varepsilon_j, \|F_{A_j}\|_{B(k-1)} < \varepsilon_j, A_j - A_j(\pm\infty) \in \mathcal{T}_{jA_j(\pm\infty)}\}.$$

**Proposition 6.2.12** *For an anti-self-dual connection $A$ and any $\varepsilon_j > 0$, there exist $T_0 > 0$ and a path $a_j : [T_0, \infty) \to \mathcal{L}_{jA_j(\infty)}$ such that $A_j(t) + a_j(t) \in \mathcal{R}^\pi_j$ is deformed-flat connection for all $t \geq T_0$ with $a_j(t)$ sufficiently small $\|a_j\|_{B(k-1)} < \varepsilon_j$. Similarly for $t \leq -T_0$.*

Proof: By Lemma 6.2.10, the path $A_j : [T_0, \infty) \to \mathcal{L}^1_{jA_j(\infty)}$ is in the Banach subspace $\mathbf{B}(k-1) \cap U_{\delta_j\varepsilon_j,+}$. The deformed-flat equation $(f^\pi_j(A_j(t) + a_j(t)) = 0)$

$$D^\pi_{A_j(t)}a_j(t) + *(a_j(t) \wedge a_j(t)) + *F_{A_j(t)} + gradh_\pi(A_j(t)) = 0 \quad (6.12)$$

has a solutions which is parametrized by $\ker D^\pi_{A_j(t)}$ from Lemma 6.1.2, Lemma 6.2.10 and Lemma 6.2.11 and $\|a_j\|_{B(k-1)} < \varepsilon_j$. ∎

The above $A_j(t) + a_j(t)$ is deformed-flat for every $|t| > T_0 = T$ and $a_j(t)$ is $\delta_j$-decay in the $t$-direction from the definition of the norm $B(k)$.

6.2.1.d **Changing metrics**

Let $H$ be a hanbdlebody $\#D^2 \times S^1$ with boundary isomorphic to $\Sigma$, $\partial H = \Sigma$. By inserting $3g - 3$ disjoint disks $D^2_a \subset H$ which cut the handlebody into 3-balls, we obtain a decomposition $\{\Sigma_\varepsilon : 0 \leq \varepsilon \leq 1\}$ of $\Sigma = \Sigma_1$ by shrinking the curves $\partial D^2_a$ to curves with smaller radii in $D^2_a$. At the end of the deformation $\Sigma_0$, we obtain a union of 2-spheres with 3 marked points, and filling in these 2-spheres by 3-balls we recover the handlebody $H$.

For genus $g$ handlebody $H$, the core $\mathbf{C}$ of the handlebody $H$ consists of $2g - 2$ vertices and $3g - 3$ arcs. Choose a Riemannian metric $g_H$ on $H$ such that on $H$ is a product metric near the boundary. Let $Y_1, Y_2$ both have such a metric and let the shrinking metrics $g^\varepsilon_i, i = 1, 2$ on the handlebodies $Y_1, Y_2$ deforming then into its core. As we shrink the metrics on $Y_1$ and $Y_2$, the $A_j(t)$ on the handlebody corresponds to become (deformed-)flat.

For each handle $S^1 \times D^2$, the shrinking metric

$$g_{H\delta} = d^2\theta + \delta_0^2 dg_{D^2}, \tag{6.13}$$

with $\theta$ coordinate of $S^1$ and patching $g_{H\delta}$ smoothly together along $S^2$ as in [25] §3.2 (ii) (d), the resulting metric is $g_{j\delta}$ on $Y_j$ handlebody (see figure 1). In order to patch metric $g_{j\delta}$ with $g_\varepsilon = g^\varepsilon$ (as in §5.1), we need to take an average metric along the overlap, which is similar construction in [25] §3.3. The overlap region is a small annulus $U_1 = \Sigma_\delta \times [-1 + S^{-1}\varepsilon, -1 + S\varepsilon]$ on $Y_1$ side, $U_2 = \Sigma_\delta \times [1 - S\varepsilon, 1 - S^{-1}\varepsilon]$ on $Y_2$ side. Here $S(>1)$ is another parameter (to be fixed later in the proof) with $S\varepsilon < 10^{-2}$ (say), the metric on $\Sigma_\delta$ is shrinking along $D^2$-direction, rather than $S^1$-direction. We only discuss the patching on $Y_1$ side, the other side is same.

Define $f_1 : U_1(\subset Y_1) \to \Sigma_\varepsilon \times [-1 + S^{-1}\varepsilon, -1 + S\varepsilon]$ by $f_1(x,s) = (x, -s - 2 + S^{-1}\varepsilon + S\varepsilon)$, $\Sigma_\varepsilon$ is the surface $\Sigma$ with metric $g_\varepsilon$. The linear inversion map $f_1$ taking the surface $\Sigma_\delta \times \{-1 + S^{-1}\varepsilon\}$ to $\Sigma_\varepsilon \times \{-1 + S\varepsilon\}$ induces an orientation-reversing diffeomorphism from $U_1$ to $\Sigma_\varepsilon \times [-1 + S^{-1}\varepsilon, -1 + S\varepsilon]$. The map $f_2$ is given by $f_2(x,s) = (x, -s + 2 - S^{-1}\varepsilon - S\varepsilon)$ from $U_2$ to $\Sigma_\varepsilon \times [1 - S\varepsilon, 1 - S^{-1}\varepsilon]$. Thus we have our homology sphere $Y$ to be

$$Y_1 \#_{f_1} (\Sigma \times [-1 + S^{-1}\varepsilon, 1 - S^{-1}\varepsilon]) \#_{f_2} Y_2,$$

where the annuli $U_j$ are identified by $f_j$.

We first extend the metric $g_{j\delta}$ on $Y_j$ by a wrapped metric on the annulus $U_j$. Let $\phi_j$ be a monotone cutoff function satisfying:

$$\phi_j(s) = 1, |s| \geq 1 - S^{-1}\varepsilon; \quad \phi_j(s) = \varepsilon^2, |s| \leq 1 - S\varepsilon.$$

The metric $g_{j\delta}$ on $U_j$ is given by

$$ds^2 + s^2(\phi_j(s)d\theta^2 + \delta_0^2 dg_{\partial D^2}^2), \tag{6.14}$$

in terms of each handle.

**Definition 6.2.13** *The Riemannian metric $g_Y$ on $Y$ is defined as follows:*

1. On $Y_j$, $g_Y = g_{j\delta}$ on $Y_j, j = 1, 2$;

2. On $U_j$, $g_Y = \chi_j f_j^*(g_{j\delta}) + (1-\chi_j)g_\varepsilon = \chi_j g_{j\delta} + (1-\chi_j)f_j^*(g_\varepsilon)$, because of the linearity of $f_j$.

3. On $\Sigma \times [-1+S\varepsilon, 1-S\varepsilon]$, $g_Y = g_\varepsilon$.

Here $\chi_j \in C^\infty[0, \infty)$ satisfies

$$\chi_j(s) = 1, |s| \geq 1 - S^{-1}\varepsilon; \quad \chi_j(s) = 0, |s| \leq 1 - S\varepsilon.$$

More specifically, we will define $\chi_j$ as follows. Let $s_0 = -1 + S^{-1}\varepsilon, s_1 = -1 + S\varepsilon$. Then we define a function

$$F(s) = (\frac{-s + s_0 + s_1}{s})^2 \phi_1(-s + s_0 + s_1) - \varepsilon^2. \tag{6.15}$$

The function $F(s)$ is an monotone increasing function with

$$F(s_0) = (\frac{s_1^2}{s_0^2} - 1)\varepsilon^2 < 0, \quad F(s_1) = \frac{s_0^2}{s_1^2} - \varepsilon^2 > 0.$$

There exists a unique $s_* \in (s_0, s_1)$ such that $F(s_*) = 0$. Then we define the cutoff $\chi_1$ to be a smoothing function of $\chi$:

$$\chi(s) = \frac{F(s)}{F(s_0)}, s_0 \leq s \leq s_*; \quad \chi(s) = 1, s \leq s_0; \quad \chi(s) = 0, s \geq s_*. \tag{6.16}$$

I.e $\chi_1 = \chi * \rho$, where $\rho(s)$ is a mollifier (see [20] §7.2) with $\|\chi_1\|_{C^0} \leq \|\chi\|_{C^0}$.

The metric $g_Y$ on $\Sigma \times [-1+S\varepsilon, 1-S\varepsilon]$ is the one for any $0 < \varepsilon < \varepsilon_1$ (in Theorem 6.1.3). We will choose $\min\{\delta_3, \delta_4\}$(see §6.2.1.e below) and then pick $\delta_0 = \varepsilon \leq \min\{\delta_3, \delta_4, \varepsilon_1\}$. The metric on the homology 3-sphere $Y$ is the one for $\delta_0, \varepsilon$ sufficiently small so that one can deform an anti-self-dual connection into a holomorphic curve.(see figure 2)

**Lemma 6.2.14** *For any $\varepsilon > 0$, there exists $S_0 > 1$ such that for all $1 < S \leq S_0$ with $S_0\varepsilon < 10^{-2}$, we have*

$$\|g_Y - g_\varepsilon\|_{C^0} < \varepsilon, \quad \|g_Y - g_{j\delta}\|_{C^0} < \varepsilon, \quad on \ U_j.$$

Proof: By definition (6.2.13), we calculate the $C^0$-norm of $g_Y - g_\varepsilon, g_Y - g_{j\delta}$ on the annulus region $U_j$.

$$(g_Y)_{ss} = (g_\varepsilon)_{ss}, \quad (g_Y)_{\partial D^2} = (g_\varepsilon)_{\partial D^2} \ (\text{for } \varepsilon = \delta).$$

$$(g_Y)_{\theta\theta} = \chi_j(s)(f_j(s)^2 \phi_j(f_j(s))) + (1-\chi_j)s^2\varepsilon^2, \quad \text{on } U_j.$$

Then we have

$$(g_Y)_{\theta\theta} - (g_\varepsilon)_{\theta\theta} = \frac{1}{\varepsilon^2}\chi_1(s)F(s)(g_\varepsilon)_{\theta\theta}. \tag{6.17}$$

¿From the definition of $\chi_1(s)$ and monotonicity of $F(s)$, we obtain

$$\begin{aligned} \|g_Y - g_\varepsilon\|_{C^0} &\leq 4\left|\frac{s_1^2}{s_0^2} - 1\right| \\ &\leq C\max\{S^2 - 1, |S^{-2} - 1|\} \cdot \varepsilon. \end{aligned}$$

Similarly, we get $(g_Y)_{\theta\theta} = \chi_j(s)s^2\phi_j(s) + (1-\chi_j(s))f_j^2(s)\varepsilon^2$.

$$(g_Y)_{\theta\theta} - (g_{1\delta})_{\theta\theta} = (1-\chi_1(s))(-1 + (\frac{-s+s_0+s_1}{s})^2\frac{\varepsilon^2}{\phi_1(s)})(g_{1\delta})_{\theta\theta}. \tag{6.18}$$

Therefore

$$\begin{aligned} \|g_Y - (g_{1\delta})_{\theta\theta}\|_{C^0} &\leq 4\left|\frac{s_0^2}{s_1^2} - 1\right| \\ &\leq C\max\{S^2 - 1, |S^{-2} - 1|\} \cdot \varepsilon. \end{aligned}$$

By choosing $S_0$ close to 1 enough to make

$$C\max\{S_0^2 - 1, |S_0^{-2} - 1|\} \leq \frac{1}{2},$$

we thus prove the lemma. ∎

The above lemma tells us that we may glue $Y_j$ and $\Sigma_\varepsilon \times [-1+S^{-1}\varepsilon, 1-S^{-1}\varepsilon]$ on the tiny overlap region for $1 < S \leq S_0$. Now we fix the parameter $S$ for the region. For forms $u$ supported on $\Sigma_\varepsilon \times [-1+S^{-1}\varepsilon, 1-S^{-1}\varepsilon]$, we have

$$\frac{1}{2}\|u\|_{L^p(g_\varepsilon)} \leq \|u\|_{L^p(g_Y)} \leq 2\|u\|_{L^p(g_\varepsilon)}. \tag{6.19}$$

(Similar for $g_{j\delta}$) Let $\pi_+^g$ be the projection onto self-dual 2-forms with respect to the metric $g$. Note that $\pi_+^g$ is a continuous map with respect to the metrics, (c.f. [10] and [25] )

$$\|\pi_+^{g_Y} - \pi_+^{g_\varepsilon}\| \leq C\|g_Y - g_\varepsilon\|_{C^0}. \tag{6.20}$$

### 6.2.1.e Deforming on the compact piece $Y_j \times [-T, T]$

To deform $A_j(t)$ into $\mathcal{R}_j^\pi$ for $|t| \leq T$, we need to shrink the metrics on handlebodies. This gives a control of the $L_k^p$-norm of the curvature of the anti-self-dual connection $A_j(t)$ on $Y_j \times [-T, T]$ for $j = 1, 2$. Then the Kuranishi technique pushes $A_j(t)$ into the deformed-flat connection path. Combining with Lemma 6.2.1, we have derived a deformation from an anti-self-dual connection to a connection $A_j(t) + a_j(t) \in \mathcal{R}_j^\pi$ for all $t \in \mathbf{R}, j = 1, 2$.

As in [24], we choose the framed loops $\gamma = \{\gamma_j\}_{1 \leq j \leq m}$ so that $m = 2g$ and the homotopy classes of $\gamma_1(S^1 \times 0), \cdots, \gamma_{2g}(S^1 \times 0)$ are the set of generators in $\pi_1(\Sigma)$. Given such a choice, the space $L_m = \Pi^m SU(2)/SU(2)$ can be identified with the repreentation space of $\Sigma - \{p_0\}$ the punctured Riemann surface with based point $p_0$. In order to understand the deformed-flat connection, we need to explicit to describe the deformed-flat connections on one of the solid tori $\gamma_i(S^1 \times D^2)$. Denote $\lambda_i = \gamma_i(S^1 \times 0)$ and $\mu_i = \gamma_i(0 \times \partial D^2)$.

Let $h_i : SU(2) \to \mathbf{R}$ be a smooth function which is invariant under the adjoint action of $SU(2)$. Let $h$ be a function from $\mathcal{A}(Y)$ to $\mathbf{R}$ in (3.7) with $m = 2g$. The $L^2$ gradient of $h$ is

$$grad h(A) = \sum_i \nabla_i h * \mu,$$

with $\nabla_i h = h_i'(hol_{\gamma_i}(x, A))$ is the partial derivative of the lifting of $h$ to $SU(2)^{2g}$ in the direction of the $i$th factor, identified with an element of $\mathbf{su(2)}$ by virtue of the canonical bilinear form on $\mathbf{su(2)}$.

(a) The deformed-flat connection $A$ satisfies

$$f_j^\pi(A) = *F_A + grad_\pi h(A) = 0, \quad \text{on } Y_j.$$

(b) The corresponding gradient-like flow equation for an one parameter family $A_j(t)$ is

$$\frac{\partial A_j(t)}{\partial t} = f_j^\pi(A_j(t)).$$

(c) The corresponding perturbed anti-self-dual equation over $Y_j \times \mathbf{R}$ is

$$F_{A_j}^+ + \sum_{i=1}^{2g} \nabla_i h \mu^+ = 0,$$

where $\mu^+$ is the self-dual component of $\mu$ in $Y \times \mathbf{R}$ with product metric.

**Lemma 6.2.15** *Let $A_j \in \mathcal{R}_j$ be a deformed-falt connection.*

1. *Then $A_j$ is flat on $Y_j \setminus (\cup_{i=1}^m \gamma_i(S^1 \times D^2))$.*

2. *If the holonomy around longitude $\lambda_i$ is $\exp(i\beta_i)$ under suitable trivialization, then the holonomy around meridian $\mu_i$ is $\exp(i\alpha_i)$ related by $\alpha_i = f_i'(\beta_i)$, where $f_i$ is a smooth $2\pi$ periodic function.*

3. *Restriction to $Y_j \setminus (\cup_{i=1}^m \gamma_i(S^1 \times D^2))$ gives the one-to-one correspondence between $\mathcal{R}_j$ and the gauge equivalence classes of flat connections over $Y_j \setminus (\cup_{i=1}^m \gamma_i(S^1 \times D^2))$ such that the holonomies $\alpha_i, \beta_i$ are related in (2).*

4. $A_j|_{\gamma_i(S^1 \times D^2)}$ *is gauge equivalent to a connection matrix form*

$$\overline{A_j}(\theta, d) = \exp(i\beta_i) d\theta + A_{j1} dx_1 + A_{j2} dx_2,$$

*where $(\theta, x_1, x_2)$ are the coordinates on $S^1 \times D^2$ and $A_{j1}, A_{j2}$ depend only on $h_i$ and $\exp(it) = \begin{pmatrix} e^{it} & 0 \\ 0 & e^{-it} \end{pmatrix}$*

Proof: (1) follows from the definition of $h$. Note that $h_i(hol_{\gamma_i}(x, A_j))$ is determined by a real valued function $f_i$ on $[0, 2\pi]$ by

$$hol_{\gamma_i}(x, A_j) = \exp(it), \quad f_i(t) = h_i(\exp(it)).$$

Any flat $SU(2)$ connection over this torus boundary reduces to $U(1)$ reducible connection so is determined by $\alpha_i, \beta_i$. These give a complete descrption of the gauge equivalence classes fo flat connections. (2) and (3) are proved in [5] Lemma 4. Note that when the holonomy around $\lambda_i$ is ad-trivial, the connection is flat and the results follows more easily. (4) follows from (2) and (3). ∎

The perturbed anti-self-dual connection $A_j(t)$ on $\gamma_i(S^1 \times D^2), j = 1, 2$ can be transformed into $A_j'(t) = \exp(i\beta_i) d\theta + B_j(t)$ by gauge transformation $g_i(t)$:

$$\frac{\partial g_i(t)}{\partial \theta} + (A_j(t))_\theta g_i(t) = i\beta_i, \quad g_i(t)(0) = Id. \tag{6.21}$$

The $B_j(t)$ is a matrix of 1-forms on $0 \times D^2$. This is an inhomogenouse ordinary differential equation, which has a solution varying smoothly with respect to $t \in [-T, T]$.

**Lemma 6.2.16** *For any $0 < \varepsilon_4 \leq \varepsilon_2$, there exist a metric $g_{j\delta_*}$ on $Y_j$ and an element $A_j \in \mathcal{R}_j$ such that for $0 < \delta \leq \delta_*$ and $t \in [-T, T]$*

$$\|A_j(t) - A_j\|_{L^2_k(Y_j, g_{j\delta})} < \varepsilon_4, \quad \|f_j^\pi(A_j(t))\|_{L^2_{k-1}(Y_j, g_{j\delta})} < \varepsilon_4.$$

Proof: Note that on $\gamma_i(S^1 \times D^2)$ the deformed-flat connection $A_j$ can be written as $\exp(i\beta_i)d\theta + A_{j1}dx_1 + A_{j2}dx_2$, where $A_{j1}, A_{j2}$ are related to $h_i(hol_{\gamma_i}(x, A_j))$. Then by Proposition 1.5 in [34], $grad_\pi h(A_j)$ and $h_\pi(A_j)$ both lie in $[-\varepsilon, \varepsilon]$ for $\varepsilon < \varepsilon_4/2$. Take $\varepsilon_4 = \min\{\varepsilon_0/2, \varepsilon_2/2\}$ for $\varepsilon_0$ in Proposition 1.5 [34] and $\varepsilon_2$ in Proposition 6.2.8. Then we have

$$\|A_{j1}dx_1 + A_{j2}dx_2\|_{L^2_k(Y_j, g_{j\delta})} < \frac{\varepsilon_4}{2}.$$

The connection matrix $A'_j(t) - A_j(t) = B_j(t)$ is a matrix of 1-forms on the disk $D^2$ and the perturbed ASD $A_j(t)$ is smooth on the $Y_j \times \mathbf{R}$, in particular $\|A'_j(t)\|_{C^k(Y_j \times [-T,T])} \leq M_j$. Thus shrinking the metric $g_{j\delta}$ on $S^1 \times D^2$ gives that

$$\begin{aligned}\|A'_j(t) - A_j\|_{L^2_k(Y_j, g_{j\delta})} &\leq \|A_{j1}dx_1 + A_{j2}dx_2\|_{L^2_k(Y_j, g_{j\delta})} + \|B_j(t)\|_{L^2_k(Y_j, g_{j\delta})} \\ &< \frac{\varepsilon_4}{2} + M_j Vol(Y_j, g_{j\delta})\end{aligned}$$

Thus choosing $g_{j\delta_3}$ small enough such that $M_j Vol(Y_j, g_{j\delta_3}) < \varepsilon_4/2$. Then first estimate follows.

Note $f_j^\pi(A_j(t)) = *F(A_j(t)) + grad_\pi h(A_j(t))$. We have that $\|grad_\pi h(A_j(t))\|_{L^2_{k-1}} < \varepsilon$ from [34]. The curvature term $F(A_j(t))$ has pointwise estimates from Proposition 4.1.1 (2), so there is a uniformly bounded constant $N_j$ such that

$$\|F(A_j(t))\|_{L^\infty(Y_j \times [-T,T])} \leq N_j,$$

therefore $\|F(A_j(t))\|_{L^2_{k-1}(Y_j, g_{j\delta})} \leq N_j Vol(Y_j, g_{j\delta})$. Choose $\delta_4$ such that $N_j Vol(Y_j, g_{j\delta}) < \varepsilon_4/2$, then we can fix the shrinking metric $g_{j\delta}$ for $\delta \leq \delta_* = \min\{\delta_3, \delta_4\}$. The results follows by adding the above two inequalities. ∎

**Theorem 6.2.17** *Given an (perturbed) anti-self-dual connection $\{A_j(t)\}$ and $\varepsilon_4 > 0$, there exists a $\delta_0 = \min\{\delta_3, \delta_4, \varepsilon_1\}$ such that for any $t \in [-T, T]$, $0 < \delta < \delta_0$, the connection $A_j(t)$ can be deformed into a smooth deformed-flat connection, i.e. $A_j(t) + a_j(t) \in \mathcal{R}_j$ with $\|a_j(t)\|_{L^2_k(Y_j, g_{j\delta})} < \varepsilon_4$.*

Proof: This is just a corollary of Theorem 6.1.2 by using Lemma 6.2.9 and Lemma 6.2.16. ∎

**Remarks:** (1) Theorem 6.2.17 asserts that $A_j(t)$ can be deformed into $\mathcal{R}_j$, not necessary into $\mathcal{R}_j^\pi$. It may happen that for some $t \in [-T, T]$, $A_j(t)$ be deformed into an reducible deformed-flat connection, rather than an element in $\mathcal{R}_j^\pi = (f_j^\pi)^{-1}(0) \cap \mathcal{B}_j^*$ (irreducibles).

(2) $\{a_j(t)\}$ for $t \in [-T, T]$ is in $L^2_{k,\delta}(Y_j \times [-T, T])$, hence also in $L^p_{1,\delta}(Y_j \times [-T, T])$. The other deformed piece of $\{a_j(t)\}$ for $|t| > T$ are also with enough regularity since at $t \to \pm\infty$, the elements in $\mathcal{R}_j^\pi$. So $\{A_j(t) + a_j(t)\}$ for $t \in \mathbf{R}$ is in Banach space $L^p_{1,\delta}(Y_j \times \mathbf{R})$.

### 6.2.2 Deforming ASD into $\mathcal{R}_\Sigma$ on $Y_0 = \Sigma \times [-1 + S^{-1}\varepsilon, 1 - S^{-1}\varepsilon]$

In the above, we have put an anti-self-dual conneciton in the right boundary condition. But $A + \{a_j(t)\}$ is not defined over $Y_0$. Using the normal components of $a_j(t)$ vanishing along the boundary $\Sigma$, we extend $a_j(t)$ inside $Y_0$ a little bit by the following cutoff functions.

**Definition 6.2.18** *Define a smooth function $\eta_j(s), j = 1, 2$ such that*

$$\eta_j(s) = 1, |s| \leq 1 - S^{-1}\varepsilon, \quad \eta_j(s) = 0, |s| \geq 1 - S\varepsilon.$$

*Then the connection $\overline{A} = A + \eta_1(s)a_1(t) + \eta_2(s)a_2(t)$ is a well-defined $SU(2)$ connection over $Y \times \mathbf{R}$*

Such a connection $\overline{A}$ has the following properties: $\overline{A}|_{Y_j \times \{t\}} \in \mathcal{R}_j, j = 1, 2$; $\overline{A}$ is an ASD on $(\Sigma \times [-1 + S\varepsilon, -1 + S\varepsilon]) \times \mathbf{R}$, $\overline{A}(t, s)$ is not necessary a flat connection on $\Sigma \times (t, s)$. Now we further deform $\overline{A}$ over $Y_0 \times \mathbf{R}$ to be in $\text{Map}[\mathbf{R} \times [-1 + S^{-1}\varepsilon, 1 - S^{-1}\varepsilon], \mathcal{R}_\Sigma]$. Here we state the standard Kuranishi map for Riemann surface $\Sigma$.

Let $P_\Sigma \to \Sigma$ be a principal $SU(2)$ bundle and a flat connection $\alpha$ (In our case here, $P_\Sigma = \Sigma \times SU(2)$ a trivial bundle). Let $\mathcal{G}_\Sigma$ be the $L^2_2$ gauge transformations of $P_\Sigma$ and $\mathcal{A}_\Sigma$ the $L^2_1$ connectionss space. Denote $S_\alpha = \alpha + \ker d_\alpha \cap L^2_1(\Omega^1(\Sigma, adSU(2)))$ the slice to the action of $\mathcal{G}_\Sigma$ on $\mathcal{A}_\Sigma$ at $\alpha$, $\pi_\alpha$ the $L^2$ projection of $L^2(\Omega^2(\Sigma, adSU(2))) \to H^2(\Sigma, ad\alpha)$.

**Proposition 6.2.19** *Let $\alpha$ be a flat connection on $P_\Sigma$. Then there exist*

1. *a $stab(\alpha)$-invariant neighborhood $V$ of $0$ in $H^1(\Sigma, adSU(2))$*

2. a $\mathcal{G}_\Sigma$-invariant neighborhood $U$ of $\alpha$ in $\mathcal{A}_\Sigma$

3. a stab($\alpha$)-real analytic embedding $\phi : V \to S_\alpha \cap U$ whose differential at $0$ is the natural inclusion of $H^1(\Sigma, ad\alpha)$ into $\ker(d_\alpha)^* \cap \Omega^1(\Sigma, adSU(2))$

4. The Kuranishi map $\Phi : V \to H^2(\Sigma, ad\alpha)$ is given by

$$\Phi(w) = -\frac{1}{2}\pi_\alpha([w \wedge w]).$$

The proof of Proposition 6.2.19 is quite standard and can be found in [27] for example.

We have given the properties of flat connections on $Y_j$, $j = 0, 1, 2, \emptyset$ in Theorem 6.2.4 and Proposition 6.2.19. In [11], there are only irreducible and nondegenerated flat connections over $\Sigma$ due to the nontrivial $w_2$ condition. But we have to include into our discussion of possible reducible flat connections over the trivial bundle.

**Lemma 6.2.20** *Let $p \geq 2$. Then there exist constants $\delta_5 > 0$ and $C_{12} > 0$ such that for every connection $A \in \mathcal{A}_\Sigma$ with*

$$\|F_A\|_{L^p(\Sigma, g_\varepsilon)} \leq \delta_5,$$

*there is an estimate*

$$\|(u, \phi)\|_{L^p_1(g_\varepsilon)} = \|u\|_{L^p_1(\Sigma)} + \|\phi\|_{L^p_1(\Sigma)} \leq C_{12}\|D_A(u, \phi)\|_{L^p(\Sigma, g_\varepsilon)},$$

*for all $(u, \phi) \in (\ker D_A)^\perp$ in $\Omega^1 \oplus \Omega^0(\Sigma, adSU(2))$ and $D_A(u, \phi) = (*d_A u + d_A^* \phi, d_A u)$.*

Proof: For every flat connection $A$, $F_A = 0$ and $\ker D_A = H_A^1 \oplus H_A^0(\Sigma, adSU(2))$, the estimate holds. Given a flat connection $\alpha$, there exist constants $\delta_5 > 0$ and $C_{12} > 0$ such that the estimate holds for a $\mathcal{G}_\Sigma$-invariant neighborhood $U_\alpha$ in $\mathcal{A}_\Sigma$ and $\Gamma_\alpha$- invariant neighborhood $V_{\alpha, \delta_5}$ of $0$ in $H^1(\Sigma, ad\alpha)$, where

$$V_{\alpha, \delta_5} \oplus U_{\alpha, \delta_5} = \{(u, \phi)|A = \alpha + u, \|u\|_{L^p_1} + \|\phi\|_{L^p_1} < \delta_5, \|F_A\|_{L^p} < \delta_5\}. \qquad (6.22)$$

Suppose the contrary. Then there would have a sequence $A_n \in \mathcal{A}_\Sigma$ and $(u_n, \phi_n) \in (\ker D_{A_n})^\perp$ such that for $\delta_5(n) \to 0$,

$$\frac{1}{n}\|(u_n, \phi_n)\|_{L^p_1} > \|D_{A_n}(u_n, \phi_n)\|_{L^p}. \qquad (6.23)$$

By Uhlenbeck's compactness theorem, there exist a subsequence (still call) $A_n \in \mathcal{A}_\Sigma$ and $g_n \in \mathcal{G}_\Sigma$ such that $g_n^*(A_n)$ converges to a flat connection $\alpha$ in $L^p$-norm. So for $n$ large enough, $\|A_n - \alpha\|_{L^p_1} <$ $\delta_5$, then there is a $L^2$-projection map which is an injective from $\ker D_\alpha$ to $\ker D_{A_n}$ (as in Lemma 6.2.3, see also §7.1.2 in [10] for construction). Since both are closed, finite dimensional with the same rank for $\ker D_\alpha$ and $\ker D_{A_n}$, it gives a way to identify all the spaces $(\ker D_{A_n})^\perp$ for $n$ large with $(\ker D_\alpha)^\perp$. Now we normalize the subsequence $(u_n, \phi_n)$ in $(\ker D_\alpha)^\perp$ (after the identification) so that $\|(u_n, \phi_n)\|_{L_1^p} = 1$,

$$\frac{1}{n} = \frac{1}{n}\|(u_n, \phi_n)\|_{L_1^p} > \|D_{A_n}(u_n, \phi_n)\|_{L^p}.$$

and $(u_n, \phi_n)$ has a weak limit $(u_\alpha, \phi_\alpha) \neq (0, 0)$ in $L_1^p$, it follows

$$D_\alpha(u_\alpha, \phi_\alpha) = 0, \quad (u_\alpha, \phi_\alpha) \in (\ker D_\alpha)^\perp.$$

Then the contradiction proves the lemma. ∎

**Corollary 6.2.21** *There exists $0 < \varepsilon_4' \leq \varepsilon_4$ such that for $0 < \varepsilon \leq \varepsilon_4'$, we have*

$$\|F_A\|_{L^p(\Sigma, g_Y)} \leq 2\delta_5,$$

$$\|(u, 0)\|_{L^q(\Sigma, g_Y)} \leq 8C_{12}\|D_A(u, 0)\|_{L^p(\Sigma, g_Y)},$$

*for all $(u, 0) \in (\ker D_A)^\perp$ and $1/4 + 1/q \leq 1/p$.*

Proof: Over the annulus region, the metric has been changed to $C^0$-close metric $g_Y$ of $g_\varepsilon$. Thus the first inequality follows from (6.20). But the estimate for $D_A$ is no longer true for $(u, \phi)$ since $D_A$ contains $d_A^*$ which requires the derivate estimate of the metric. Fortunately, the deformation for flat connection is in the space of $\Omega^1(\Sigma, adSU(2))$, i.e. we can take $\phi = 0$. Now

$$\begin{aligned}
\|D_A(u, 0)\|_{L^p(\Sigma, g_\varepsilon)} &= \|(*d_A u, d_A u)\|_{L^p(\Sigma, g_\varepsilon)} \\
&\leq \|(*_{g_Y} d_A u, d_A u)\|_{L^p(\Sigma, g_\varepsilon)} + \|(*_{g_Y} - *_{g_\varepsilon})d_A u\|_{L^p(\Sigma, g_\varepsilon)} \\
&\leq 2\|*_{g_Y} d_A u\|_{L^p(\Sigma, g_Y)} + 2(1 + C\varepsilon)\|d_A u\|_{L^p(\Sigma, g_Y)}.
\end{aligned}$$

Choosing $\varepsilon_4' \leq \varepsilon_4$ with $2C\varepsilon_4' < 1/2$, we have

$$\|(u, 0)\|_{L^q(\Sigma, g_Y)} \leq 8C_{12}\|D_A(u, 0)\|_{L^p(\Sigma, g_Y)},$$

from the Lemma 6.2.20 and changing metrics (6.20). ∎

In order to apply Lemma 6.1.2 to $\overline{A}$ on the $Y_0 = \Sigma \times [-1+S^{-1}\varepsilon, 1-S^{-1}\varepsilon]$, we use the shrinking metric $g_Y$. For an ASD connection (balanced one), we have a pointwise estimate of the curvature term in Proposition 4.1.1 (2), but $\overline{A}$ is further perturbed on $\Sigma \times ([-1+S^{-1}\varepsilon, -1+S\varepsilon] \cup [1-S\varepsilon, 1-S^{-1}\varepsilon])$ from an perturbed ASD connection $A$.

**Lemma 6.2.22** *For the connection $\overline{A}$ and $2 \leq p \leq 3$, there exists a constant $C_{13}$ such that*

$$\|F_{\overline{A}}\|_{L^p(Y \times \mathbf{R}, g_Y)} \leq C_{13}.$$

*Furthermore, there are a $\gamma_1 > 0$ and $\gamma_1 > \gamma/2$ and constant $C_{14}$, so that*

$$sup|F_{\overline{A}}|_{y,t} \leq C_{14} e^{-\gamma_1 |t|}.$$

Proof: The curvature term can be calculated as follows:

$$F_{\overline{A}} = F_A + d_A(\eta_1 a_1 + \eta_2 a_2) + (\eta_1 a_1 + \eta_2 a_2) \wedge (\eta_1 a_1 + \eta_2 a_2). \tag{6.24}$$

We estimate the third term first.

$$\begin{aligned}
\|(\eta_1 a_1 + \eta_2 a_2) \wedge (\eta_1 a_1 + \eta_2 a_2)\|_{L^p(Y \times \mathbf{R}, g_Y)} &= \|(\eta_1 a_1 + \eta_2 a_2)\|^{1/2}_{L^{2p}(Y \times \mathbf{R}, g_Y)} \\
&\leq (\|a_1\|_{L^{2p}(Y_1' \times \mathbf{R}, g_Y)} + \|a_2\|_{L^{2p}(Y_2' \times \mathbf{R}, g_Y)})^{1/2} \\
&\leq (2C_{15}\|a_1\|_{L^2_k(g_{1\delta})} + 2C_{15}\|a_2\|_{L^2_k(g_{2\delta})})^{1/2} \\
&\leq C_{16} \varepsilon_4^{1/2}.
\end{aligned}$$

The first inequality is from the support of $\eta_i$ and

$$Y_1' = Y_1 \cup (\Sigma \times [-1+S^{-1}\varepsilon, -1+S\varepsilon]), \quad Y_2' = Y_2 \cup (\Sigma \times [1-S\varepsilon, 1-S^{-1}\varepsilon]). \tag{6.25}$$

The second follows from changing the metric $g_Y$ into $g_{j\delta}$ and the Sobolev imbedding theorem, and the last from Theorem 6.2.17.

The second term in the curvature of $\overline{A}$ is

$$d_A(\eta_1 a_1 + \eta_2 a_2) = [\eta_1'(s) a_1 + \eta_2'(s) a_2] + [\eta_1 d_A a_1 + \eta_2 d_A a_2].$$

¿From the support of $\eta_i$ and the compactness of balanced 1-dimensional moduli space, we have

$$\|[\eta_1 d_A a_1 + \eta_2 d_A a_2]\|_{L^p(Y \times \mathbf{R}, g_Y)} \leq C_{17} \varepsilon_4,$$

from Theorem 6.2.17 and Lemma 6.2.14. The estimate of $\sum_{j=1}^{2} \eta_j'(s)a_j(t)$ can be done by using the exponential decay property of $a_i(t)$ from Proposition 6.2.12 and Theorem 6.2.16 :

$$
\begin{aligned}
\|\eta_1'(s)a_1(t)\|_{L^p(g_Y)} &= (\int_{-1+S^{-1}\varepsilon}^{-1+S\varepsilon} |\eta_1'(s)|^p ds \cdot \int_{\Sigma \times \mathbf{R}, g_Y} |a_1(t)|^p)^{1/p} \\
&\leq C_{18}\varepsilon^{-1+1/p}(\int_{\Sigma \times \mathbf{R}, g_Y} |a_1(t)|^p)^{1/p} \\
&\leq C_{18}'\varepsilon^{-1+1/p}Vol(\Sigma, g_Y)^{1/p} \\
&\leq C_{18}''\varepsilon^{-1+3/p}
\end{aligned}
$$

Combining the two terms, we obtain

$$\|(\chi_1 a_1 + \chi_2 a_2)\|_{L_1^p(Y \times \mathbf{R}, g_Y)} \leq C_{19}\varepsilon^{-1+3/p}, \tag{6.26}$$

for $0 < \varepsilon < \varepsilon_4$.

The first term $F_A$ is in $L^p$ for $p \geq 2$. Hence we have the $L^p$ bound for the curvature of $\overline{A}$. Note that $\gamma = \gamma(\lambda_0) > 0$ depends continuously on the smallest value of absolute eigenvalues of $D_a$ (see [9]). Over compact manifold $Y$, the eigenvalues of $D_a$ and its compact perturbation $D_{a'}$ are sufficiently close to each other. Hence $\gamma_1 = \gamma(\lambda_0 + \varepsilon')$ so that $\gamma_1 > \gamma/2$ for $\varepsilon'$ sufficiently small. ∎

**Remark:** The Lemma 6.2.22 recaptured the $\delta$-decay property for the deformed ASD $\overline{A}$. Hence the deformation on handlebodies still live in the correct space.

**Lemma 6.2.23** *For any $\delta_5 > 0$, there exists $\varepsilon_5 > 0$ such that for all $0 < \varepsilon < \varepsilon_5$ we have that*

$$\|F_{\overline{A}}\|_{L^p(\Sigma, g_Y|_\Sigma)} < \delta_5.$$

Proof: The shrinking metric $g_Y$ on the surface $\Sigma$ gives the volume close to $Vol(\Sigma, g^\varepsilon) = O(\varepsilon^2)$. Hence by Lemma 6.2.22

$$
\begin{aligned}
\|F_{\overline{A}}\|_{L^p(\Sigma, g_Y)} &\leq C_{14}e^{-\gamma_1|t|}(Vol(\Sigma, g^\varepsilon))^{1/p} \\
&\leq C_{20}e^{-\gamma_1|t|}\varepsilon^{2/p}
\end{aligned}
$$

Choosing $\varepsilon_5$ so small that $C_{20}e^{-\gamma_1|t|}\varepsilon^{2/p} \leq \delta_5/2$, we obtain the desired estimate. ∎

**Proposition 6.2.24** *There exists $\varepsilon_5 > 0$ such that for all $0 < \varepsilon < \varepsilon_5$, $\overline{A}(s,t)$ can be deformed into $\mathcal{R}_\Sigma$ for all $(s,t) \in [-1,1] \times \mathbf{R}$. The element $\overline{A}(s,t) + b(s,t)$ has the property $\|b(s,t)\|_{L_1^p(\Sigma)} < \varepsilon_5$ and $b(s,t) \in \Omega^1(\Sigma, adSU(2))$.*

Proof: This is solving the equation on $\Omega^1(\Sigma, adSU(2)) \oplus \{0\}$

$$D_{\overline{A}}(b, \phi) + (b \wedge b, 0) + (F_{\overline{A}}, 0) = (0, 0). \tag{6.27}$$

Hence the result follows from inverse function therorem with the aid of Lemma 6.2.23 and Corollary 6.2.21. ∎

### 6.2.3 Deforming the perturbed ASD into a holomorphic curve

Now we have deformed an ASD $A$ over $Y \times \mathbf{R}$ to $A_d = \overline{A} + b(s,t)$ which satisfies

$$A_d|_{Y_j} \in \mathcal{R}_j, \quad j = 1, 2; \quad A_d|_{\Sigma \times \{(s,t)\}} \in \mathcal{R}_0, \quad (s,t) \in [-1 + S^{-1}\varepsilon, 1 - S^{-1}\varepsilon] \times \mathbf{R}.$$

Furthermore $A_d \in \mathcal{A}_{1,\delta}^p$ for $3 \geq p > 2$ (recall the definition in §2.2). Previousely we have that the ASD $A$ is automatically a solution of $\overline{\partial}_J u = 0$, but this is no longer true for $A_d$. We will first estimate that the deformed ASD $A_d$ is not far away from being an anti-self-dual connection. Thus the smallest eigenvalue estimate also holds for the self-duality operator twisted by $A_d$. Using the comparsion of self-duality operator and Cauchy-Riemann operator, we get the smallest eigenvalue estimate for the Cauchy-Riemann operator twisted by $A_d$. We also show that the connection $A_d$ is an almost holomorphic curve since $A_d$ is an almost anti-self-dual connection. Then inverse function theorem provides us the solution for the deformed holomorphic curve equation.

**Proposition 6.2.25** *For any $\delta_6 > 0$, there exists a $\varepsilon_6 > 0$ such that for $0 < \varepsilon < \varepsilon_6$, $2 < p < 3$, we obtain the following inequlities with respect to $g_Y$-metric:*

$$\|A - A_d\|_{L_{1,\delta}^p} < \delta_6; \quad \|F_{A_d}^+\|_{L_{0,\delta}^p} < \delta_6.$$

Proof: Note that $A_d - A = \chi_1 a_1 + \chi_2 a_2 + b$. In Lemma 6.2.22 (6.27), we have for $2 < p < 3$ and $\varepsilon < \varepsilon_3$,

$$\|\chi_1 a_1 + \chi_2 a_2\|_{L_{1,\delta}^p} \leq C_{19} \varepsilon^{-1 + 3/p}.$$

By Proposition 6.2.12, we have
$$\|b\|_{L^p_{1,\delta}(Y\times[T,\infty))} \leq \varepsilon_j.$$

By Proposition 6.2.24, we get
$$\|b\|_{L^p_{1,\delta}(Y\times[-T,T])} \leq C_{21}\varepsilon_5.$$

Thus we have the first inequality $\|A - A_d\|_{L^p_{1,\delta}} < C_{19}\varepsilon^{-1+3/p} + \varepsilon_j$. Now

$$F^+_{A_d} = d^+_A(\chi_1 a_1 + \chi_2 a_2) + (\chi_1 a_1 + \chi_2 a_2) \wedge (\chi_1 a_1 + \chi_2 a_2)_+$$
$$+ d^+_A b + (b \wedge b)_+ + [(\chi_1 a_1 + \chi_2 a_2), b]_+. \tag{6.28}$$

¿From the estimates in Lemma 6.2.22 and equation (6.20), we obtain

$$\|d^+_A(\chi_1 a_1 + \chi_2 a_2) + (\chi_1 a_1 + \chi_2 a_2) \wedge (\chi_1 a_1 + \chi_2 a_2)_+\|_{L^p_{0,\delta}} < C_{22}\varepsilon^{-1+3/p}.$$

$$\begin{aligned}
\|d^+_A b + (b \wedge b)_+\|_{L^p_{0,\delta}} &\leq \|d^+_A b\|_{L^p_{0,\delta}} + \|(b \wedge b)_+\|_{L^p_{0,\delta}} \\
&\leq \|b\|_{L^p_{1,\delta}(A)} + \|b\|^{1/2}_{L^{2p}_{0,\delta/2}} \\
&\leq C_{21}\varepsilon_5 + C'_{23}\|b\|^{1/2}_{L^p_{1,\delta}} \\
&\leq C_{23}\varepsilon_5.
\end{aligned}$$

The last term in $F^+_{A_d}$ can be also estimated by Hölder inequlity and previouse two estimates. Choose $\varepsilon_6$ such that
$$\max\{C_{19}\varepsilon^{-1+3/p} + \varepsilon_j, C_{22}\varepsilon^{-1+3/p} + C_{23}\varepsilon_5\} \leq \delta_6,$$

we get the curvature estimate. ∎

In the above, we have showed that $A_d$ is in the neighborhood of anti-self-dual moduli space

$$U_{\delta_6,\mathcal{M}_{Y\times\mathbf{R}}} = \{B \in \mathcal{B}_{Y\times\mathbf{R}} | \text{exists } A \in \mathcal{M}_{Y\times\mathbf{R}} \text{ such that } \|A - B\|_{L^p_{1,\delta}} < \delta_6, \|F^+_B\|_{L^p_{0,\delta}} < \delta_6\}.$$

**Lemma 6.2.26** *There exists $\delta_7$ such that for $0 < \delta_6 < \delta_7$ there is a constant $C_{24}$ independent of $\delta_6$ such that for all $B \in U_{\delta_6,\mathcal{M}_{Y\times\mathbf{R}}}$*

$$\|u\|_{L^p_{1,\delta}(Y\times\mathbf{R})} \leq C_{24}\|(d^+_B)^*\|_{L^p_{0,\delta}(Y\times\mathbf{R})}. \tag{6.29}$$

Proof: This is Lemma 3.2.11 and Lemma 3.2.12 in [25]. ∎

Now the restriction of self-duality operator from $\Omega^1 \oplus \Omega^0 \oplus \Omega^0(\Sigma, adSU(2))$ to $\Omega^1(\Sigma, adSU(2)) \oplus \{0\} \oplus \{0\}$ will give the smallest eigenvalue for the Cauchy-Riemann operator.

**Corollary 6.2.27** *There exists $\delta_7$ such that for $0 < \delta_6 < \delta_7$ there is a constant $C_{24}$ independent of $\delta_6$ such that*

$$\|\xi\|_{L^p_{1,\delta}(\Sigma \times \mathbf{R})} \leq C_{24} \|(E_{A_d})^* \xi\|_{L^p_{0,\delta}(\Sigma \times \mathbf{R})}.$$

The following Lemma shows that the deformed ASD connection $A_d$ is almost holomorphic curve.

**Lemma 6.2.28** *There exists $\varepsilon_6 > 0$ such that for $3 > p > 2$, $0 < \varepsilon < \varepsilon_6$ and $\delta_6$ in Corollary 6.2.27,*

$$\|\overline{\partial}(A_d)\|_{L^p_{0,\delta}} \leq \delta_6,$$

*where $\overline{\partial}(A_d) = \frac{\partial A_d}{\partial t} - d_{A_d}\Psi + *_s(\frac{\partial A_d}{\partial s} - d_{A_d}\Phi)$.*

Proof: From the ASD connection $A = A_0 + \Phi ds + \Psi dt$, we have that

$$A_d = A_0 + \Phi ds + \Psi dt + (\chi_1 a_1 + \chi_2 a_2) + b.$$

The holomorphic curve term $\frac{\partial A_d}{\partial t} - d_{A_d}\Psi + *_s(\frac{\partial A_d}{\partial s} - d_{A_d}\Phi)$ is one of the components in $F^+_{A_d}$. Thus the result follows from Proposition 6.2.25. ∎

Floer described explicite charts of $\mathcal{P}^p_{1,\delta}(a,b)$ and trivializations of $T_u\mathcal{P}^p_{1,\delta}(a,b)$ in Theorem 3 of [14], where

$$E_u : T_u\mathcal{P}^p_{1,\delta}(a,b) \to L^p_{0,\delta}(u^* T\mathcal{R}_\Sigma),$$

and the deformed holomorphic curve equation is

$$f(\xi) = \overline{\partial}(A_d) + E_{A_d}\xi + N_{A_d}(\xi).$$

The nonlinear term $N_{A_d}(\xi)$ satisfies the estimate

$$\|N_{A_d}(\xi) - N_{A_d}(\eta)\|_{L^p_{0,\delta}} \leq C(A_d, \delta)\|\xi - \eta\|_{L^p_{1,\delta}}(\|\xi\|_{L^p_{1,\delta}} + \|\eta\|_{L^p_{1,\delta}}), \tag{6.30}$$

where the constant $C(A_d, \delta)$ depending only on $\|\nabla A_d\|_{L^\infty}$ and $\delta$. The estimate 6.2.3 has been obtained in [15] (4.6) for $\delta = 0$. The weighted Sobolev norm estimate follows from weighted Hölder inequality and Lemma 7.2 in [26] (see also Theorem 3.3.6 in [25] for the instanton case). Now we apply the inverse function theorem (Lemma 6.1.2) to

$$f(\xi) = \overline{\partial}(A_d) + E_{A_d}\xi + N_{A_d}(\xi),$$

with $f(0) = \overline{\partial}(A_d), Df(0) = E_{A_d}$ with the bounded right inverse from Corollary 6.2.24, and $E = L^p_{1,\delta} \cap T_{A_d}\mathcal{P}^p_{1,\delta}(a,b), F = L^p_{0,\delta}(A_d^*T\mathcal{R}_\Sigma)$.

**Theorem 6.2.29** Let $A \in \mathcal{M}(a,b)$ with $dim\mathcal{M}(a,b) = 1$ and $\delta_7$ in Corollary 6.2.27. Then if $0 < \varepsilon < \varepsilon_6$ and $2 < p < 3$, we can deform $A$ to a smooth holomorphic curve $A_d + \xi \in \mathcal{P}^p_{1,\delta}(a,b)$, where $\|\xi\|_{L^p_{1,\delta}}$ is sufficiently small less than $\delta_7$.

Proof: Using Lemma 6.2.28 and Corollary 6.2.27, we have

$$\|G_{A_d}\overline{\partial}(A_d)\|_{L^p_{1,\delta}} \leq C_{25}\|F^+_{A_d}\|_{L^p_{0,\delta}} \leq C_{25}\delta_6,$$

where $G_{A_d}$ is the right inverse operator of $E_{A_d}$. Also

$$\|G_{A_d}N_{A_d}(\xi) - G_{A_d}N_{A_d}(\eta)\|_{L^p_{1,\delta}} \leq C_{25}C(A_d,\delta)\|\xi - \eta\|_{L^p_{1,\delta}}(\|\xi\|_{L^p_{1,\delta}} + \|\eta\|_{L^p_{1,\delta}}).$$

We apply Lemma 6.1.2 with $\delta_1 = (8C_{25}C(A_d,\delta))^{-1}$ for the metric $(\Sigma, g_\varepsilon)$. So $A_d + \xi$ is an holomorphic curve over $\mathcal{R}_\Sigma$ with $\|\xi\|_{L^p_{1,\delta}}$ small and is smooth by Proposition 2.2.17 for elements in $dim\mathcal{M}^p_{1,\delta}(a,b) = 1$. ∎

Theorem 6.2.29 provides an injective map

$$\mathcal{T}_\varepsilon : \hat{\mathcal{M}}_{g_Y}(a,b) \to \hat{\mathcal{M}}_J(a,b),$$

for $0 < \varepsilon < \varepsilon_6$ and $dim\mathcal{M}(a,b) = 1, dim\mathcal{M}_J(a,b) = 1$. Now we take $\varepsilon_7 = \min\{\varepsilon_1, \varepsilon_6\}$ for $\varepsilon_1$ in Theorem 6.1.3. Thus for $0 < \varepsilon < \varepsilon_7$ we obtain

$$T_\varepsilon : \hat{\mathcal{M}}_J(a,b) \to \hat{\mathcal{M}}_{g_\varepsilon}(a,b), \tag{6.31}$$

$$\mathcal{T}_\varepsilon : \hat{\mathcal{M}}_{g_Y}(a,b) \to \hat{\mathcal{M}}_J(a,b). \tag{6.32}$$

with the same metric on the tube $\Sigma \times [-1 + S^{-1}\varepsilon, 1 + S^{-1}\varepsilon]$.

Note that the map $T_\varepsilon$ has no requirement on the metrics of two handlebodies, so we can just take the metrics needed in the process of deforming $T_\varepsilon$. But the metric on $\Sigma \times [-1 + S^{-1}\varepsilon, 1 + S^{-1}\varepsilon]$ has been changed slightly. I.e $A$ is an ASD with repsect to $g_\varepsilon$, then $A$ is an almost ASD with respect to $g_Y$ since $\|g_Y - g_\varepsilon\|_{C^0} < \varepsilon$. So $d_A^{+g_Y}$ has a right inverse from the right inverse of $d_A^{+g_\varepsilon}$ by Lemma 3.2.12 in [25]. Thus we can deform $g_\varepsilon$-ASD $A$ into a $g_Y$-ASD $A + a$, and vise verse. So we have an orientation preserving bijective map

$$U : \hat{\mathcal{M}}_{g_\varepsilon}(a,b) \to \hat{\mathcal{M}}_{g_Y}(a,b),$$

for $\mu(a) - \mu(b) = 1$. Now we have two injective maps

$$U \circ T_\varepsilon : \hat{\mathcal{M}}_J(a,b) \to \hat{\mathcal{M}}_{g_Y}(a,b), \quad T_\varepsilon : \hat{\mathcal{M}}_{g_Y}(a,b) \to \hat{\mathcal{M}}_J(a,b).$$

The inverse function theorem provides the uniqueness for each map, i.e. $U \circ T_\varepsilon$ and $T_\varepsilon$ are injective. Since $\hat{\mathcal{M}}(a,b)$ and $\hat{\mathcal{M}}_J(a,b)$ are compact for the balanced moduli spaces, it follows that $U \circ T_\varepsilon$ and $T_\varepsilon$ are bijective maps between the two balanced moduli spaces.

For $u \in \hat{\mathcal{M}}_J(a,b)$, we have the following inequalities from the inverse function theorem

$$\|U \circ T_\varepsilon u - u\|_{L^p_{1,\delta}} \leq \delta'_7,$$

$$\|T_\varepsilon(U \circ T_\varepsilon u) - U \circ T_\varepsilon u\|_{L^p_{1,\delta}} \leq \delta_7.$$

So this shows $T_\varepsilon(U \circ T_\varepsilon)$ is very close to the indentity map. Hence our method to avoid the ontoness of $T_\varepsilon$ gives the same result as showing that every ASD connection is in the image of the deformed map $T_\varepsilon$ up to gauge transformation in [11].

Let $O_J(a,b), O(a,b)$ be the spaces of orientations of the determined line bundles for Cauchy-Riemann operators and self-duality operators, respectively. The following proposition is Proposition 10.2 in [11].

**Proposition 6.2.30** *For every pair $a, b \in \mathcal{R}_Y^*$, there is a natural bijection*

$$\tau : O_J(a,b) \to O(a,b).$$

Now we have our main theorem.

**Theorem 6.2.31** *For a Heegaard decomposition $(Y; Y_1, Y_2; Y_0)$ with genus $g \geq 3$. We have the following natural isomorphism*

$$HF_*(Y; Z) \cong HF_*^{sym}(\mathcal{R}_1, \mathcal{R}_2; \mathcal{R}_0).$$

Proof: Proposition 3.2.2 gives a natural isomorphism (3.15) between the two chain complexes. Choose coherent orientation $\sigma_J \in O_J(a,b)$ and consider the induced coherent orientation $\sigma \in O(a,b)$. These two $\sigma_J, \sigma$ determine the orientations of moduli spaces $\mathcal{M}(a,b)$ and $\mathcal{M}_J(a,b)$, also the balanced moduli spaces. These orientations are invariant under Floer's glueing maps in the symplectic case [15] and in the instanton case [13] (see the proof of Proposition 10.3 in [11]). Now fix $a, b \in \mathcal{R}^*(Y)$ and consider the map

$$U \circ T_\varepsilon : \hat{\mathcal{M}}_J(a,b) \to \hat{\mathcal{M}}_{g_Y}(a,b).$$

The induced map on the space of orientations agrees with the map $\tau$ in Proposition 6.2.27. Hence $U \circ T_\varepsilon$ is orientation preserving.

For $\dim \mathcal{M}_J(a,b) = 1$, $\dim \mathcal{M}(a,b) = 1$, $U \circ T_\varepsilon$ is bijective of finite sets from the balanced moduli space $\hat{\mathcal{M}}_J(a,b)$ to the balanced moduli space $\hat{\mathcal{M}}_{g_Y}(a,b)$ for $0 < \varepsilon < \varepsilon_7$ sufficiently small, and also preserves the signs for each element in $\hat{\mathcal{M}}_J(a,b)$ and $\hat{\mathcal{M}}_{g_Y}(a,b)$. I.e. the same argument in the proof of Theorem 10.1 in [11] works. Hence the Floer boundary maps $\partial_J^{sym}$ and $\partial$ agree for $0 < \varepsilon < \varepsilon_7$ sufficiently small. This identifies the two Floer homologies. ∎

**Corollary 6.2.32** *Let $Y$ be an integral homology 3-sphere and $(Y; Y_1, Y_2; Y_0)$ be its Heegaard decomposition. The symplectic Floer homology $HF_*^{sym}(\mathcal{R}(Y_1), \mathcal{R}(Y_2); \mathcal{R}(Y_0))$ is independent of the Heegaard decompositions. The Casson invariant of the intergral homology 3-sphere is independent of the Heegaard decomposition (see [1]).*

Yale University, New Haven, CT 06520, U.S.A.

Ronnie Lee, rlee@math.yale.edu

Weiping Li, wli@math.yale.edu; wli@math.okstate.edu